\documentclass[preprint,12pt]{elsarticle}
\usepackage{graphicx} 
\graphicspath{{./figs/}}
\usepackage{float}
\usepackage{morefloats}
\usepackage{mathtools}
\usepackage{amssymb,amsfonts,amsmath}
\newcommand{\mbf}[1]{\mathbf{x}} \newcommand{\eps}{\varepsilon}

\newcommand{\omr}{\mathbf{f}}
\newcommand{\prm}{\mathbf{p}}
\newcommand{\ps}{\mathbf{p}^*}
\newcommand{\peff}{p_\mathrm{eff}}
\newcommand{\keff}{k_\mathrm{eff}}
\newcommand{\Ord}{\mathcal{O}}

\begin{document}
\begin{frontmatter}
\title{Manifold learning for parameter reduction}
\author[add1]{Alexander Holiday\corref{cor2}}
\author[add2]{Mahdi Kooshkbaghi\corref{cor2}}
\author[add2]{Juan M. Bello-Rivas}
\author[add1]{C. William Gear}
\author[add3]{Antonios Zagaris\corref{cor1}}
\ead{antonios.zagaris@asml.com}
\author[add1,add2,add4]{Ioannis G. Kevrekidis\corref{cor1}}
\ead{yannisk@jhu.edu}
\address[add1]{Chemical and Biological Engineering, Princeton University, USA}
\address[add2]{The Program in Applied and Computational Mathematics (PACM), 
Princeton University, USA}
\address[add3]{Wageningen Bioveterinary Research, Wageningen UR, The Netherlands}
\address[add4]{Department of Chemical and Biomolecular Engineering, John Hopkins University, USA}
\cortext[cor2]{These two authors contributed equally to this work.}
\cortext[cor1]{Corresponding author}

\begin{keyword}
model reduction, data mining, diffusion maps,
data driven perturbation theory, parameter sloppiness
\end{keyword}
%
\begin{abstract}
Large scale dynamical systems (e.g. many nonlinear coupled differential equations) 
can often be summarized in terms of only a few state variables (a few equations), a
trait that reduces complexity and facilitates exploration of behavioral 
aspects of otherwise intractable models. High model dimensionality and complexity
makes symbolic, pen--and--paper model 
reduction tedious and impractical, a difficulty addressed by recently developed frameworks that computerize reduction. Symbolic work has the benefit, however, of identifying both reduced state variables and parameter combinations that matter most (\emph{effective parameters}, ``inputs"); whereas current computational reduction schemes leave the parameter reduction aspect mostly unaddressed. 
As the interest in mapping out and 
optimizing complex input--output relations keeps growing, it becomes clear that
combating the curse of dimensionality also
requires efficient schemes for input space exploration and reduction. 
Here, we explore systematic, data-driven parameter reduction by means of {\em effective parameter identification}, starting from current nonlinear manifold-learning techniques enabling state space reduction. Our approach aspires to extend the data-driven determination of effective state variables with the data-driven discovery of {\em effective model parameters}, and thus to accelerate the exploration of high-dimensional parameter spaces associated with complex models.
\end{abstract}

\end{frontmatter}

\section{Introduction}
Our motivation lies in the work of Sethna and coworkers on model sloppiness \cite{gutenkunst2007universally}, as well as in related ideas and studies on parameter non-identifiability \cite{raue_structural_2009}, active subspaces \cite{constantine_active_2014} and more. These authors investigate a widespread phenomenon, in which large ranges of model parameter values (inputs) produce nearly constant model predictions (outputs). This behavior, termed \emph{sloppiness} and observed in complex dynamic models over a wide range of fields, has been exploited to derive simplified models \cite{transtrum2010nonlinear,transtrum2014model}. Additional motivation comes from our interest in model {\em scaling} and \emph{nondimensionalization}, time-honored ways to reduce complexity but often more closely resembling an art than definite algorithms.

One extreme case of sloppiness, termed {\em parameter non-identifiability}, arises when model predictions depend solely on a reduced number of parameter combinations. In such a setting, the parameter space is foliated by lower-dimensional sets along which those combinations, and hence also the resulting {\em observables} (the outputs), retain their values. In such circumstances, it is neither possible nor desirable to infer parameter values from observations; the parameters are said to be \emph{non-identifiable}. One should, instead, re-parameterize the model with a reduced number of \emph{identifiable, effective parameters} and, if desired, use those to explore the model input--output structure. Such identifiability analysis decomposes parameter space \emph{globally} on the basis of model response, yet its symbolic
nature can make it cumbersome and highly sensitive to small perturbations: even a minute dependence on certain parameter combinations can destroy the invariance of the decomposition. (Computational) sensitivity analysis is more robust, as it weighs the degree by which parameter combinations affect response; however, it is inherently not global in parameter space, as it uses a (local) linearization. We attempt to reconcile and fuse these two perspectives into an entirely data-driven, {\em nonlinear} framework for the identification of {\em global effective parameters}.

\begin{figure}[!h]
\centering
\begin{tabular}{cc}
\includegraphics[width=0.47\textwidth]{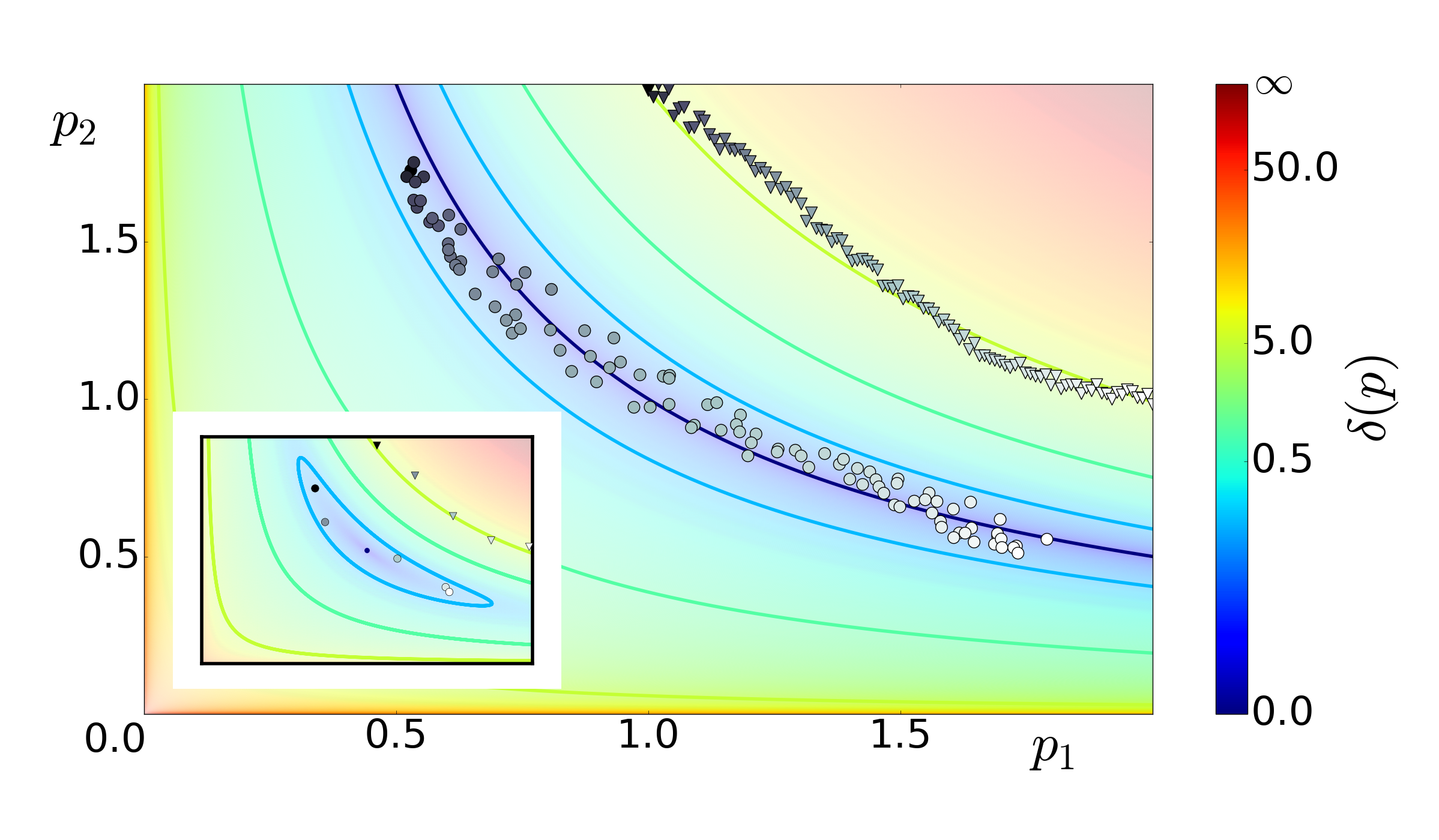}&
\includegraphics[width=0.51\textwidth]{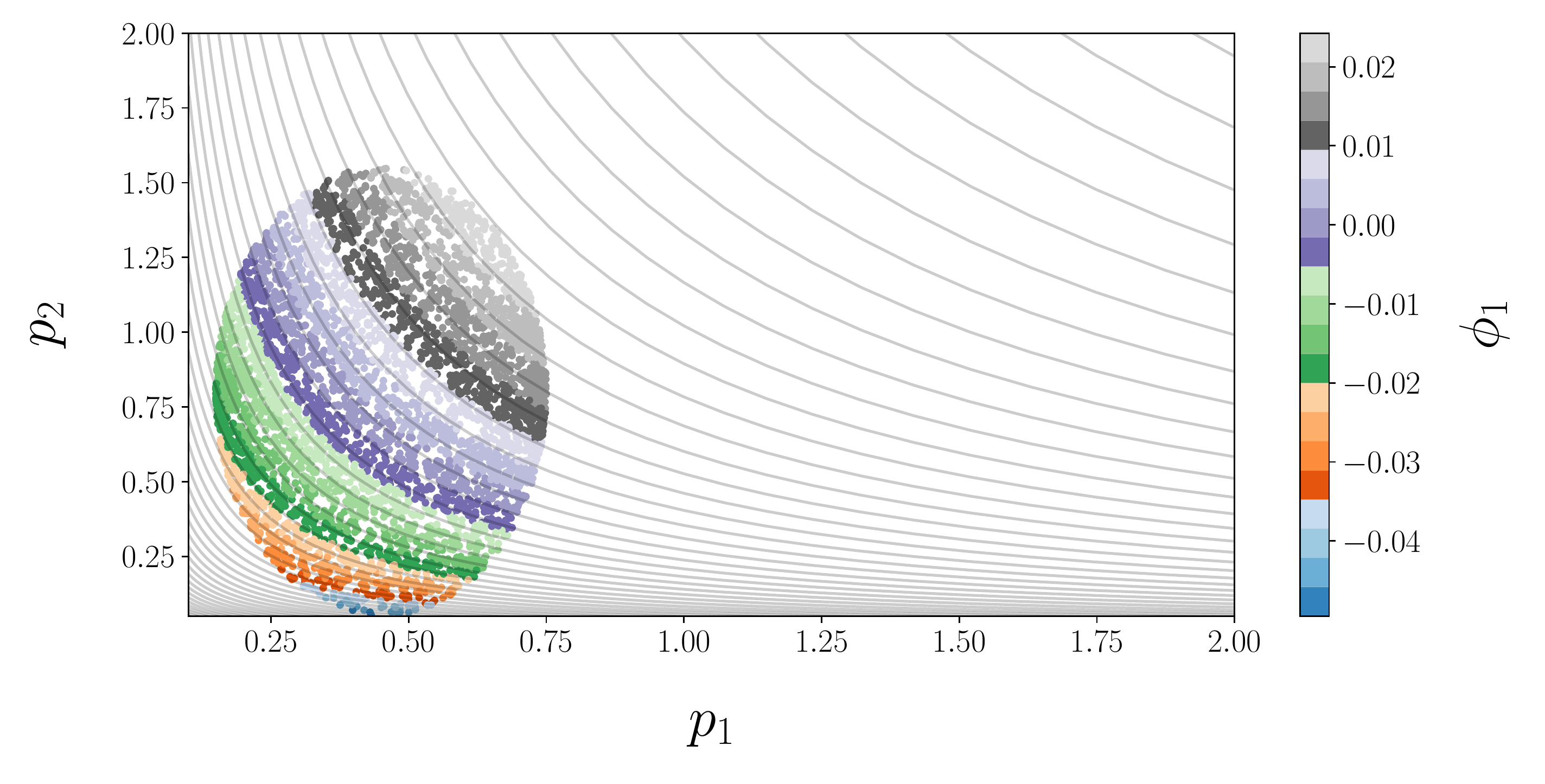}\\
(a) & (b)
\end{tabular}
\caption{
Exact and learned (global) parameter space foliations for the model $\omr_\eps(p_1,p_2) = (p_1 p_2 + 2\eps(p_1-p_2) , \ln(p_1 p_2) , (p_1 p_2)^2)$. The combination $\peff = p_1 p_2$ is an effective parameter for the unperturbed ($\eps=0$) model, since $\omr_0 =$const. whenever $p_1 p_2 =$const. (a) Level sets of the cost function $\delta(\prm) = \|\omr_\eps(\prm) - \omr^*\|$ for the unperturbed (main) and perturbed (inset) model and for data $\omr^*=(1,0,1)$ corresponding to $(p_1,p_2)=(1,1)$. Level sets of $\peff$ can be learned by data fitting: feeding various initializations (triangles) to a gradient descent algorithm for the unperturbed problem yields, approximately, the hyperbola $\peff = 1$ (circles; colored by initialization). This behavior persists qualitatively for $\eps=0.2$ despite the existence of a unique minimizer $\ps$, because $\delta(\prm)$ remains within tolerance over extended \emph{almost neutral sets} around $\ps$ that approximately trace the level sets of $\peff$. (b) Learning $\peff$ by applying DMAPS, with an output-only-informed metric (see the SI), to input--output data of the unperturbed model. For $\eps=0$, points on any level curve of $\peff$ are indistinguishable for this metric, as $\omr_0$ maps them to the same output. DMAPS, applied to the depicted oval point cloud, recovers those level curves as level sets of the single leading nontrivial DMAPS eigenvector $\phi_1$.
}
\label{fig:non-id} 
\end{figure}

To fix ideas, we consider the caricature model of Fig.~\ref{fig:non-id}, given as an explicit vector function of two parameters, $\omr_0(p_1,p_2) = (p_1 p_2 , \ln(p_1 p_2) , (p_1 p_2)^2)$. Given access to input--output information (black-box function evaluation) but no formulas, one might not even suspect that only the single parameter combination $\peff = p_1 p_2$ matters. Fitting the model to data $\omr^* = (1,0,1)$ in the absence of such information, one would find an entire curve in parameter space that fits the observations. A data fitting algorithm based only on function evaluations could be ``confused'' by such behavior in declaring convergence. As seen in Fig.~\ref{fig:non-id}(a), different initial conditions fed to an optimizer with a practical fitting tolerance $\delta\approx10^{-3}$ (see figure caption for details) converge to many, widely different results tracing a level curve of $\peff$. The subset of good fits is effectively $1-$D; more importantly, and moving beyond the fit to this particular data, the entire parameter space is foliated by such $1-$D curves (neutral sets), each composed of points indistinguishable from the model output perspective. Parameter non-identifiability is therefore a {\em structural feature} of the model, not an artifact of optimization. The appropriate, \emph{intrinsic} way to describe parameter space for this problem is through the effective parameter $\peff$ and its level sets.
Consider now the inset of Fig.~\ref{fig:non-id}(a), corresponding to the perturbed model $\omr_\eps(p_1,p_2) = \omr_0(p_1,p_2) + 2 \eps (p_1-p_2 , 0 , 0)$ and fit to the same data. Here, the parameters \emph{are} identifiable and the minimizer $(p_1,p_2)$ unique: a perfect fit exists. However, the foliation observed for $\eps=0$ is loosely remembered in the shape of the residual level curves, and the optimizer would be comparably ``confused'' in practice. It is such model features that provided one of the original motivations in the work of Sethna and coworkers; in their terminology, this model is \emph{sloppy}. The presence of lower-dimensional, almost neutral parameter sets (``echoed'' in the elongated closed curves in the inset) increases disproportionately the importance of certain parameter combinations and reduces accordingly the number of independent, \emph{effective} system parameters.

Our goal is to extract a useful {\em intrinsic} parameterization of model parameter space ({\em input space}) solely from input--output data. As we shall see, this parameterization may vary across input space regimes and, in the context of ODEs, we will associate that variation with the classical notions of regular and singular perturbations using explicit examples. For the time being, a pertinent question concerns the purely data-driven identification of the sloppy structure in Fig.~\ref{fig:non-id}. One answer is given by the manifold-learning technique we choose to work with in this paper: diffusion maps (DMAPS; see SI and e.g.~\cite{coifman_diffusion_2006}). If a given dataset in a high-dimensional, ambient Euclidean space lies on a lower-dimensional manifold, then the DMAPS objective is to parameterize it in a manner reflecting the intrinsic geometry (and thus also dimension) of this underlying manifold. In our case, we work with the space of input--output combinations, where each data point consists of parameter values {\em and} the resulting observations. DMAPS turns the dataset into a weighted graph and models a diffusion process (random walk) on it. The graph weights determine the transition probabilities between points and depend solely on an application-driven understanding of data \emph{closeness} or \emph{similarity}. Typically, DMAPS base this similarity measure on the Euclidean distance in the ambient space; yet, for our applications in most of this paper, this similarity will be informed solely by {\em output observations}. The dataset is parameterized, finally, by eigenvectors of the corresponding Markov matrix, relating in turn to a (discretized) eigenproblem for the Laplace--Beltrami operator on the underlying manifold \cite{lafon2004diffusion}; one may perceive here an analogy with Singular Value Decomposition in classical Principal Component Analysis (PCA) \cite{jolliffe1986principal}. In our input--output setting, DMAPS coordinatizes the low-dimensional manifold hosting the dataset. Both the effective parameters and the observables are now {\em functions} on this low-dimensional manifold, therefore both the input space {\em and} (what in sloppiness terminology is called) the model manifold are {\em jointly} described in terms of this  intrinsic, common parameterization based on leading diffusion modes.

As a concrete example, consider randomly sampling the input space of our model above, i.e. a $(p_1,p_2)-$parallelogram $[0,a]\!\times\![0,b]$ as in Fig.~\ref{fig:non-id}, and using as our pairwise similarity measure the Euclidean distance between points in this input space. Applying DMAPS to that dataset recovers the sampled parallelogram, i.e. DMAPS correctly identifies the dimension of the underlying manifold and coordinatizes it using two diffusion eigenmodes. For this simple shape, the leading (nontrivial, independent) eigenmodes assume the form $\phi_1(x,y) = \cos(\pi p_1/a)$ and $\phi_j(x,y) = \cos(\pi p_2/b)$, where the index $j$ of the first eigenfunction independent of $\phi_1$ depends on $a/b$. This parameterization maps the $(p_1,p_2)-$rectangle bijectively to the $(\phi_1,\phi_2)-$domain $[-1,1]\times[-1,1]$, so DMAPS recovers the original parameterization up to an invertible nonlinear transformation. Our main idea here is to retain sampling of the input (parameter) space but use, instead, a similarity measure {\em (also) informed by the output}, i.e. by the \emph{model response} at the sampled parameter points. As a first but meaningful attempt for the unperturbed example above, we work with the {\em output-only} similarity measure $\|\omr_0(\prm) - \omr_0(\prm')\|$ between parameter settings  $\prm$ and $\prm'$. 
\footnote{We warn the reader that output-only similarity measures may be inappropriate for general input--output relations (e.g. \cite{CoifmanHirn}), such as bifurcation diagrams, in which several behaviors may coexist for a single input. We illustrate this further below, using a system with inputs/outputs that do not maintain a one--to--one correspondence.}

In the context of our example, the output-only similarity measure ensures that only parameter values lying on {\em distinct level sets} of $\peff$ are seen as distinct. Because of this, our chosen similarity measure immediately reveals the {\em effective} parameter space to be $1-$D, as in Fig.~\ref{fig:non-id}(b). Coloring the points by the first DMAPS mode $\phi_1(p_1,p_2)$ confirms that this data-driven procedure ``discovers'' sloppiness. Data points having different parameter settings (different ``genotypes'') but the same output (same ``phenotypes'') are found as level sets of the first nontrivial DMAPS eigenfunction $\phi_1$ on the dataset, obtained in turn by our black-box simulator and without recourse to the explicit input--output relation $\omr_0$. Additionally, the decomposition of parameter space into ``meaningful'' and ``neutral'' parameter combinations can be performed using a small local sample, possibly resulting from a short local search -- e.g. a few gradient descent steps, or local brief simulated annealing runs. This type of {\em local} decomposition can prove valuable to the optimization algorithm, as it reveals local directions that are fruitful to explore and others (along neutral sets) that preserve model predictions (goodness of fit). These latter ones may, in turn, become useful later in {\em multi-objective optimization}, where one optimizes {\em additional} objectives along level sets discovered during optimization of the {\em original} one~\cite{silver2017mastering}.
It is precisely the preimages, in parameter space, of the level sets of the first meaningful DMAPS coordinate $\phi_1$ that correspond to the neutral parameter foliation.

The remainder of the paper is structured as follows: In Section~2, we use a simple, linear, $2-$D, singularly perturbed dynamical system to bring forth the components of our data-driven framework, while retaining the connection with sloppiness terminology. Readers unfamiliar with DMAPS may want to start with the brief relevant material in the SI. The main result we illustrate in that section is the connection between singular perturbation dynamic phenomenology and data-driven detection of (what one might consider as) loss of observed dimensionality. This occurs here simultaneously in both state (model output) and parameter (model input) space. We also explore the transition region between unperturbed and singularly perturbed regimes and, finally, contrast ``data-driven {\em singular} perturbation detection'' with ``data-driven {\em regular} perturbation detection'' through another simple--yet informative--caricature. In Section~3, we move beyond caricatures to other prototypes. In Section~3a, we explore a simple kinetic example with two sloppy and one meaningful nonlinear parameter combination, readily discovered by DMAPS.
This brings up the important issue of physical understanding: the correspondence between input combinations uncovered through data mining and 
physically meaningful parameters.
That model also enables comparison of analytical and data-driven approaches (QSSA, \cite{Bodenstein}). Section~3b uses the time-honored, textbook example of Michaelis--Menten--Henri enzyme kinetics to show something we found surprising: how data-driven computations may discover parameter scalings (in this case, an alternative nondimensionalization) that better characterize the boundaries of the singular perturbation regime. Section~3 concludes with the discussion of an important subject, namely non-invertible input--output relations. We elucidate that issue using a classical chemical reaction engineering literature example, connecting the Thiele modulus (parameter) and the effectiveness factor (model output) for transport and reaction in a catalyst pellet. An important connection between the dynamics of our \emph{measurement process} and our data-mining framework arises naturally in this context. In Section~4 we summarize, and also bring up analogies with and differences from the active subspace literature: the ``effective parameters'' discovered by DMAPS are {\em nonlinear} generalizations of linear active subspaces. We conclude with a discussion of shortcomings, as well as possible extensions and enhancements of our approach.

\section{Singularly/Regularly Perturbed Prototypes}

\begin{figure}[!htp]
\centering
\begin{tabular}{cc}
\includegraphics[width=0.43\textwidth]{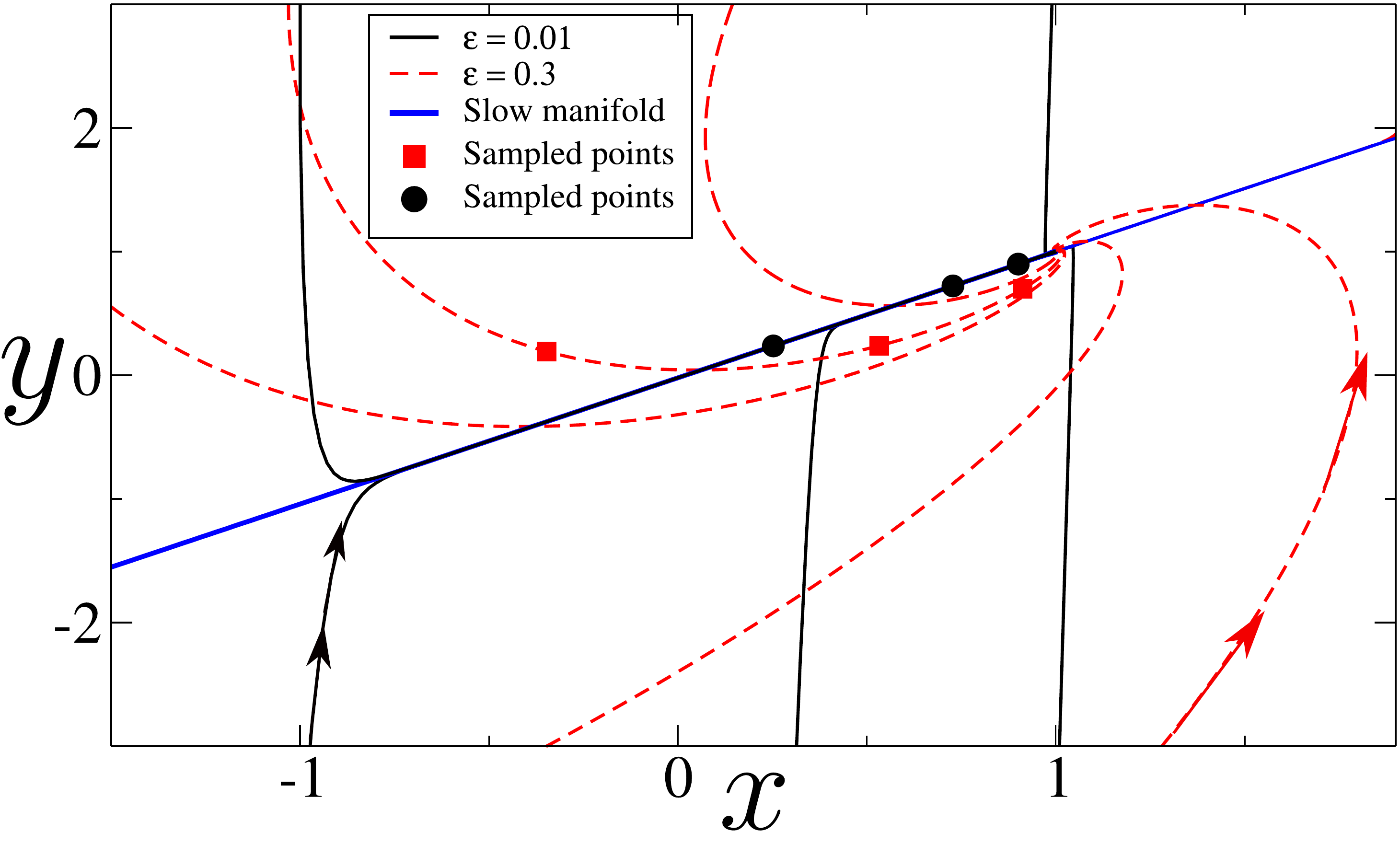} & 
\includegraphics[width=0.49\textwidth]{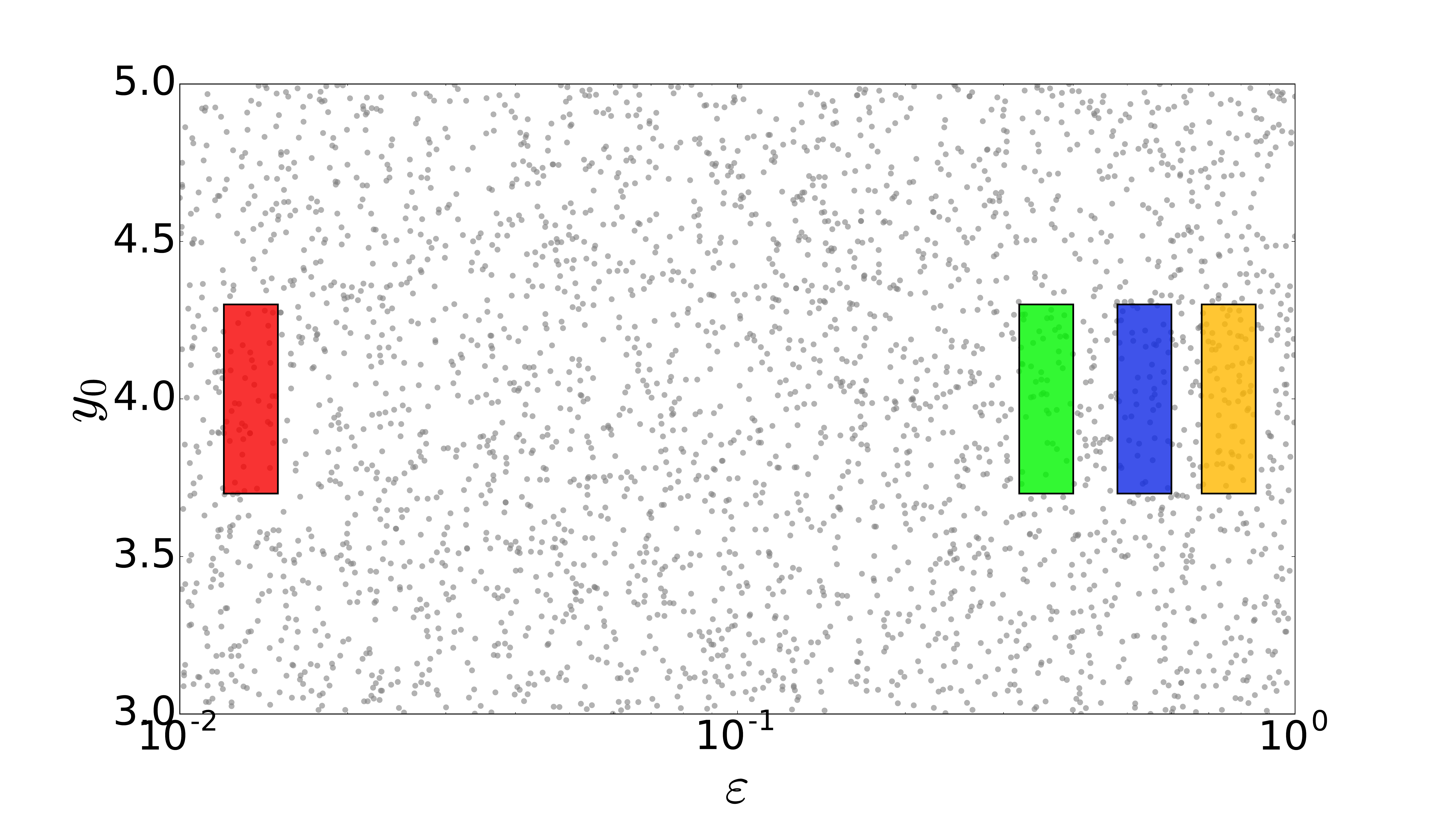} \\
(a) & (b)\\
\includegraphics[width=0.49\textwidth]{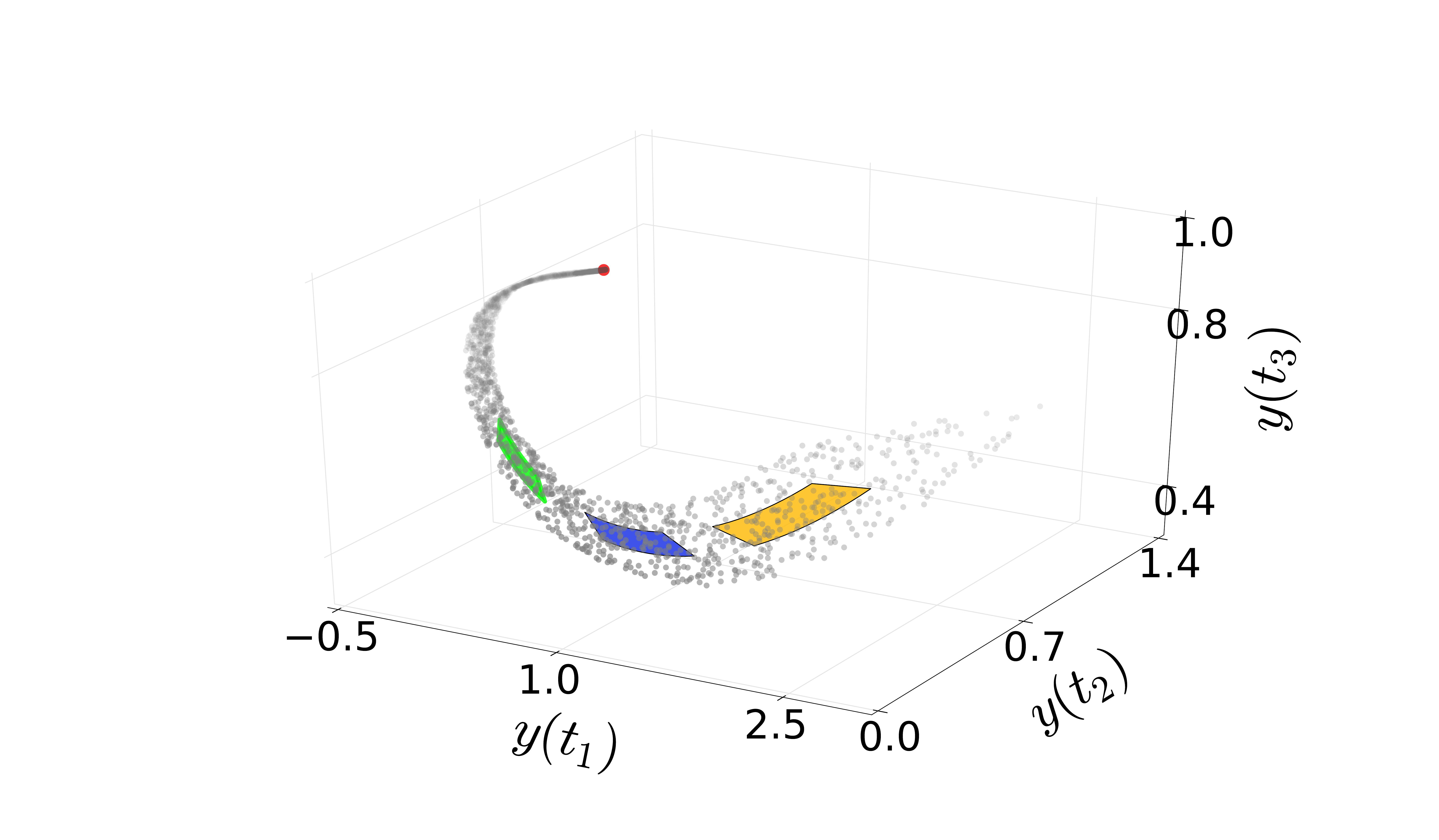}&
\includegraphics[width=0.46\textwidth]{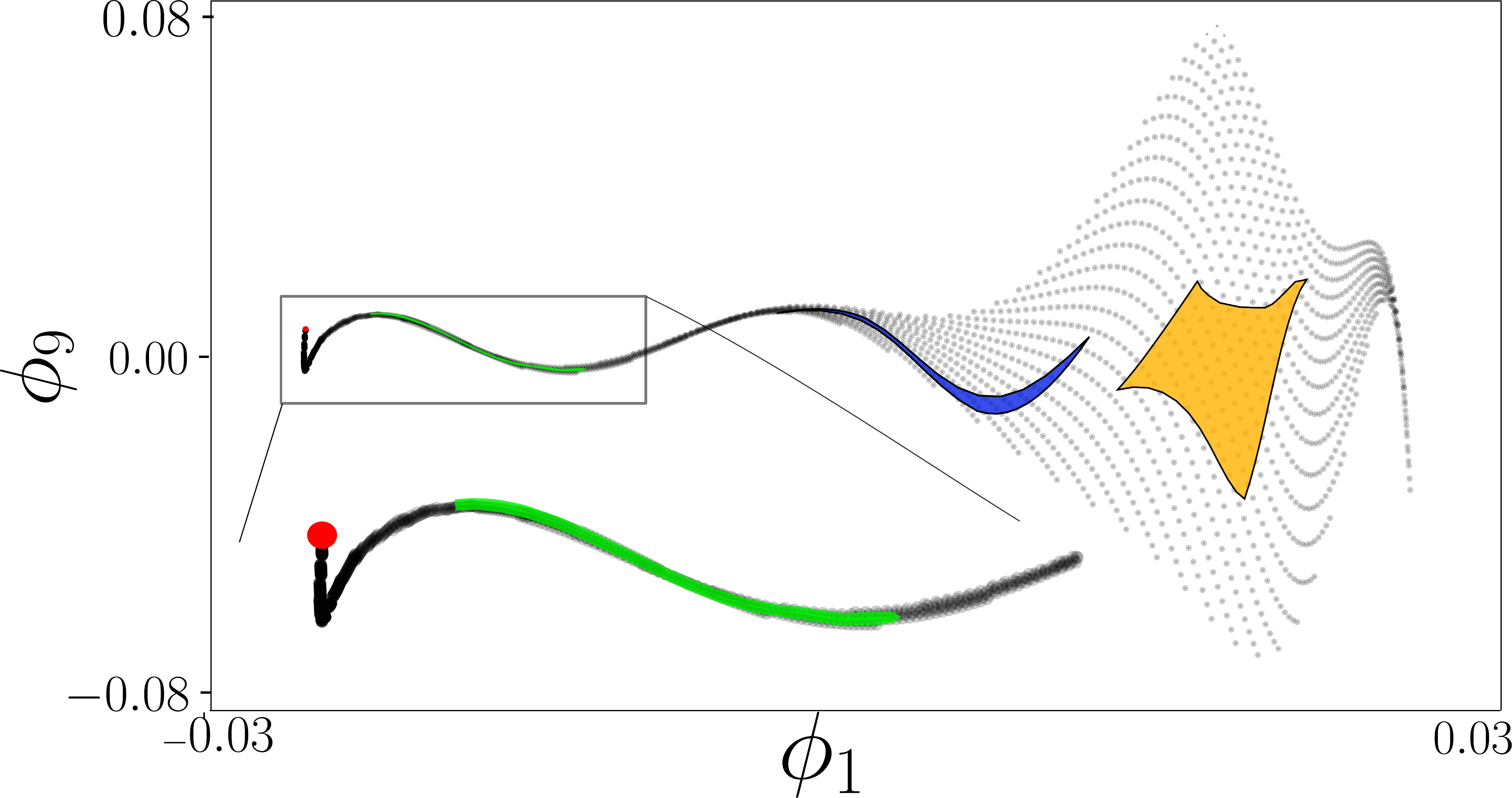}\\
(c) & (d)
\end{tabular} 
\caption{
(a) Phase portrait of \eqref{model-singpert} for a range of initial conditions $(x_0,y_0)$ and two representative $\eps-$values. 
Solid points mark states at the monitoring times $(t_1,t_2,t_3)=(0.5,1,1.5)$ 
for trajectories starting at $x_0 = -1$. 
For $\eps=0.01$, the points lie close to the slow subspace and appear $y_0-$independent; 
for $\eps=0.3$, instead, they lie off it and vary appreciably with $y_0$. 
(b) A sample of the model input space, overlaid with distinct rectangular patches. (c) Mapping of the input sample of panel~(b) to the $3-$D output space. 
The images of the random sample outline part of the model manifold, while those of the patches show the dimensionality reduction due to the singularly perturbed structure of the model. (d) Mapping of the input sample in DMAPS coordinates. The transformations from (b) to (c--d) are discussed in the text.
}
\label{fig:sing-pert-patch}
\end{figure}
To fix ideas and definitions, we start with a dynamical model
\begin{equation}
\dot{\mathbf{x}}(t \vert \prm)
=
\mathbf{v}(\mathbf{x} \vert \prm) 
\ \mbox{and} \
\mathbf{x}(t_0 \vert \prm) = \mathbf{x}_0(\prm), 
\ \mbox{where} \
\dot{\mathbf{x}}
\equiv
\frac{d\mathbf{x}}{dt}.
\label{x'=v}
\end{equation}
The vector $\mathbf{x}(t \vert \prm) \in \mathbb{R}^d$ collects the state variables at time $t$ that are observed for parameter settings $\prm \in \mathbb{R}^M$, initialization $\mathbf{x}_0(\prm)$ at $t_0$ and vectorfield $\mathbf{v}(\cdot \vert \prm) : \mathbb{R}^d \to \mathbb{R}^d$. Equation~(\ref{x'=v}) determines the system state $\mathbf{x}(t \vert \prm)$ for all times $t>t_0$, but the \emph{model output} or \emph{response} consists only of partial observations of that time course; e.g. certain state variables at specific times. {\em Observing the system} means fixing $\prm$ and initial conditions ({\em inputs}) and recording a number $N \ge M$ of scalar outputs into $\omr(\prm) \in \mathbb{R}^N$. Each input yields a well-defined output $\omr(\prm)$; as the former moves in parameter space, the latter traces out a (generically $M-$dimensional) \emph{model} manifold $\mathcal{M}$. Our \emph{data points} on this manifold are {\em input--output combinations} $(\prm,\omr(\prm)) \in \mathbb{R}^{M+N}$ and not merely the outputs; see SI. In the interest of visualization, whenever the map $\prm \mapsto \omr(\prm)$ is injective below, we only plot the projection of the model manifold on the output space $\mathbb{R}^N$.
\begin{figure}[!h]
\centering
\begin{tabular}{cc}
\includegraphics[width=0.45\textwidth,height=0.15\textheight]{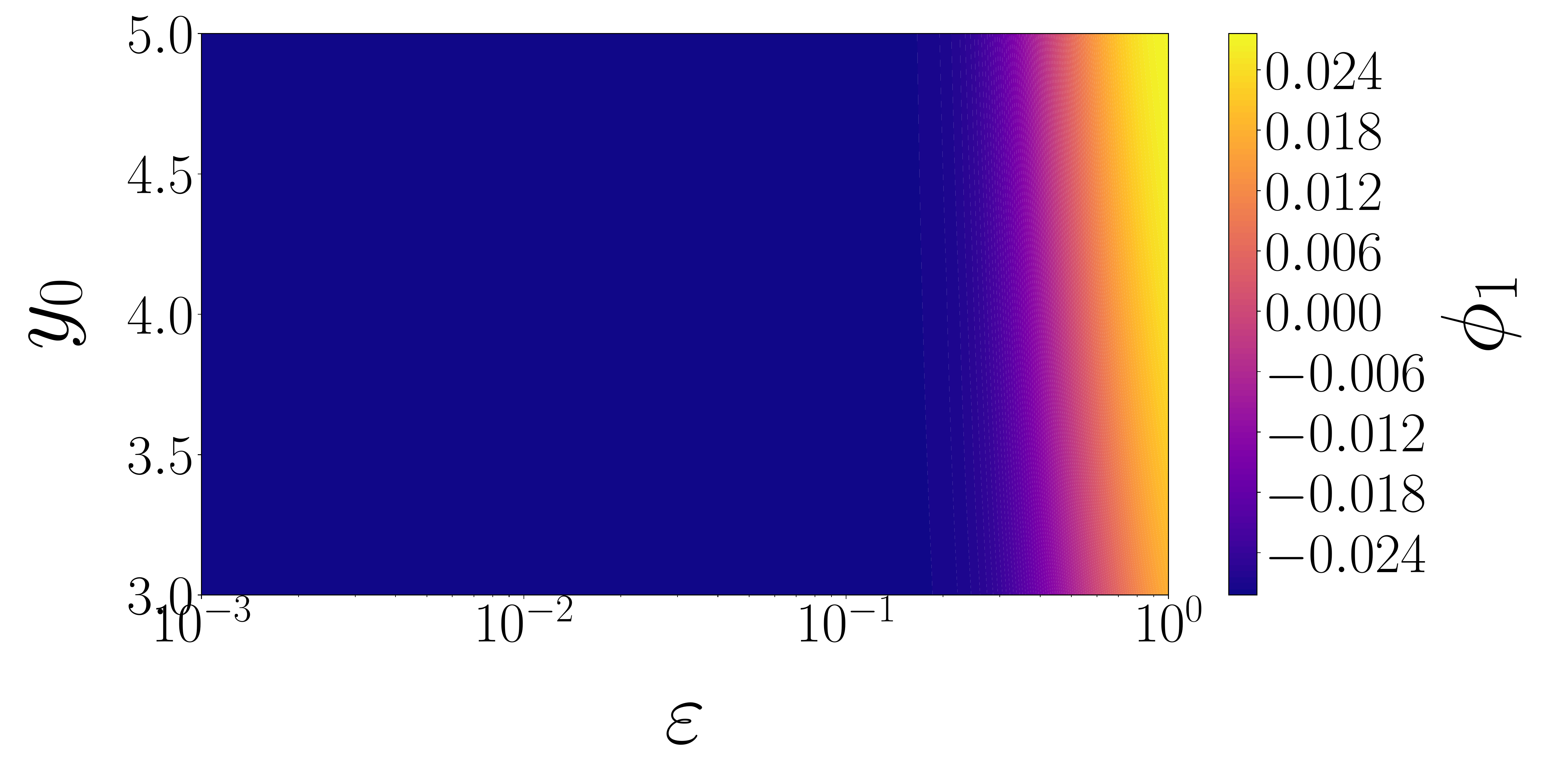}&
\includegraphics[width=0.45\textwidth,height=0.15\textheight]{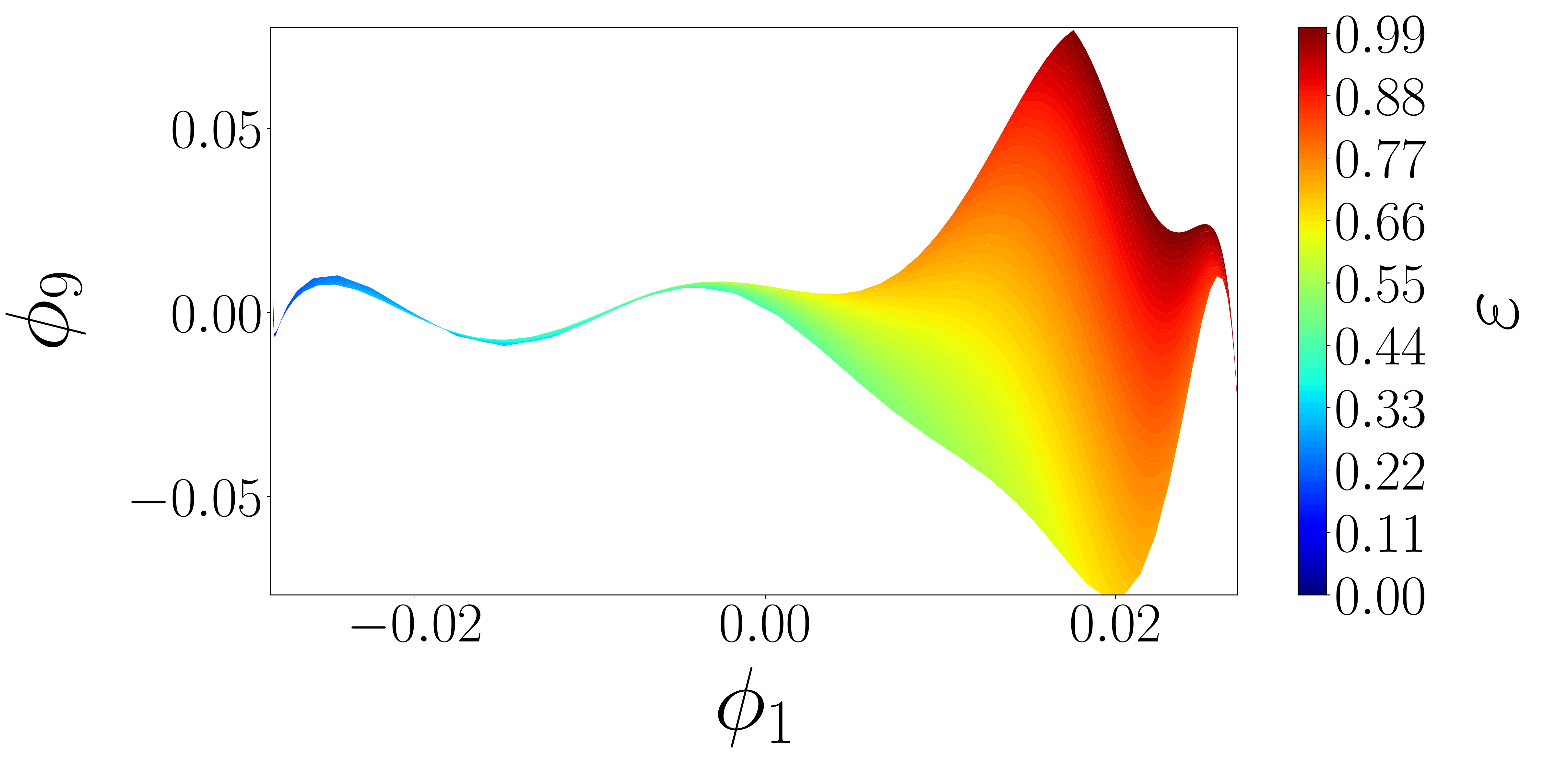}\\
(a) & (b)\\
\includegraphics[width=0.45\textwidth,height=0.15\textheight]{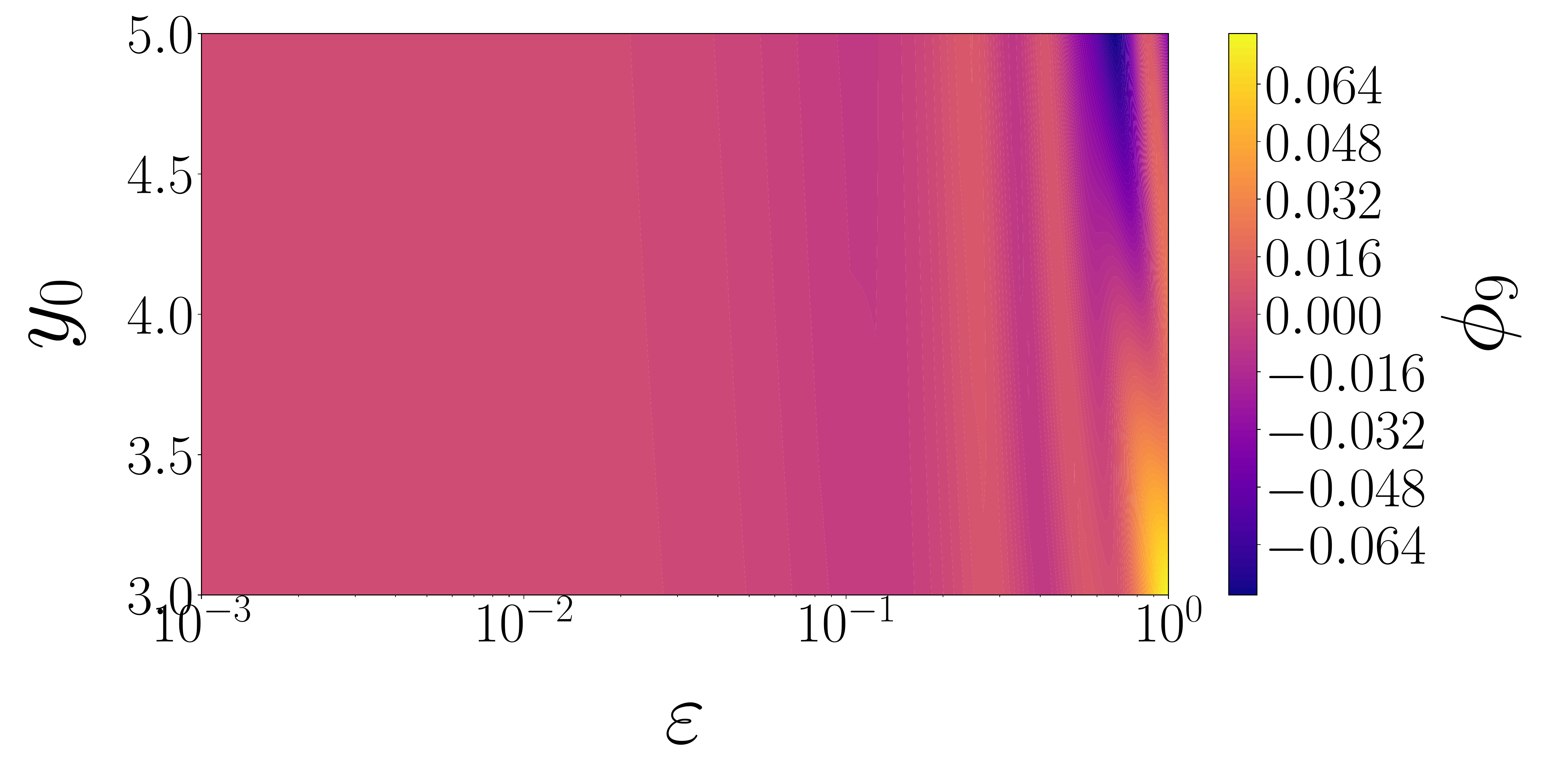}&
\includegraphics[width=0.45\textwidth,height=0.15\textheight]{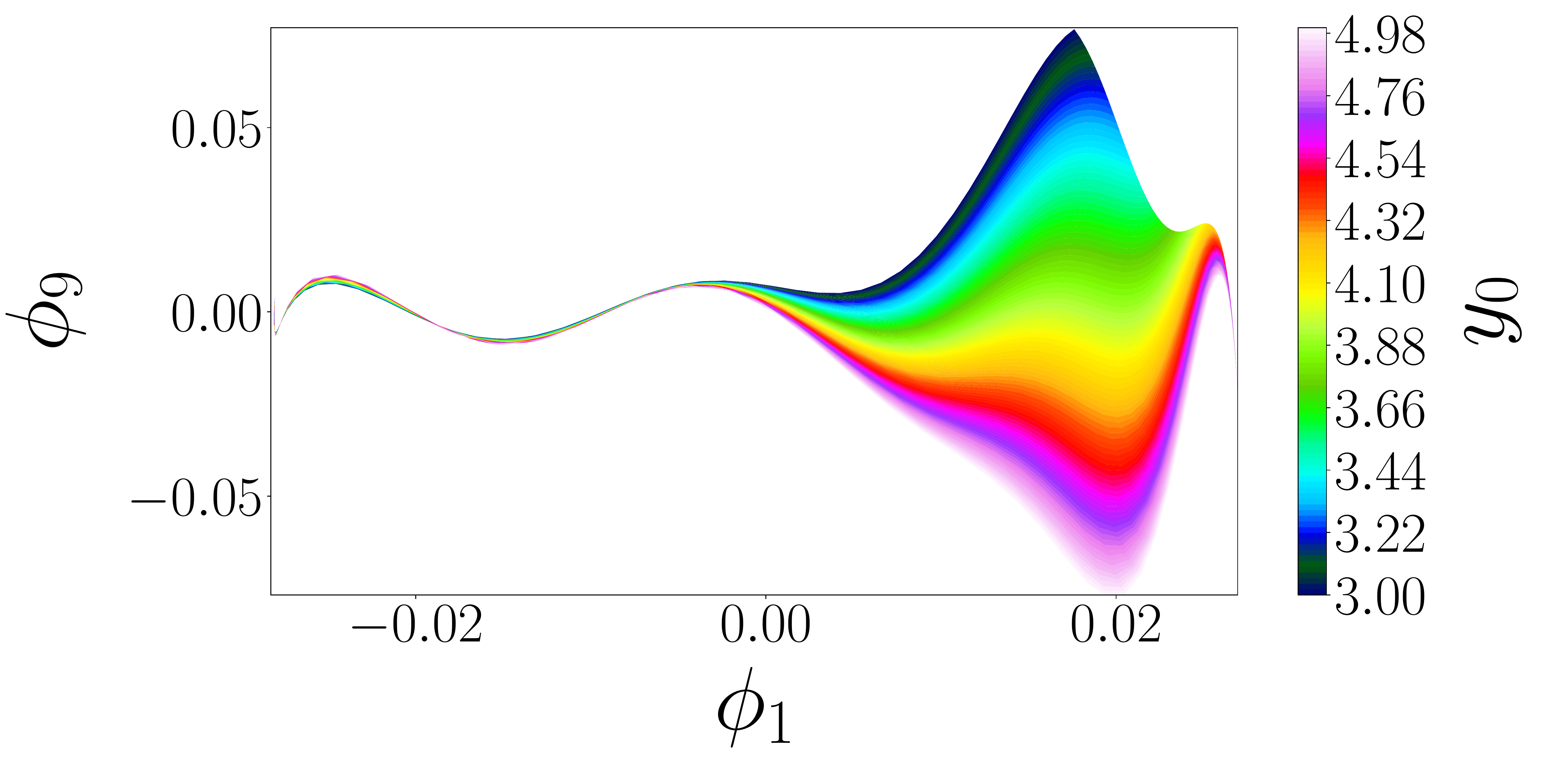}\\
(c) & (d)
\end{tabular}
\caption{
Application of DMAPS to the singularly perturbed model~\eqref{model-singpert} with an output-only-informed metric. (a, c) Input (parameter) space, coordinatized by two independent eigenmodes $\phi_1$ and $\phi_9$; (b,d) the diffusion coordinate domain (DMAPS space), coordinatized by $\eps$ and $y_0$.
All parameter settings in the singularly perturbed regime ($\eps \lessapprox 0.03$) yield effectively the same model response, $(\phi_1,\phi_9) \approx (-0.028,0)$ in diffusion coordinates, as seen by the broad monochromatic swaths at small $\eps-$values in (a,c). Intermediate $\eps-$values ($0.03 \lessapprox \eps \lessapprox 0.2$) yield an effectively $1-$D output: $\phi_9$ becomes slaved to $\phi_1$, see (b,d). Even larger $\eps-$values yield a fully $2-$D model manifold, captured by the independent color variation in all panels. The progressive decline of effective domain dimensionality is evident, as $\eps$ decreases, as is the loss of memory of the initial condition $y_0$, starting already in the $1-$D regime.
\label{fig:sing-pert-dmaps}
}
\end{figure}
To illustrate these definitions, we consider a singularly perturbed caricature chosen for its amenability to analysis,
\begin{equation}
\begin{array}{rcl}
\dot{x}
&=&
2 - x - y ,
\\
\varepsilon
\dot{y}
&=&
\hspace*{5.5mm}
x - y ,
\end{array}
\quad\mbox{with}\quad
\begin{array}{rcl}
x(0)
&=&
x_0 ,
\\
y(0)
&=&
y_0 .
\end{array}
\label{model-singpert}
\end{equation}
We also fix $x_0$ and distinct times $t_1,t_2,t_3$ (see caption of Fig.~\ref{fig:sing-pert-patch}), view both $\prm = (\eps,y_0)$ as inputs and monitor $y$; concisely, $\omr(\prm)= [y(t_1 \vert \prm) \,,\, y(t_2 \vert \prm) \,,\, y(t_3 \vert \prm)]$. The final ingredient is a \emph{metric} that provides the DMAPS kernel with a measure of closeness between different input--output combinations. For simplicity, we discuss here the {\em output-only} Euclidean metric $\|\omr(\prm) - \omr(\prm')\|$ and defer a discussion of other options to a later section. The phase portraits corresponding to two distinct $\eps-$values are plotted in Fig.~\ref{fig:sing-pert-patch}(a). For small enough $\eps$, all points on the vertical line segment $x=x_0$ ({\em fast fiber}) contract quickly to effectively the same {\em base-point} on a $1-$D invariant subspace ({\em slow subspace}) before our monitoring even begins. Memory of $y_0$ and of the boundary layer (\emph{inner solution}) is practically lost and, in the timescale of our monitoring protocol, trajectories with bounded $y_0$ shadow the evolution of that base-point and yield, with $\mathcal{O}(\eps)$ accuracy, the same output mirroring the leading order slow dynamics (\emph{outer solution}). As $\eps$ increases, the output begins to vary appreciably because the fast contraction rate decelerates and the slow invariant subspace is perturbed. However, observations still lie practically on the slow subspace and are thus insensitive to $y_0$. For even larger $\eps$, the disparity in contraction rates is relatively mild and different inputs yield visibly different trajectories; the output is jointly affected by $\eps$ and $y_0$.

This situation is evident in Fig.~\ref{fig:sing-pert-patch}(b--c), showing a randomly sampled set of inputs $\prm^{(1)} , \ldots , \prm^{(L)}$ and their simulated outputs $\omr(\prm^{(1)}) , \ldots , \omr(\prm^{(L)})$; the colored patches are meant as visual aids. The yellow patch outside the singularly perturbed regime maps into a $2-$D region of the model manifold, whereas intermediate ones (blue, green) are gradually stretched into $1-$D segments; as $\eps \downarrow 0$, or $\log(\eps) \to -\infty$, 
the effective model manifold dimensionality cascades from two to one to zero. In the $1-$D part of the model manifold and over the scales we consider, $\eps$ informs model output much more strongly than $y_0$. As $\eps\downarrow0$, the output trajectory approaches a well-defined limit -- the \emph{leading order outer solution} -- and all inputs are mapped to within $\mathcal{O}(\eps)$ of a \emph{parameter-free} output: the model manifold ``tip.'' This is evident in the red patch deep inside the singularly perturbed regime, demonstrating the {\em joint} reduction in state and in parameter space dimensionality for the scales of interest: first, the evolution law involves {\em a single} state variable, with the other slaved to it algebraically; and second, all small enough $\eps-$values produce at leading order the same, practically $y_0-$independent output.

\begin{figure}[!htp]
\centering
\begin{tabular}{cc}
\includegraphics[height=0.15\textheight]{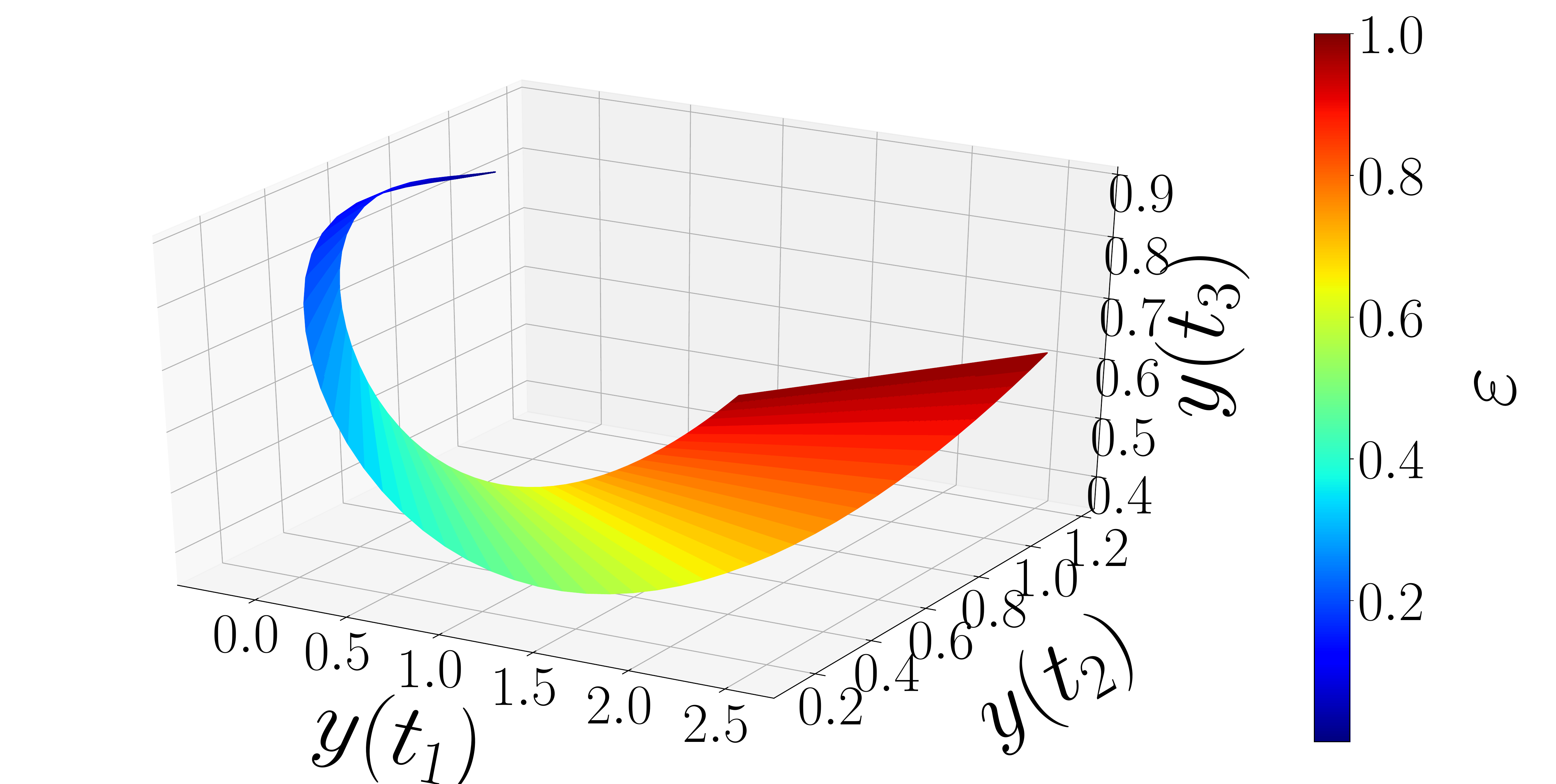}&
\includegraphics[height=0.15\textheight]{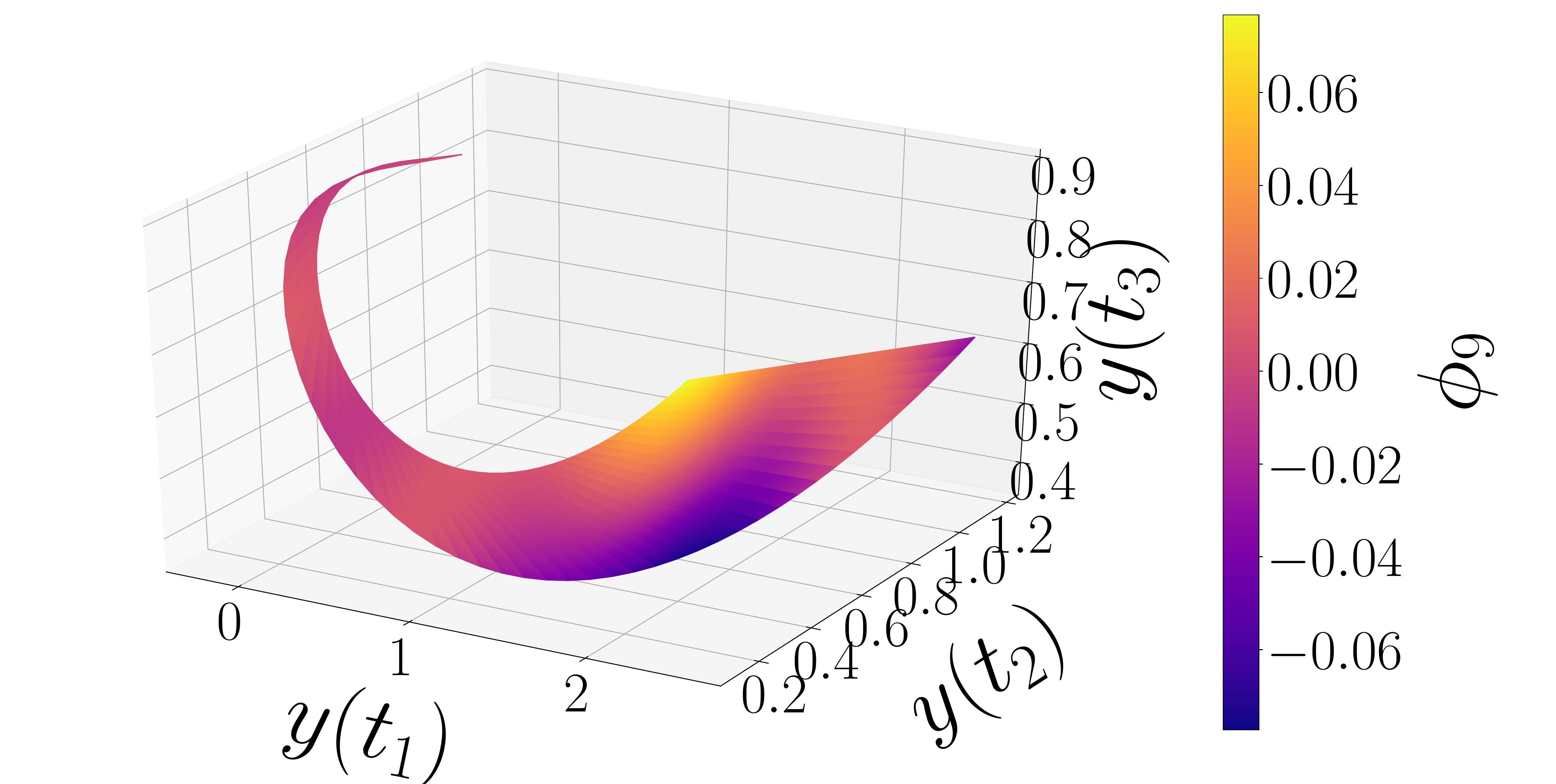} \\
(a) & (b)\\
\includegraphics[height=0.15\textheight]{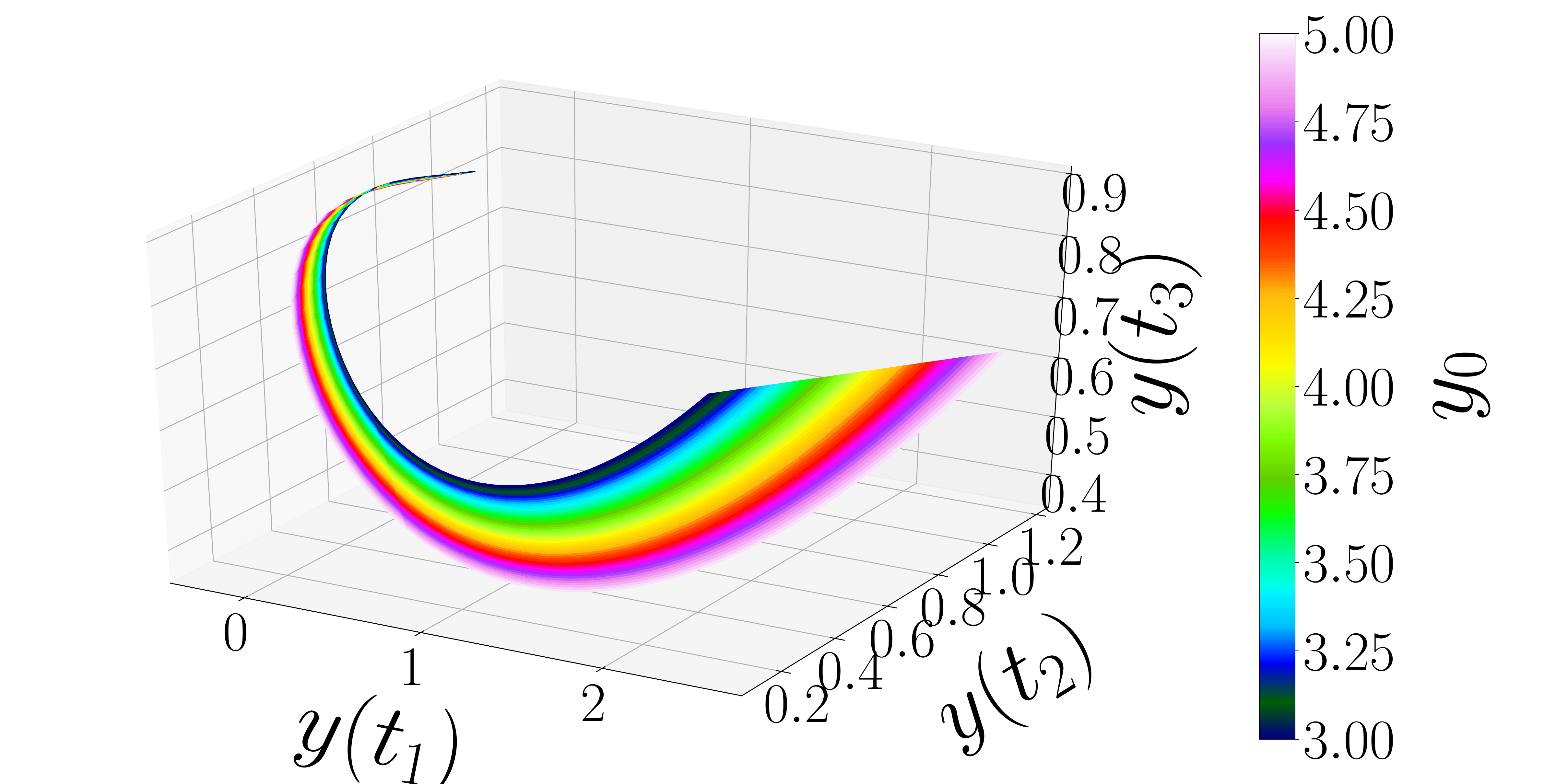} &
\includegraphics[height=0.15\textheight]{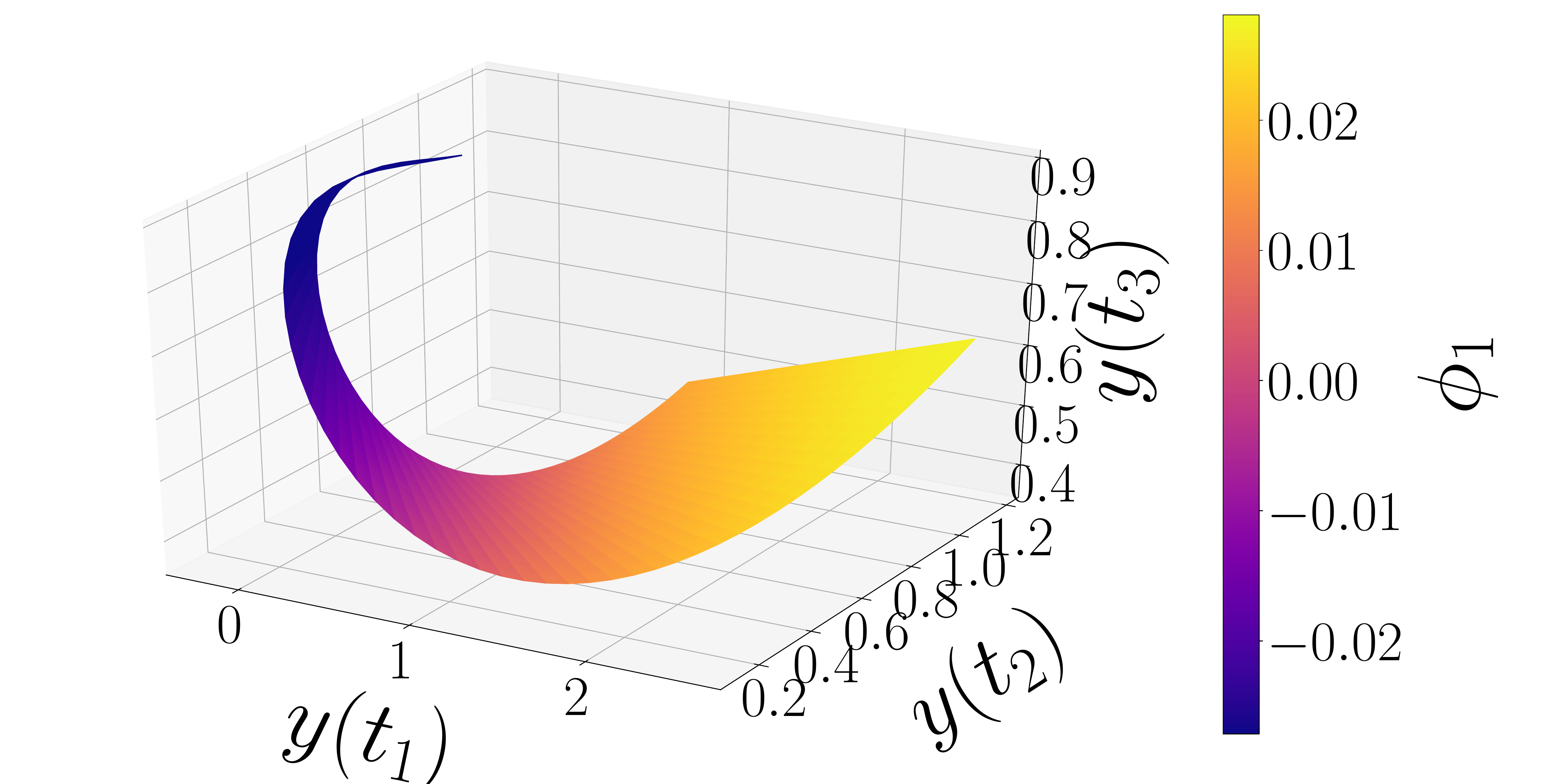} \\
(c) & (d)\\
\end{tabular}
\caption{
(a$-$b) Model manifold colored by the model parameters (inputs). For large $\eps$, the model manifold is evidently coordinatized by $(\eps,y_0)$. As $\eps$ decreases, the system loses memory of the initialization $y_0$ and model responses for different $y_0$ bundle together. In the singularly perturbed regime (deep blue in panel~(a)), all memory of $y_0$ has been lost. (c$-$d) Model manifold colored by leading independent DMAPS coordinates. Evidently, $\phi_1$ tracks $\eps$ well, with the regime $\eps \ll 1$ corresponding to $\phi_1 \approx -0.028$. The coordinate $\phi_9$ is transverse to $\phi_1$ in the $2-$D regime but, as $\eps$ decreases, becomes slaved to it and dimension reduction occurs.
\label{fig:sing-pert-mm}
}
\end{figure}

To glean the information above by data mining, we apply DMAPS (see \cite{juan_m_bello_rivas_2017_581667} for the code) with an {\em output-only informed kernel} to the dataset 
and obtain the re-coordinatization $(\prm^{(\ell)},\omr(\prm^{(\ell)})) \mapsto (\phi_1^{(\ell)} , \phi_9^{(\ell)})$. Here, $\phi_1 , \phi_9 \in \mathbb{R}^L$ are independent eigenvectors of the DMAPS kernel, i.e. discretizations (on the dataset) of eigenfunctions of the Laplace--Beltrami operator defined on the model manifold $\mathcal{M}$. As such, they describe diffusive eigenmodes whose level sets endow $\mathcal{M}$ with an intrinsic, nonlinear coordinate system. The domain of that coordinate system ({\em DMAPS space}) is shown in Fig.~\ref{fig:sing-pert-patch}(d). Here also, the stretching factor increases and the dimensionality of the mapped patches cascades, as we progress into the singularly perturbed regime. The preimage, in parameter space, of that coordinate system is shown in Fig.~\ref{fig:sing-pert-dmaps}(a,c), allowing us to define distances between inputs in terms of the outputs they generate. Figure~\ref{fig:sing-pert-dmaps}(b,d) portrays complementary images, namely the coordinatization of DMAPS space in terms of the inputs $\eps$ and $y_0$. Finally, Fig.~\ref{fig:sing-pert-mm}(a--b) and Fig.~\ref{fig:sing-pert-mm}(c--d) show the model manifold colored by the inputs as well as by the diffusion eigenmodes.

These figures relate input, output and DMAPS domains to model dynamics and suffice to reproduce our earlier observations on model output.
In the $2-$D part of the DMAPS domain and of the model manifold, distinct points on the latter correspond to distinct diffusion coordinates and distinct inputs $\eps$ and $y_0$; see Fig.~\ref{fig:sing-pert-dmaps}. As $\eps$ decreases, however, the dependence on $y_0$ becomes attenuated and the output controlled by $\eps$ alone. In the terminology introduced earlier, $y_0$ becomes \emph{sloppy} and both the DMAPS domain and the model manifold transition to a $1-$D regime parameterized by $\eps$; the $y_0-$values span an ever-diminishing width. In this regime, the {\em level sets} of the eigenmodes visibly align with each other, both in input space and on the model manifold; see Fig.~\ref{fig:sing-pert-dmaps}(a,c) and Fig.~\ref{fig:sing-pert-mm}(c--d). Finally, as $\eps\downarrow0$, all parameter settings converge to the same $(\phi_1,\phi_9)-$value, as the output converges to the ``tip'' of the model manifold and of the DMAPS domain.

In summary, our output-informed application of DMAPS parameterizes the input--output combinations comprising the dataset in a manner indicative of how model inputs dictate model outputs. The parameterization applies primarily to the output component of the dataset, but it can be pulled back to yield a simultaneous, {\em consistent} re-parameterization of the input component. This showcases the main contribution in this paper: a way to intuit system properties by parameterizing the input--output relation through the geometry of the manifold that collects model inputs and model outputs, as encoded in eigenfunctions of the Laplace--Beltrami operator.
%
\subsection{Regularly perturbed prototype}
\begin{figure}[h]
\centering
\begin{tabular}{cc}
\includegraphics[width=0.45\textwidth,height=0.2\textheight]{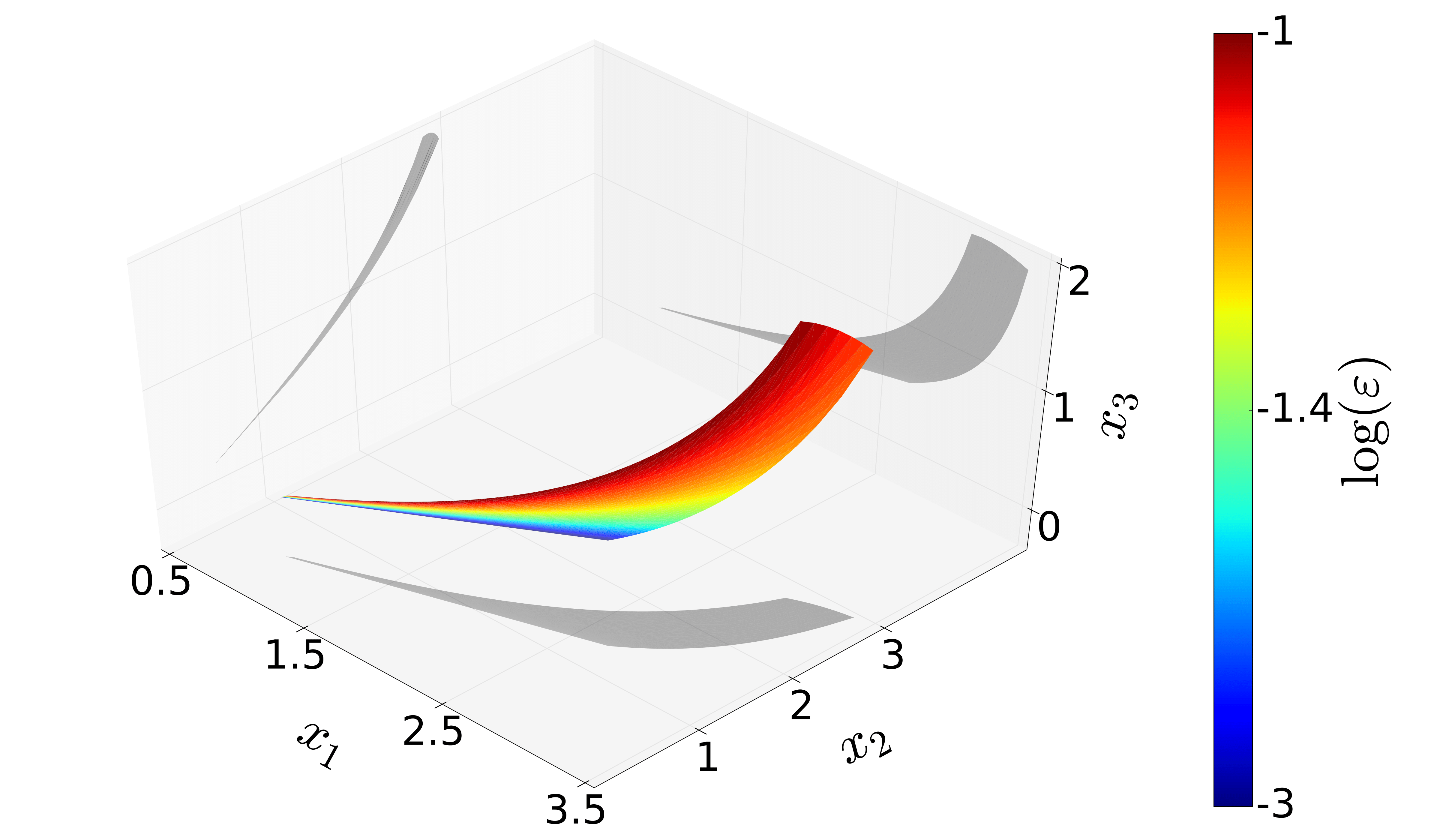}&
\includegraphics[width=0.45\textwidth,height=0.2\textheight]{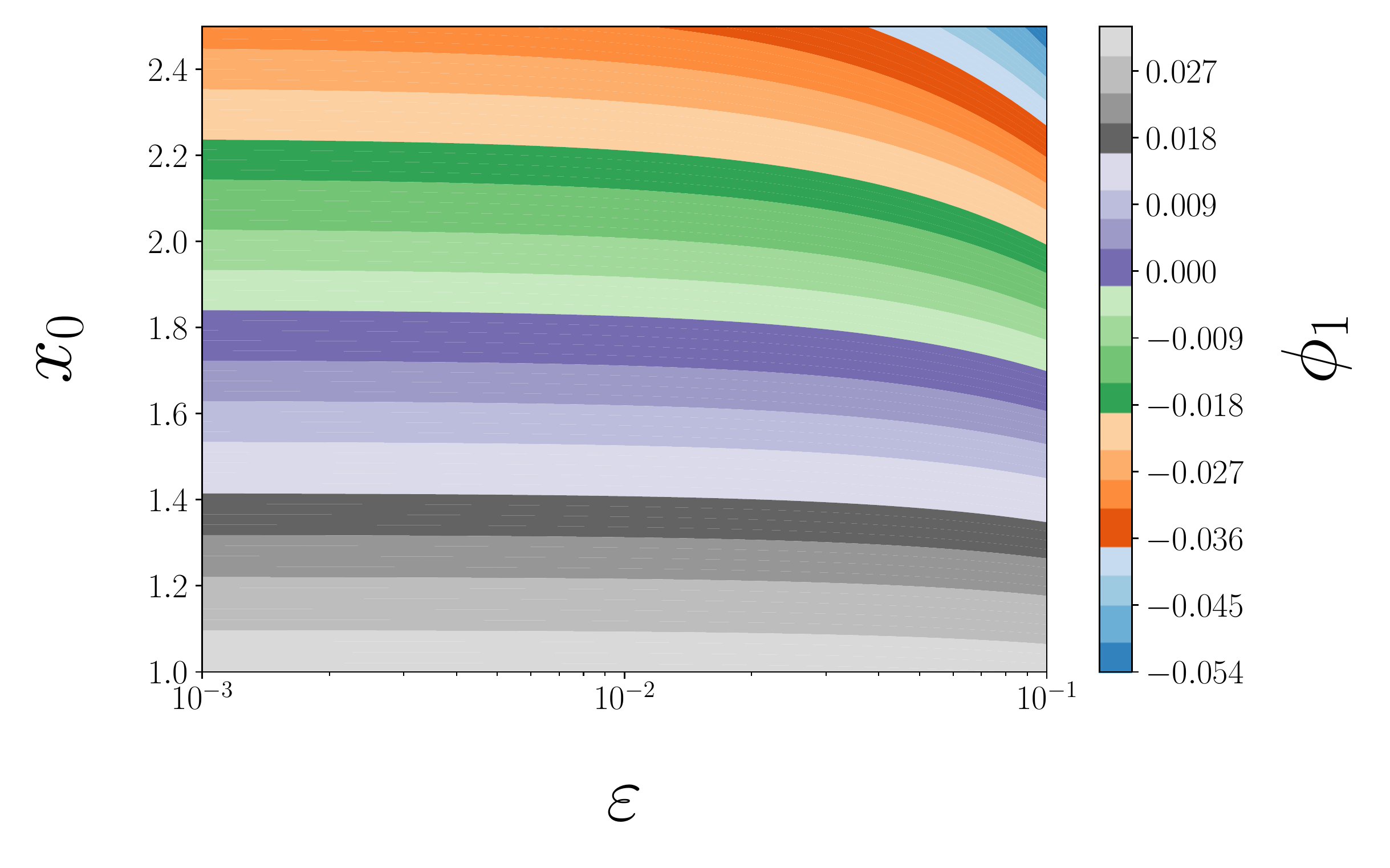} \\
(a) & (b)
\end{tabular}
\caption{
Model manifold and parameter space for the regularly perturbed system~\eqref{1D-model-regpert} with inputs $(\eps,x_0)$ and monitoring times $(t_1,t_2,t_3) = (0.25, 1.0, 1.75)$. (a) Model manifold colored by $\eps$ and its projections on the coordinate planes. For large $\eps$, the system is fully $2-$D (orange/red region). As $\eps$ decreases, the model response becomes increasingly determined by $x_0$ alone. Contrary to Fig.~\ref{fig:sing-pert-dmaps}, where the $\eps \ll 1$ regime was mapped to a limiting point, the limiting submanifold here is $1-$D (blue straight line). (b) Parameter space colored by the DMAPS coordinate $\phi_1$. DMAPS visibly captures the importance of $x_0$, as $\eps \downarrow 0$: the parameterization varies in the $x_0$ direction and remains unchanged along lines of constant $x_0$.
}
\label{fig:regpert}
\end{figure}
In the {\em singularly perturbed} prototype discussed above, we noted the {\em simultaneous} loss of output sensitivity to (certain) initial conditions and parameters, as $\eps \downarrow 0$. Additionally, we demonstrated how this system behavior can be intuited by mining input--output data with the help of DMAPS. Figure~\ref{fig:regpert}(b) shows the result of applying the same methodology to the {\em regularly perturbed} example
\begin{equation}
\dot{x} = -x + \eps x^3,
\quad\mbox{with}\
x(0) = x_0 .
\label{1D-model-regpert}
\end{equation}
Here also, we view $\prm=(\eps,x_0)$ as parameters and monitor the system state $x(t \vert \prm)$ at distinct times.

Similarly to the singularly perturbed model, the model output approaches a well-defined, limiting response in the asymptotic regime $\eps \downarrow 0$. Yet, in this case that response {\em remains strongly dependent} on $x_0$: distinct initial conditions yield distinct outputs even for $\eps \downarrow 0$, as seen plainly in Fig.~\ref{fig:regpert}. In panel~(a), the limiting edge $\eps = 0.001$ is seen to outline a $1-$D boundary of the full $2-$D model manifold, instead of a {\em point} as was the case for the singularly perturbed model. That same edge is seen to be parameterized by $\phi_1$, in panel~(b), rather than correspond to a single $\phi_1-$value. This result is clearly underpinned by the uniform convergence of the trajectory $x(t \vert \prm)$ to $x(t \vert 0,x_0) = x_0 \mathrm{e}^{-t}$, which is $\eps-$free but depends strongly on $x_0$ and defines the aforementioned $1-$D model manifold boundary. This {\em regular perturbation} behavior, and specifically the \emph{lack} of dimensionality reduction in terms of initial conditions, generalizes directly to higher state and parameter space dimensions.

\section{Beyond caricatures}
%
\subsection{The ABC model}
\begin{figure}[H]
\centering
\begin{tabular}{ccc}
\includegraphics[width=0.2\textwidth]{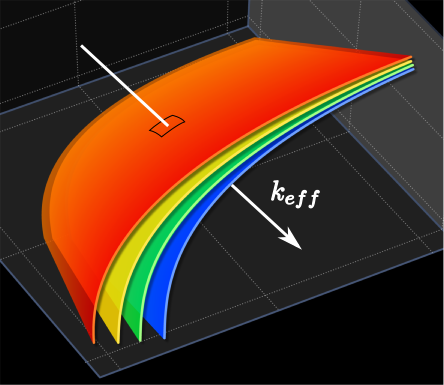}&
\includegraphics[width=0.32\textwidth]{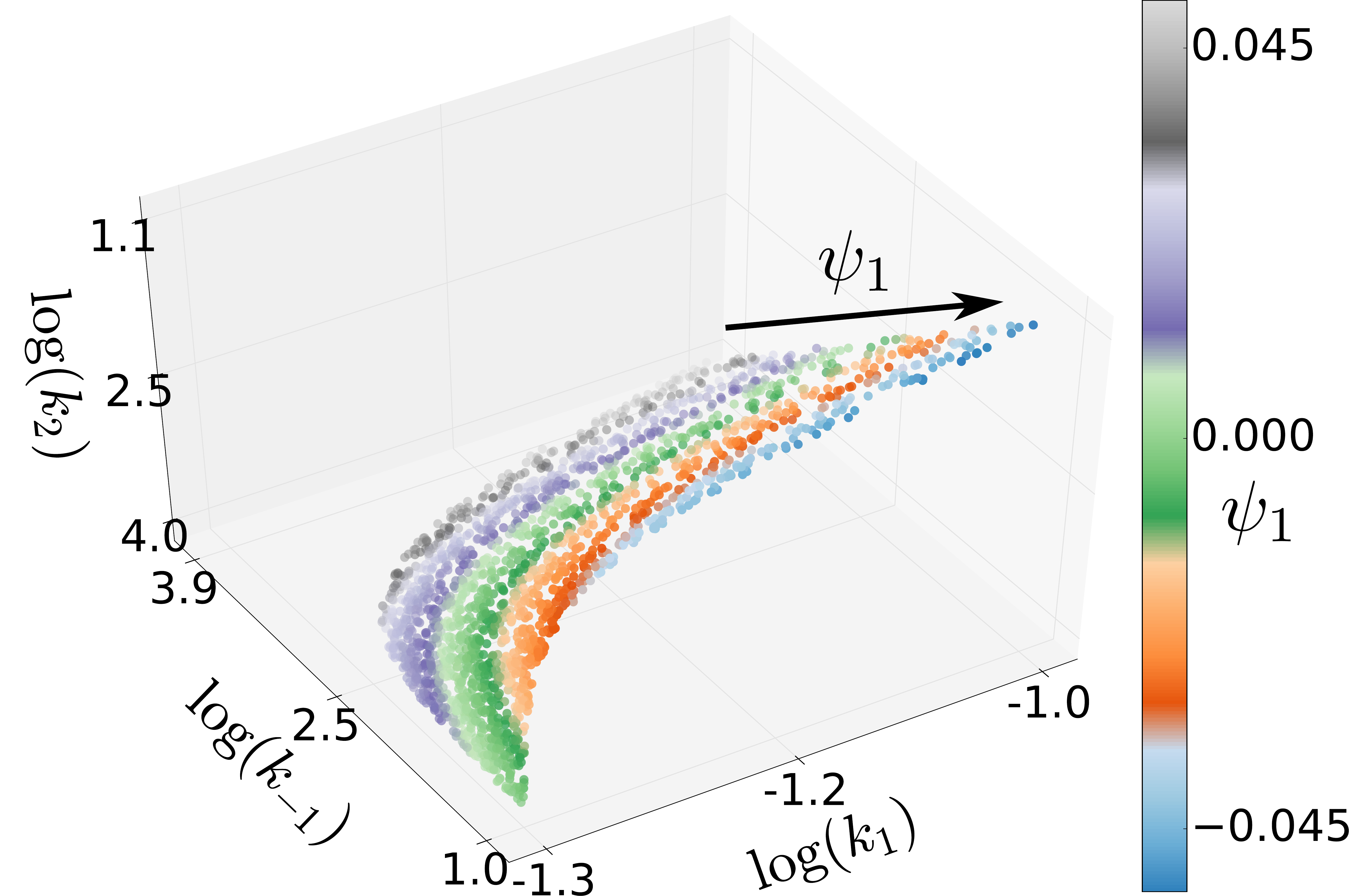}&
\includegraphics[width=0.32\textwidth]{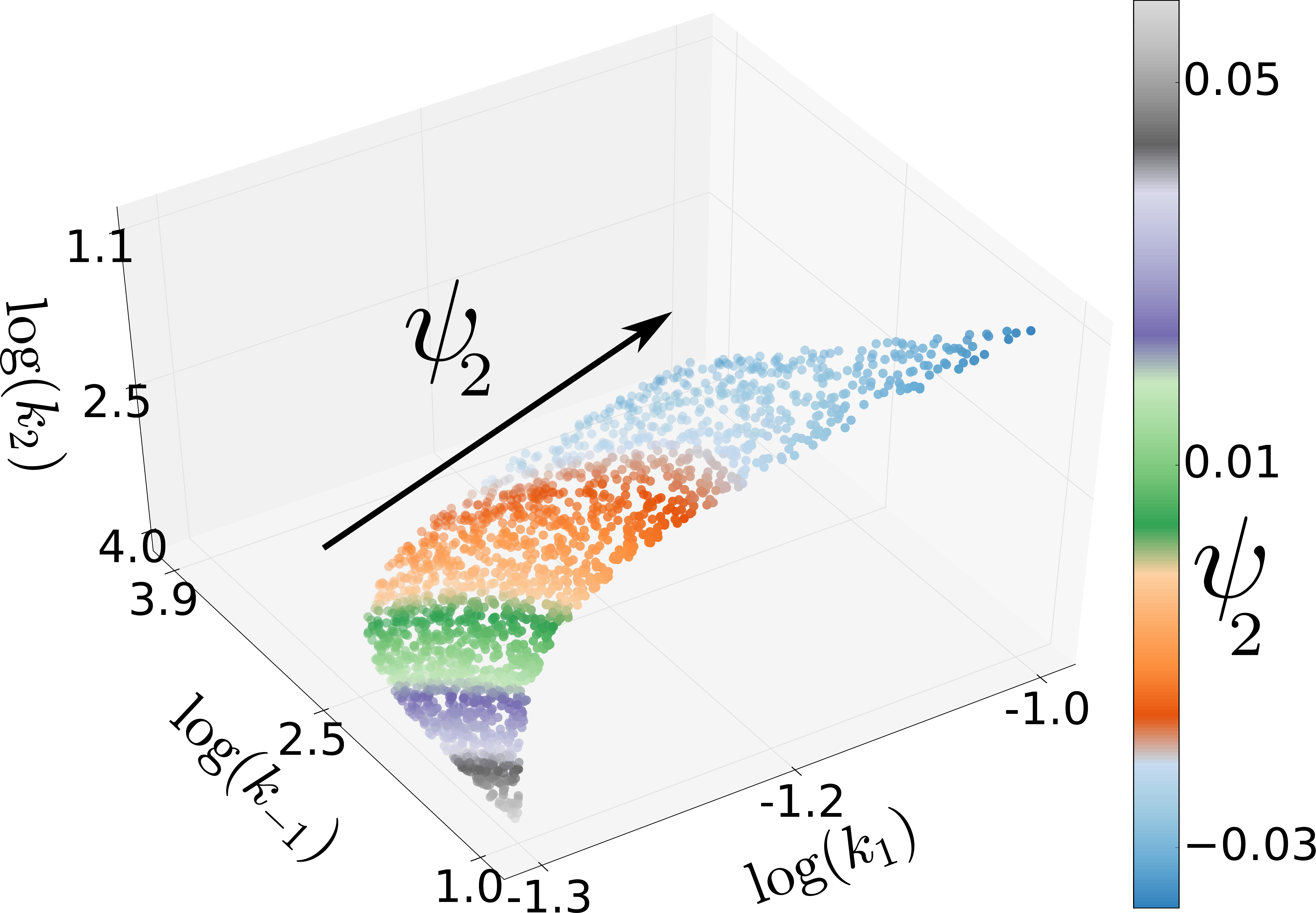}\\
(a) & (b) & (c) \\
\includegraphics[width=0.32\textwidth]{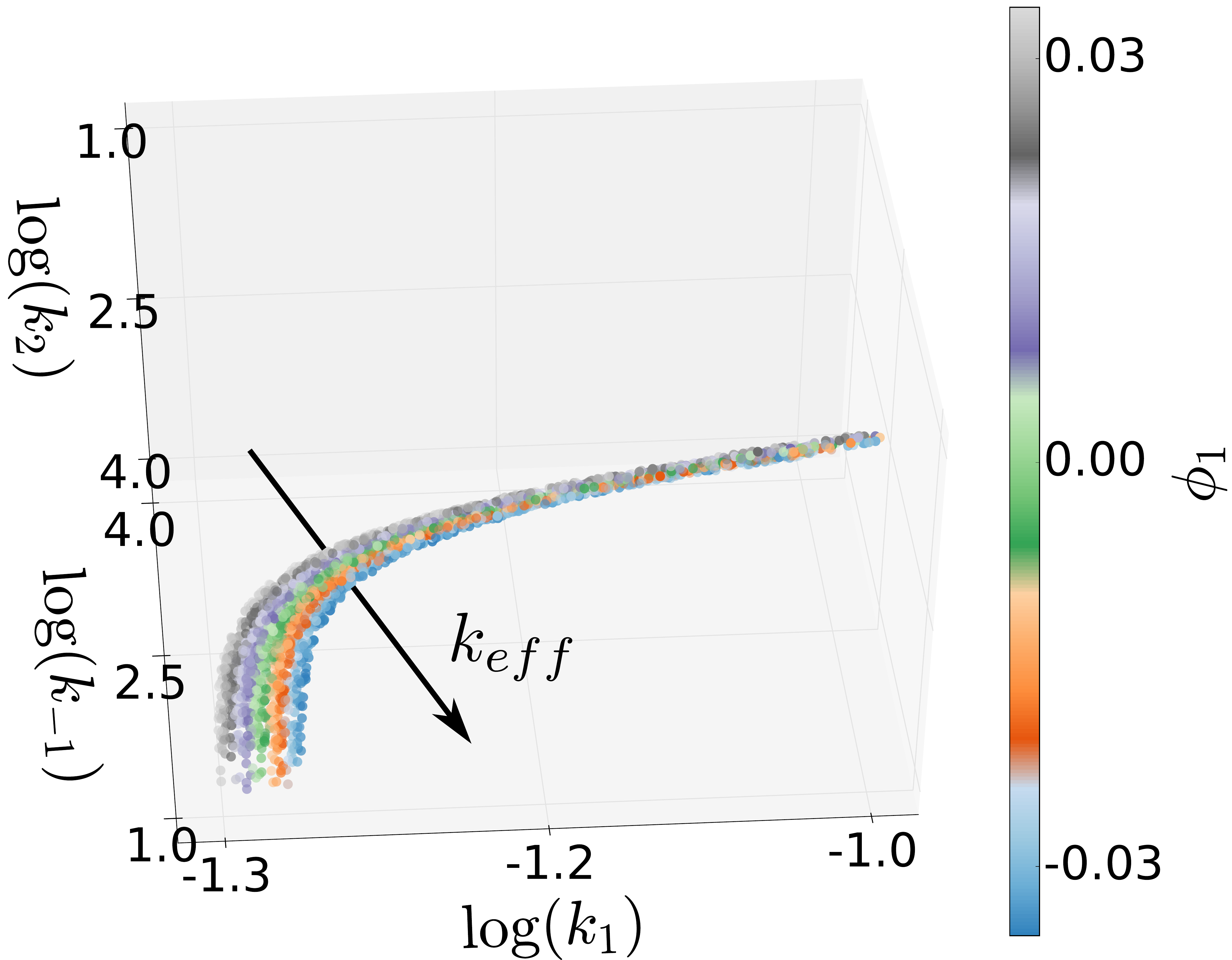}&
\includegraphics[width=0.32\textwidth]{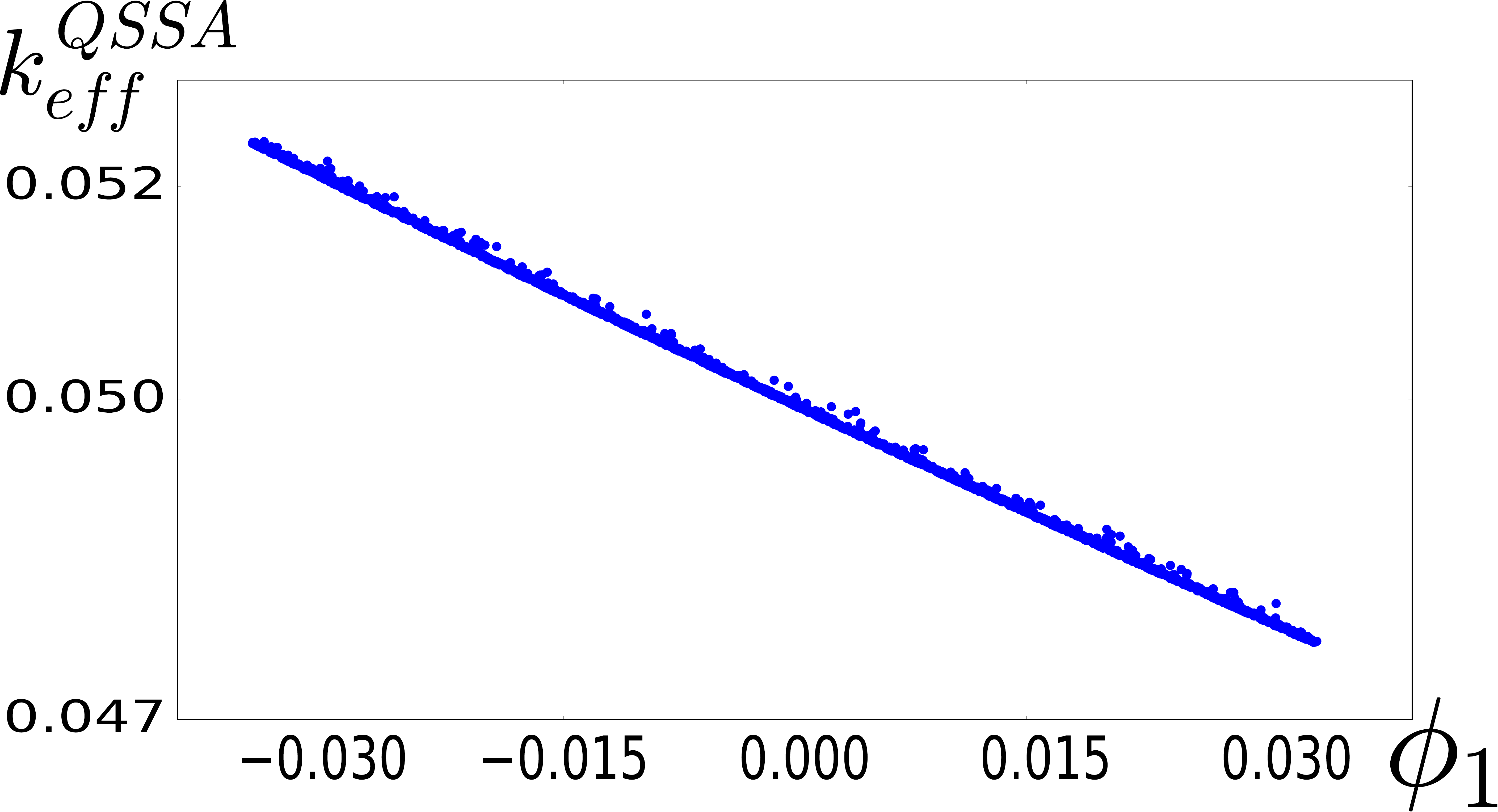}&
\includegraphics[width=0.3\textwidth]{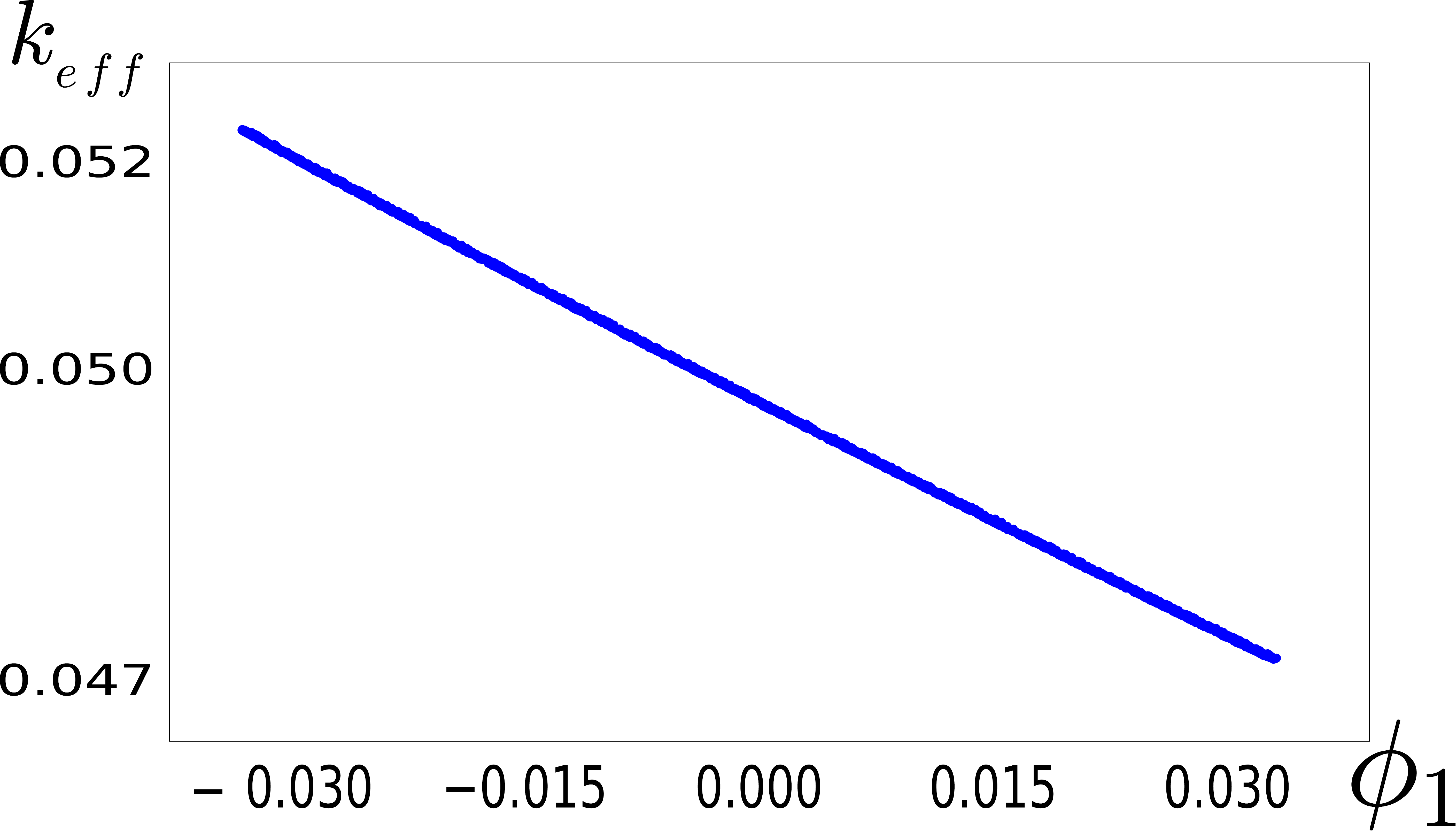}\\
(d) & (e) & (f) \\
\end{tabular}
\caption{
Data-driven detection and characterization of the effective parameter $\keff$ for \eqref{mech:abc} using DMAPS. 
All datasets were obtained by presetting the reference output $\omr^*=\omr(\ps)$ with $\ps=(10^{-1},10^3,10^3)$ 
and a specific tolerance $\delta > 0$ and, then, sampling log-uniformly a rectangular domain in input space to
retain inputs satisfying $\|\omr(\prm) - \omr^*\| < \delta$. 
(a) Illustration of level sets of $\keff$ in parameter space $(k_1,k_{-1},k_2)$. 
Equation \ref{eq:abc-redux} 
dictates that points on each same colored foil exhibit nearly identical model responses.
(b--c) Dataset for $\delta=0.001$; this is practically the $2-$D surface $\keff(\prm)=\keff(\ps)$ 
of (almost) perfect fits. An application of input-only DMAPS on this set reveals its $2-$D nature and 
coordinatizes it through the eigenfunctions $\psi_1,\psi_2$. 
(d) Dataset for $\delta = 0.1$, colored by the first output-only DMAPS eigenfunction $\phi_1$. 
DMAPS clearly discovers the effective parameter, as $\phi_1$ remains effectively 
constant on level sets of $\keff$. This striking one-to-one relation is evident in 
panel~(f), in terms of the $(\phi_1,\keff)-$coordinates. 
The same dataset plotted in $(\phi_1,\keff^\mathrm{QSSA})-$coordinates (panel e) is, 
by contrast, visibly noisier.
}
\label{fig:abc}	
\end{figure}

Having examined simple singularly and regularly perturbed models, we turn our attention to the \emph{data-driven detection} of an effective parameter in a paradigmatic chemical reaction network. We specifically consider the three-species, analytically tractable system (SI)
\begin{align}
\mathrm{A}
\xrightleftharpoons[k_{-1}]{k_1}
\mathrm{B}
\xrightarrow[]{k_2}
\mathrm{C} ,
\ \mbox{initialized with} \
 A_0=1 , B_0=C_0=0 .
\label{mech:abc}
\end{align}
The Quasi-Steady-State Approximation (QSSA; \cite[Ch.~5]{rawlings_chemical_2004}) for mechanism~(\ref{mech:abc}) reads
\begin{align}
C(t\vert\prm) = 1 - e^{-\keff^\mathrm{QSSA} t} ,
\ \mbox{where} \
\keff^\mathrm{QSSA} = \frac{k_1k_2}{(k_{-1} + k_2)} ,
\end{align}
and is valid for $k_1 \ll k_2$ \cite{rawlings_chemical_2004}. A detailed analysis, however (see SI), establishes that the approximate solution is actually 
\begin{align}
C(t\vert\prm)
=
1 - e^{-\keff t} ,
\quad\mbox{where}\
\keff
=
\frac{k_1 k_2}{k_{-1} + k_1 + k_2} .
\label{eq:abc-redux}
\end{align}
To detect this dimensionality reduction and ``discover'' $\keff$ in a data-driven manner, we view the kinetic constants as inputs, $\prm = (k_1,k_{-1},k_2)$, and monitor product concentration at preset times, $\omr(\prm) = [C(t_1 \vert \prm) , \ldots , C(t_5 \vert \prm)]$. Then, we fix a model output in the regime of applicability of \eqref{eq:abc-redux} and mine sampled parameter settings with outputs ``similar'' to that reference response. Here, we used as reference the output corresponding to parameter settings $\ps=(10^{-1},10^3,10^3)$ and measured similarity in the Euclidean sense, retaining sampled points $\prm$ satisfying $\|\omr(\prm) - \omr(\ps)\| < \delta$ for some $\delta>0$. Figure~\ref{fig:abc}(b$-$f) examines two nested such ``good datasets'', one of almost perfect fits ($\delta = 10^{-3}$; Fig.~\ref{fig:abc}(b-c)) and another of less good fits ($\delta = 10^{-1}$; Fig.~\ref{fig:abc}(d)). Data-mining the ``zero residual level set'' in 3-D parameter space with an {\em input-only} informed DMAP metric confirms its 2-D nature. The data-driven coordinatization of the full input space by output-only $\phi_1$ and input-only $(\psi_1,\psi_2)$ decomposes the space in a manner tuned to model output. A related data-processing of good fits using linear PCA was performed, e.g., as in \cite{achard_complex_2006} for a neuron model; clearly, linear PCA here would give the erroneous impression of full-dimensionality due to manifold curvature.
This result is valid in the input regime $k_1 k_2 \ll (k_1 + k_{-1} + k_2)^2$, that extends the QSSA, and $\keff$ is the effective parameter (approximately) determining the output. This expression represents a reduction of input space from $3-$D to $1-$D; the foliation of parameter space by the (nonlinear) level sets of $\keff$ is shown schematically in Fig.~\ref{fig:abc}(a).
The set of parameter settings with outputs within $\delta = 10^{-1}$ of the reference output is clearly $3-$D and visibly composed of level sets of $\keff$ spanning an appreciable $\keff$ range. An application of DMAPS with the Euclidean, output-only-informed similarity measure reveals the existence of a single effective parameter {\em without recourse to an analytic expression}. Indeed, the DMAPS coordinate $\phi_1$ traces $\keff$ accurately, see Fig.~\ref{fig:abc}(d,f). Note, for comparison, that $k_\mathrm{eff}^\mathrm{QSSA}$ is a worse predictor of model output, see Fig.~\ref{fig:abc}(e). It follows that level sets of $\phi_1$ in parameter space give (\emph{almost}) \emph{neutral sets}, i.e. level sets of $\keff$ whose points yield indistinguishable outputs. An algorithm to explore parameter space effectively would march along $\phi_1$, whereas sampling parameter inputs at constant $\phi_1$  would allow one to map out level sets of $\keff$. This can be of particular utility in multi-objective optimization~\cite{silver2017mastering}, where a \emph{second} objective can be optimized on the set $\keff = \keff(\ps)$ optimally fitting the data $\omr^*=\omr(\ps)$.  

Fig.~\ref{fig:abc}(e,f) raises the crucial issue of physical interpretation of the effective parameters discovered through data mining. Although such data-driven parameters are not expected to be physically meaningful, the user can {\em post-process their discovery} by formulating and testing hypotheses on whether they are one-to-one with (i.e., encode the same information as) physically meaningful parameters.   

\subsection{Michaelis--Menten--Henri (MMH)}
%
\begin{figure}[h]
\centering
\begin{tabular}{cc}
\includegraphics[width=0.45\textwidth]{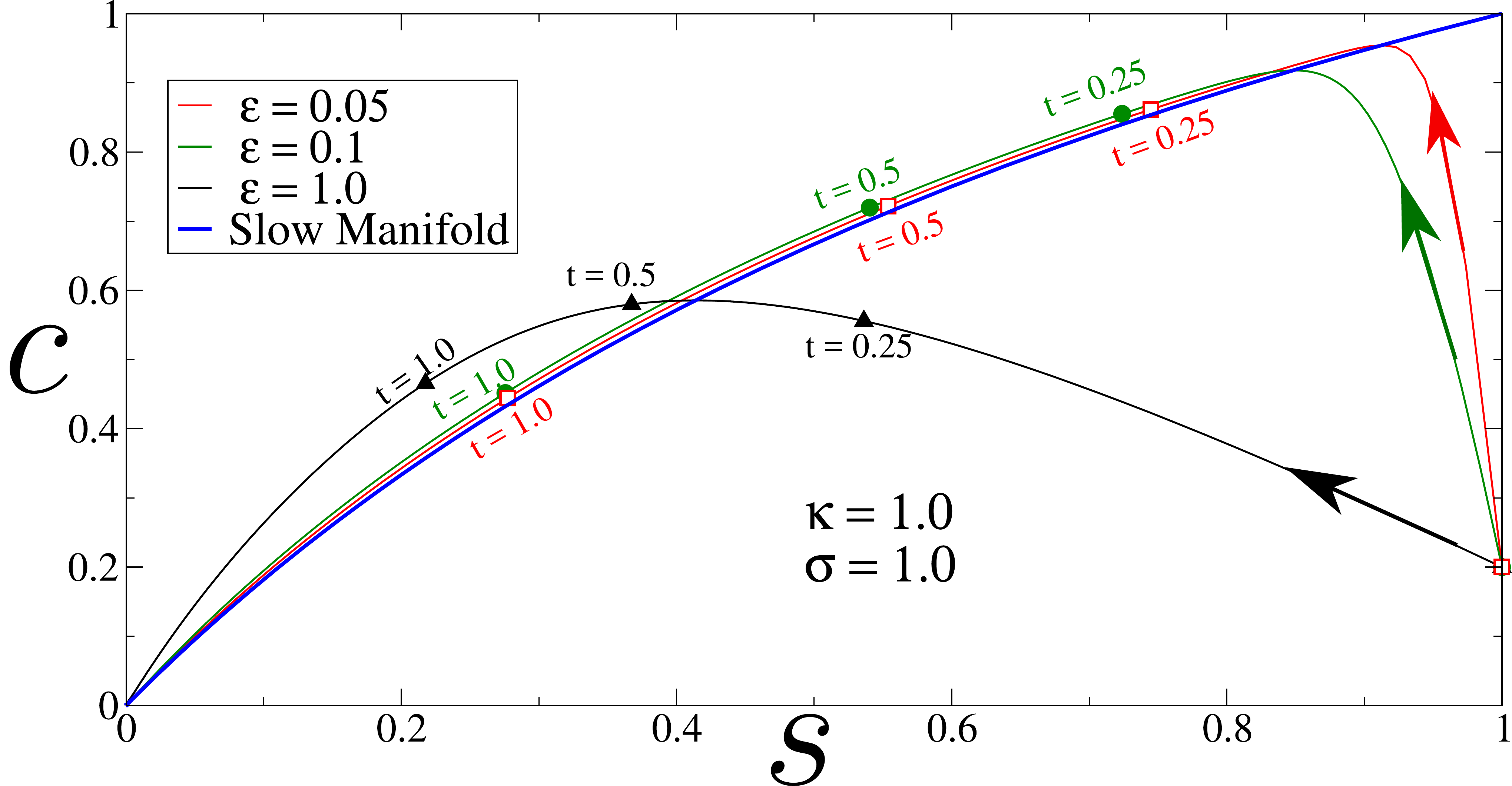}&
\includegraphics[width=0.45\textwidth]{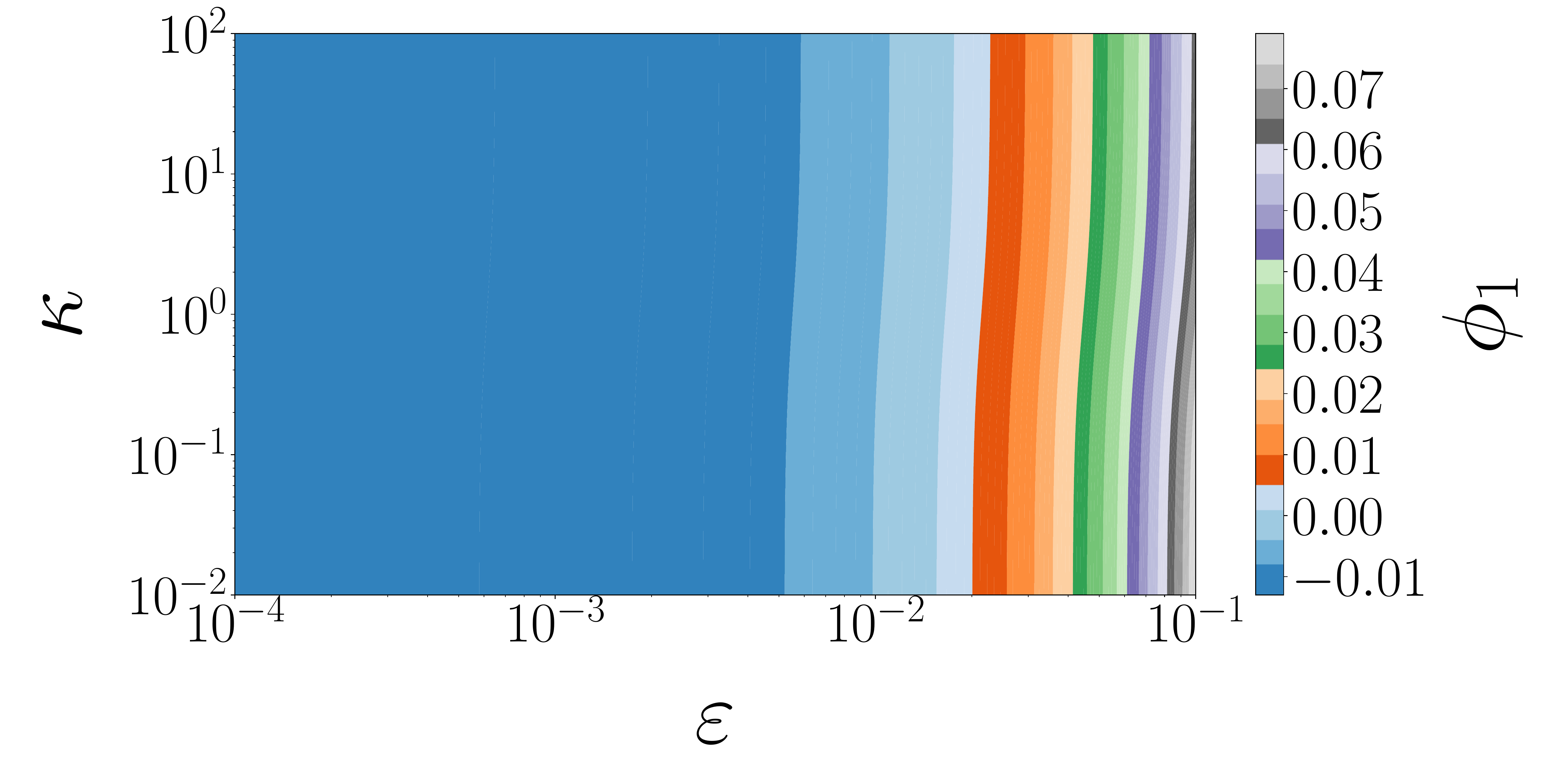}\\
(a) & (b) \\
\includegraphics[width=0.45\textwidth]{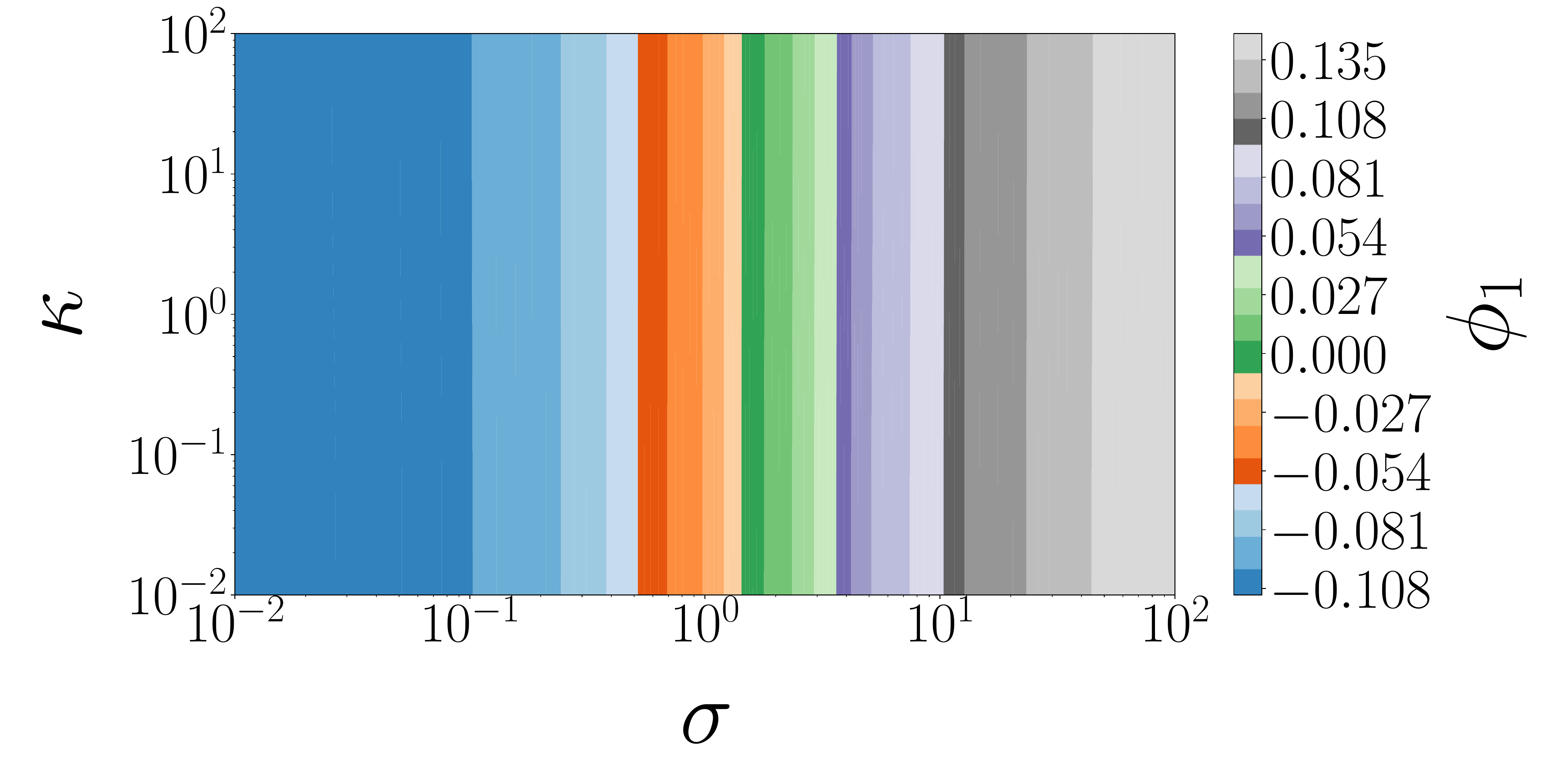}&
\includegraphics[width=0.45\textwidth]{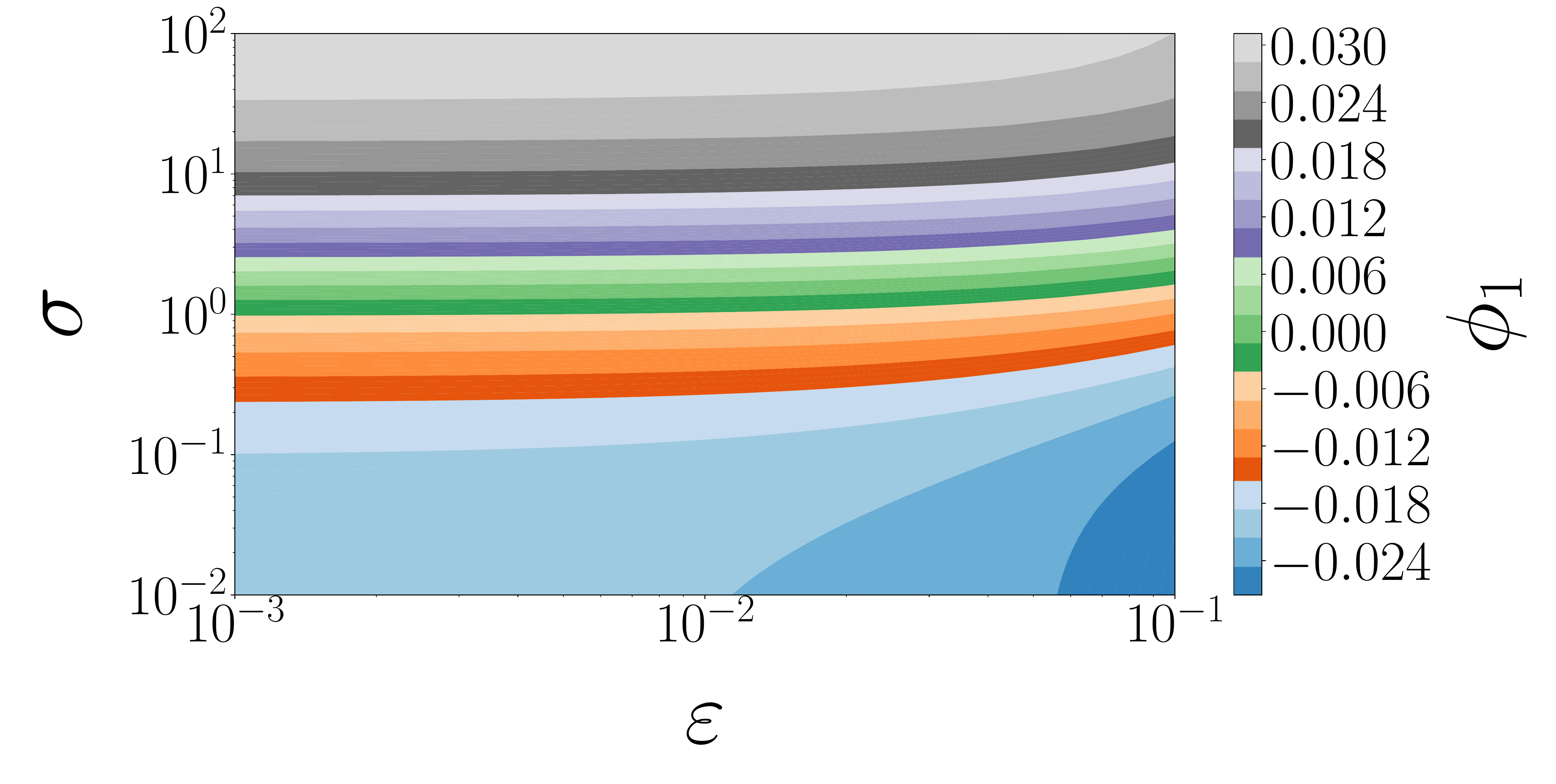}\\
(c) & (d)
\end{tabular}
\caption{
(a) Phase portrait of the rescaled MMH model~\eqref{mmh-ODE} with $\sigma=\kappa=1$. The plotted trajectories start at $(1,0.2)$ and correspond to various $\eps-$values, to illustrate the rate of attraction to the slow manifold (blue) and subsequent convergence to the origin. (b--d) Parameterization of the facets $\sigma=1$, $\eps=0.1$ and $\kappa=10$ of the $(\eps,\sigma,\kappa)-$space (input space) by the leading (output-only) DMAPS eigenvector $\phi_1$. Evidently, $\kappa$ is sloppy: the output is insensitive to it over several orders of magnitude. As $\eps\downarrow0$, the system enters an asymptotic regime whose slow, reduced dynamics is strongly informed solely by $\sigma$.
}
\label{fig:mmh-sloppy}
\end{figure}
Continuing the development of a data-driven framework to identify effective parameters, we now treat a benchmark for model reduction methods. The MMH system \cite{johnson_original_2011,michaelis_kinetik_1913} describes conversion of a substrate $\mathrm{S}$ into a product $\mathrm{P}$ through mediation of an enzyme $\mathrm{E}$ and formation of an intermediate complex $\mathrm{C}$,
\[
\mathrm{S + E}
\
\xrightleftharpoons[k_{-1}]{k_1}
\
\mathrm{C}
\
\xrightarrow{k_2}
\
\mathrm{P + E} .
\label{s2c2p}
\]
Under conditions often encountered in practice, the first reaction step reaches quickly an (approximate) chemical equilibrium and becomes rate-limiting. Product sequestration proceeds on a much slower timescale, during which the first reaction approximately maintains its quasi-steady state.

In that regime, simultaneous state {\em and} parameter space reduction is possible, as system evolution is described by a single ODE involving a subset of the problem parameters. There have been several, increasingly elaborate estimates of the parametric regime where QSSA applies, which were underpinned by different system nondimensionalizations. The first key estimate was that of \cite{heineken_mathematical_1967}, where the authors identified that regime as $E_T \ll S_T$ involving the (conserved) total amounts of enzyme, $E_T = E + C$, and substrate, $S_T = S + C + P$. In that regime, nearly all enzyme molecules become quickly bound to substrate and the complex saturates. The authors of \cite{segel_quasi-steady-state_1989} brought the kinetic constants into play and extended the regime to $E_T \ll S_T + K_M$, where $K_M = (k_{-1}+k_2)/k_1$ is the so-called Michaelis--Menten constant. This asymptotic regime extends the one of \cite{heineken_mathematical_1967} by including the case where the complex dissociates much faster than it forms.

Our goal in this section is twofold: first, to identify the effective parameter(s) informing system evolution in the asymptotic regime; and second, to show how the extended parametric region of \cite{segel_quasi-steady-state_1989} is captured in a data-driven manner by our methodology. To accomplish this in a completely automated way would necessitate using a black-box simulator for (a subset of) the dimensional state variables $S,E,C,P$ evolving in dimensional time $T$. This, in turn, would necessitate a candid discussion on tuning of monitoring times to capture the slow dynamics and how that relates to experimental/simulation data. We circumvent this issue here for brevity and focus, instead, on the equivalent, non-dimensional version in \cite{segel_quasi-steady-state_1989}. In that version, $T,S,C,E,P$ have been rescaled into dimensionless variables $t,s,e,c,p$; additionally, $e,p$ have been eliminated using the enzyme and substrate conservation laws. The result is the $2-$D ODE system
\begin{equation}
\begin{array}{rclcl}
\dot{s}
&=&
(\kappa+1)
\left[
- (1 + \sigma) s + \sigma c s + \kappa(\kappa+1)^{-1} c
\right] ,
\vspace{2mm}\\
\eps \dot{c}
&=&
(\kappa+1)
\left[
\ \ \, (1 + \sigma) s - \sigma c s - c
\right] .
\end{array}\!\!
\label{mmh-ODE}
\end{equation}
%
%
\begin{figure}[H]
\centering
\begin{tabular}{cc}
\includegraphics[width=0.45\textwidth]{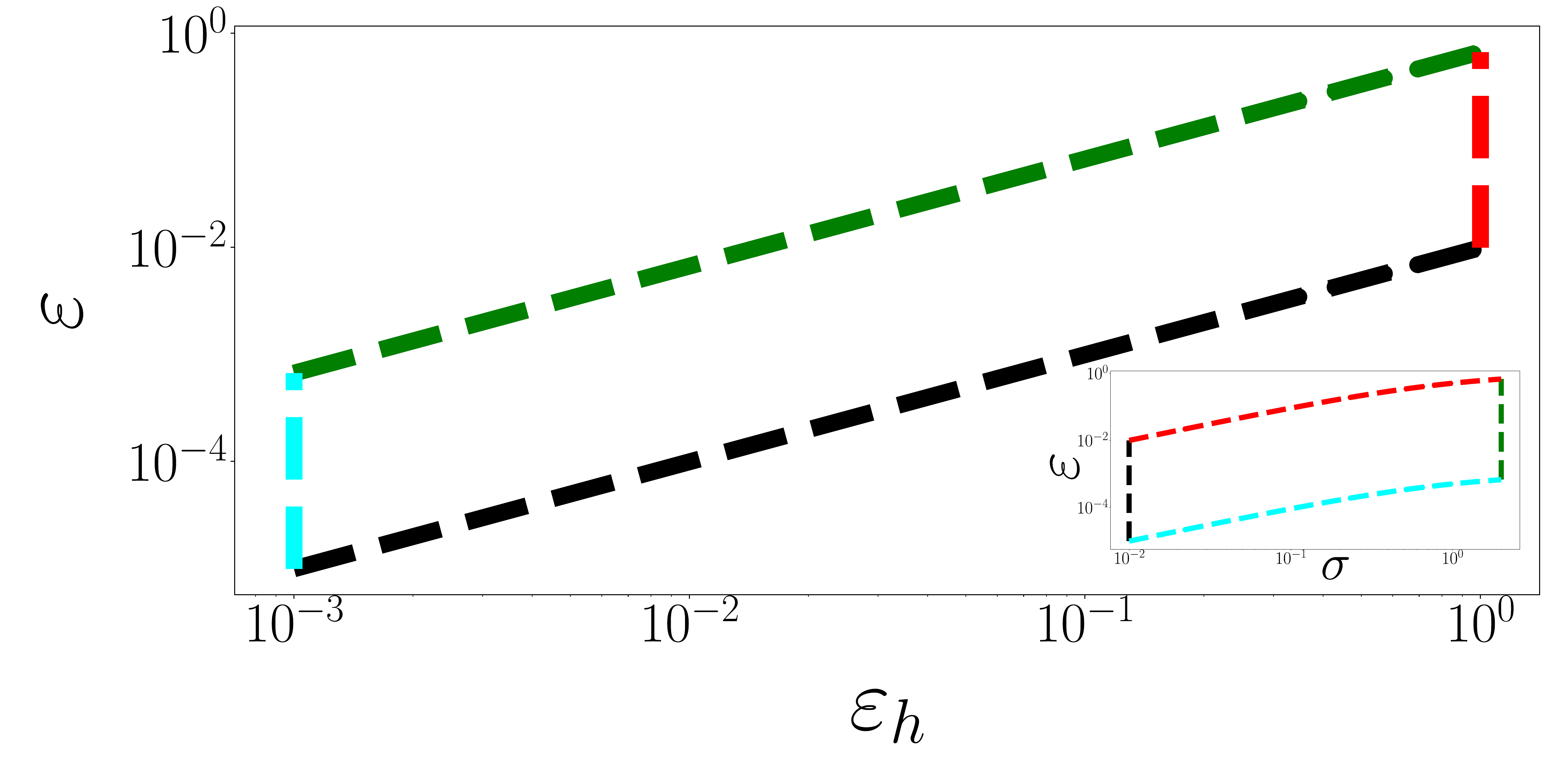} &
\includegraphics[width=0.45\textwidth]{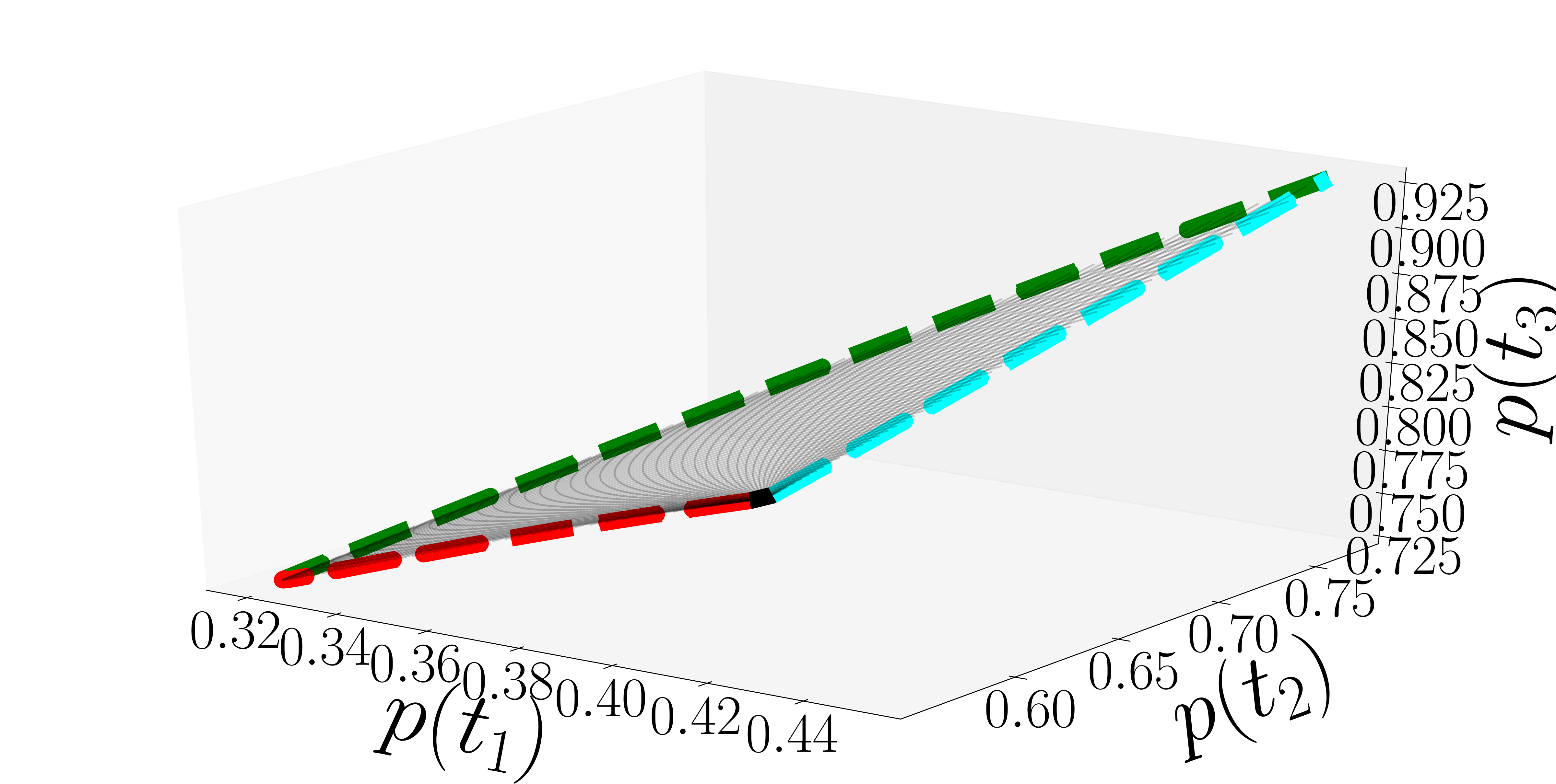}\\
(a) & (b)\\
\includegraphics[width=0.45\textwidth]{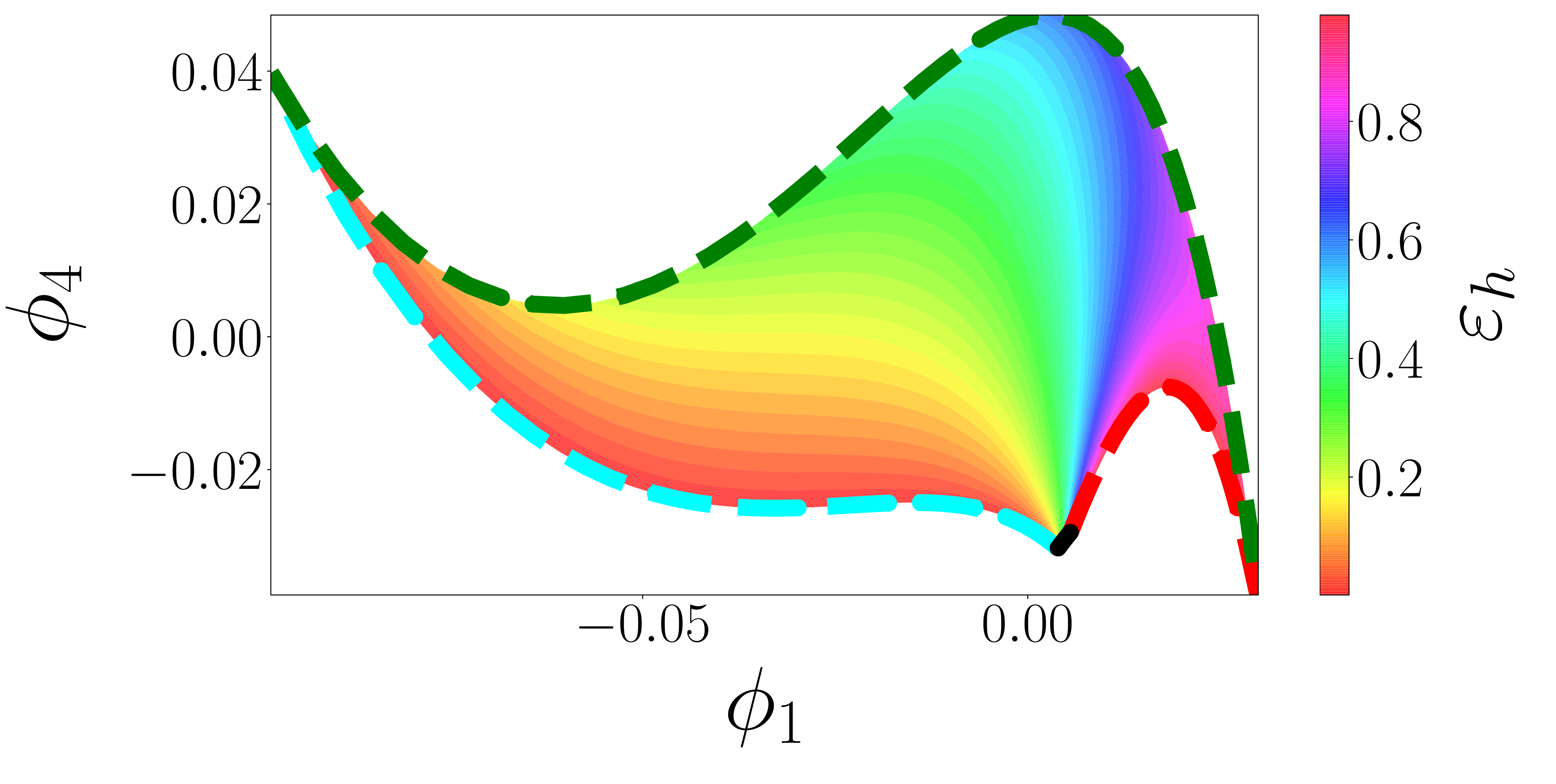}&
\includegraphics[width=0.45\textwidth]{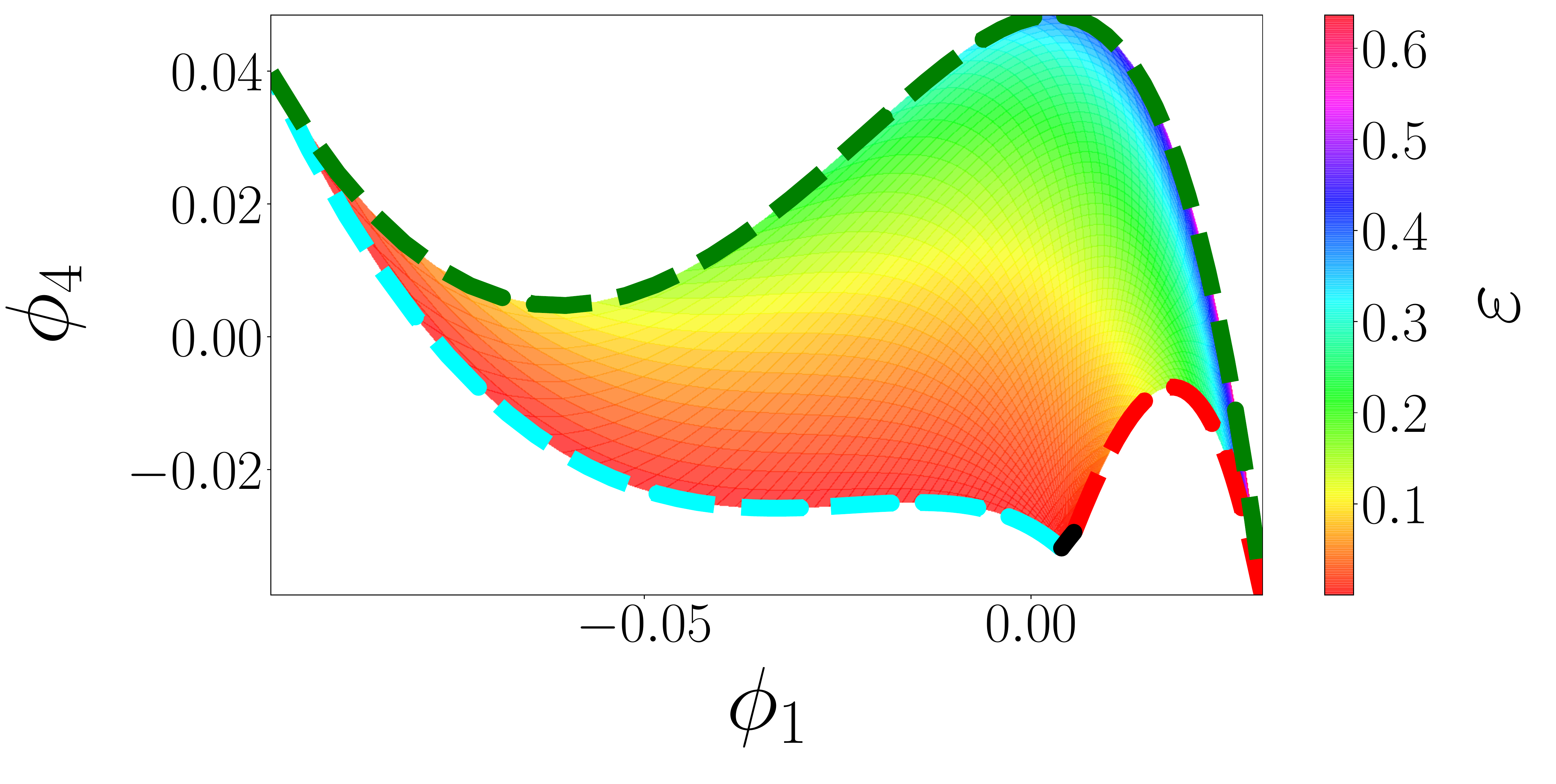}\\
(c) & (d)\\
\includegraphics[width=0.45\textwidth]{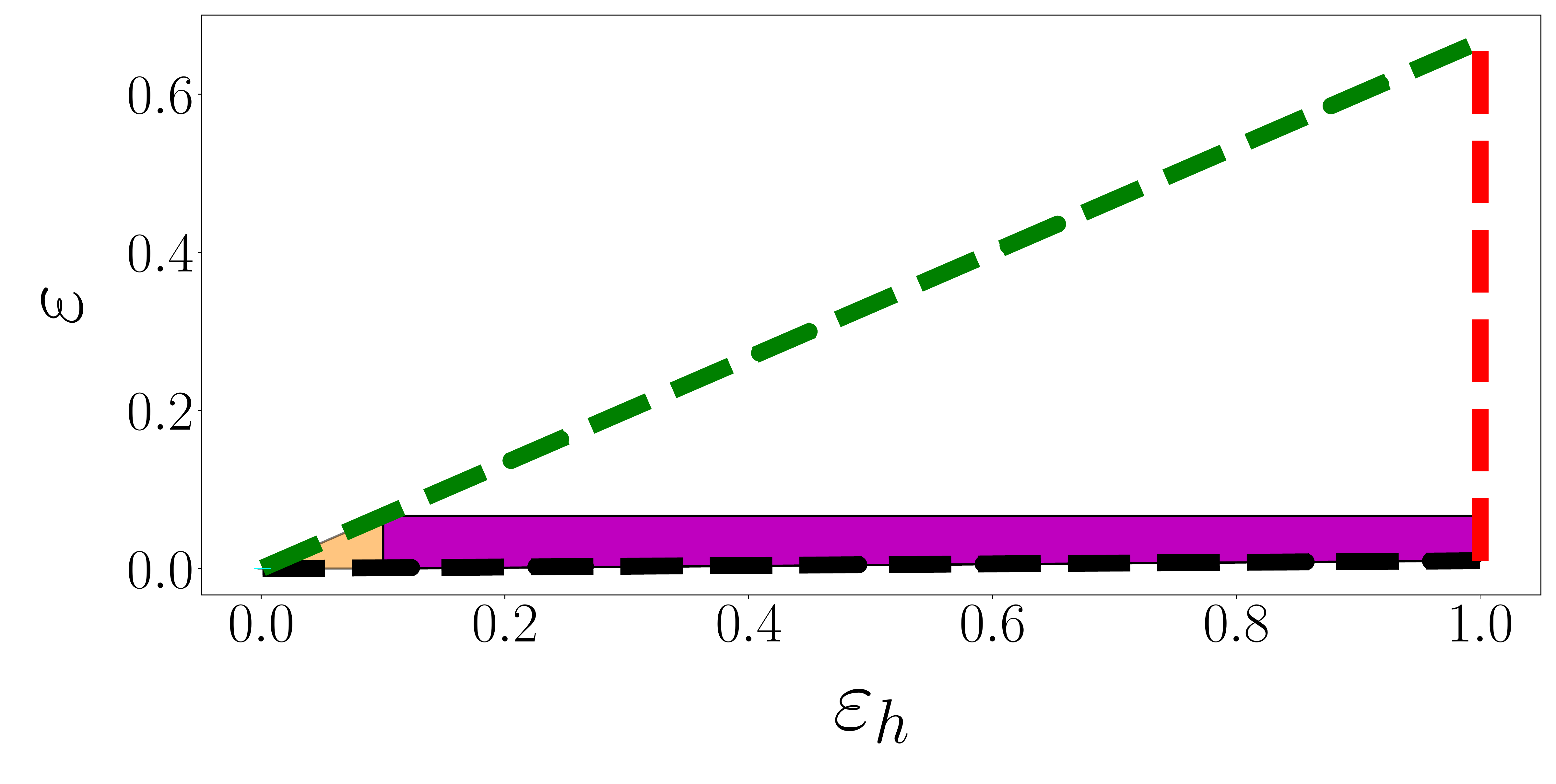}&
\includegraphics[width=0.45\textwidth]{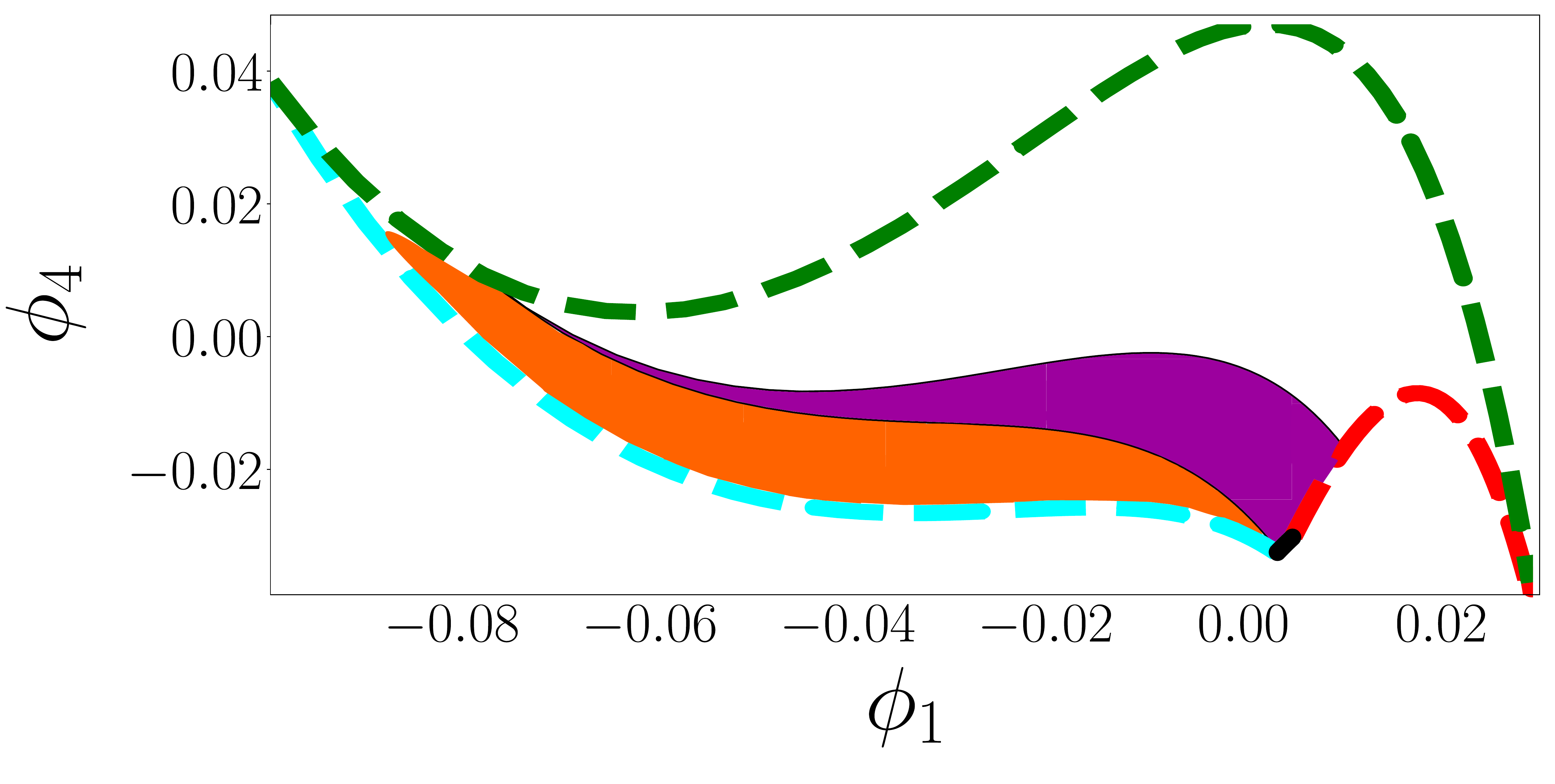}\\
(e) & (f)
\end{tabular}
\caption{
(a) Relevant parameter domains in $(\eps_h,\eps)-$space and in $(\sigma,\eps)-$space (inset; both logarithmic), related through $\sigma=1/(\eps_h/\eps-1)$. Boundaries are colored consistently across panels to help visualize the transformations between spaces.
(b) Model manifold $\mathcal{M}$ in output space, as $(\eps_h,\eps)$ vary and $\kappa=10$ is fixed. (c$-$d) Model manifold observed in DMAPS space, colored by $\eps_h$ in panel~(c) and by $\eps$ in panel~(d). (e) Similar to (a) but with uniform (not logarithmic) spacing, and with the highlighted regions $\eps_h\ll1$ (orange triangle) and $\eps\ll1$ (orange triangle and purple rectangle). (f) Image of the regions highlighted in panel~(e) in DMAPS space. 
\label{fig:mmh-maps} 
}
\end{figure}
The composite parameters here are $\eps = E_T/(S_T + K_M)$, $\sigma=S_T/K_M$ and $\kappa=k_{-1}/k_2$, and they may in principle assume any positive value. The initial conditions $(s_0,c_0)$ are arbitrary, and $e,p$ can be recovered from the rescaled conservation laws $e + \sigma c/(\sigma+1) = 1$ and $s + \eps c + p = 1$. The system is expressed in \emph{slow time}, so that quasi-steady state is achieved over an $\Ord(\eps)$ time and product sequestration occurs over an $\Ord(1)$ timescale. In this reformulation, the asymptotic regime where QSSA applies is $\eps \ll 1$, according to \cite{segel_quasi-steady-state_1989}, and $\eps_h = (1 + 1/\sigma)\,\eps \ll 1$ according to \cite{heineken_mathematical_1967}; the former plainly extends the latter. Initially, we select as our observable the rescaled complex concentration at distinct times $(t_1,t_2,t_3)=(0.5,1,1.5)$, so $\omr(\prm) = [c(t_1\vert\prm),c(t_2\vert\prm),c(t_3\vert\prm)]^\mathrm{T}$. Our parameter set is the triplet $\prm=(\eps,\sigma,\kappa)$, with $(s_0,c_0)=(1,0)$ fixed as in the original experimental setting \cite{michaelis_kinetik_1913}.

Figure~\ref{fig:mmh-sloppy}(b$-$d) demonstrates that the model response is unaffected by $\kappa$ and strongly affected by $\sigma$, with the limit $\eps \downarrow 0$ corresponding to an asymptotic regime. Further, Fig.~\ref{fig:mmh-sloppy}(d) makes it plain that the system evolution in that regime is controlled by $\sigma$. This is in stark contrast to the parameter-free reduced dynamics of caricature~\eqref{model-singpert} and agrees with theory, which predicts that the evolution of $p$ in $\Ord(1)$ timescales is dictated by the leading order problem \cite{segel_quasi-steady-state_1989} (see also SI).
\begin{equation}
p = 1 - s - \eps c ,\ \mbox{subject to} \quad \dot{s}= -c
= -(1 + 1/\sigma)\frac{s}{s + 1/\sigma}.
\label{MMH-ODE-red}
\end{equation}
On the basis of these results, we conclude that the model manifold is effectively $2-$D {\em and not $3-$D} as one might initially surmise, with the asymptotic regime $\eps \downarrow 0$ corresponding to a curve parameterized only by $\sigma$. As a corollary, the model manifold dimensionality transitions from two to one in that regime, without being further reduced to zero. This is evident in Fig.~\ref{fig:mmh-maps}(b), showing (part of) the model manifold for the setup above.

We next turn to a data-driven characterization of the asymptotic regime and relate that to the characterizations in \cite{heineken_mathematical_1967,segel_quasi-steady-state_1989}.
Using simulated trajectories of \eqref{mmh-ODE} and applying our DMAP methodology with an output-only informed metric, we coordinatize the model manifold through the independent eigenmodes $(\phi_1,\phi_4)$. Figures~\ref{fig:mmh-maps}(c--d) show that manifold in DMAPS space; the asymptotic limit is the lower-left bounding curve (light blue). We can use these diffusion coordinates to characterize the asymptotic regime as a \emph{neighborhood of that boundary}, so that the success of $\eps_h$ and $\eps$ in capturing that regime is measured by the extent their level sets track the boundary. Figures~\ref{fig:mmh-maps}(c--d) color the DMAPS domain by $\eps_h$ and by $\eps$; plainly, the $\eps-$coloring traces the domain boundary quite well, with $\eps \ll 1$ represents a {\em bona fide} neighborhood of it. Small values of $\eps_h$, on the other hand, fail to outline such a neighborhood: all level sets coalesce at the single point representing the $\eps_h-$axis (i.e. the regime $\sigma \downarrow 0$). This is made even plainer in Fig.~\ref{fig:mmh-maps}(e--f), where one sees how the $\eps_h \ll 1$ regime misses a substantial part (colored purple) of the asymptotic regime captured by $\eps \ll 1$. On account of this, we can conclude that $\eps$ is indeed a better ``small parameter'' than $\eps_h$. It is important to note that a black-box, data-driven approach can have no knowledge of $\eps$, $\eps_h$ or any other ``human'' description of the problem. What it can do, as we just saw, is enable us to {\em test human-generated hypotheses} on the data; we -- or Segel and Slemrod \cite{segel_quasi-steady-state_1989} -- are the ones generating the hypotheses.

\section{Non-invertible input-output relations}

Throughout this paper so far, we have used an output-only-informed kernel to obtain intrinsic DMAPS parameterizations of the combined input--output manifold. Our approach consisted of using eigenfunctions of the Laplace--Beltrami operator on the {\em model manifold}, and our insights about parameter (input) space came from 
how it was jointly parametrized by these eigenfunctions. The approach was useful in the data-driven study of parameter non-identifiability and even sloppiness.
We will now show that it fails dramatically when the mapping from parameter space to the model manifold is noninvertible,
i.e. when distinct, {\em isolated} parameter values produce identical model responses, $\omr(\mathbf{p}) = \omr(\mathbf{p}')$ for $\mathbf{p} \ne \mathbf{p}'$.
\begin{figure}[!h]
\centering
\includegraphics[width=0.9\textwidth]{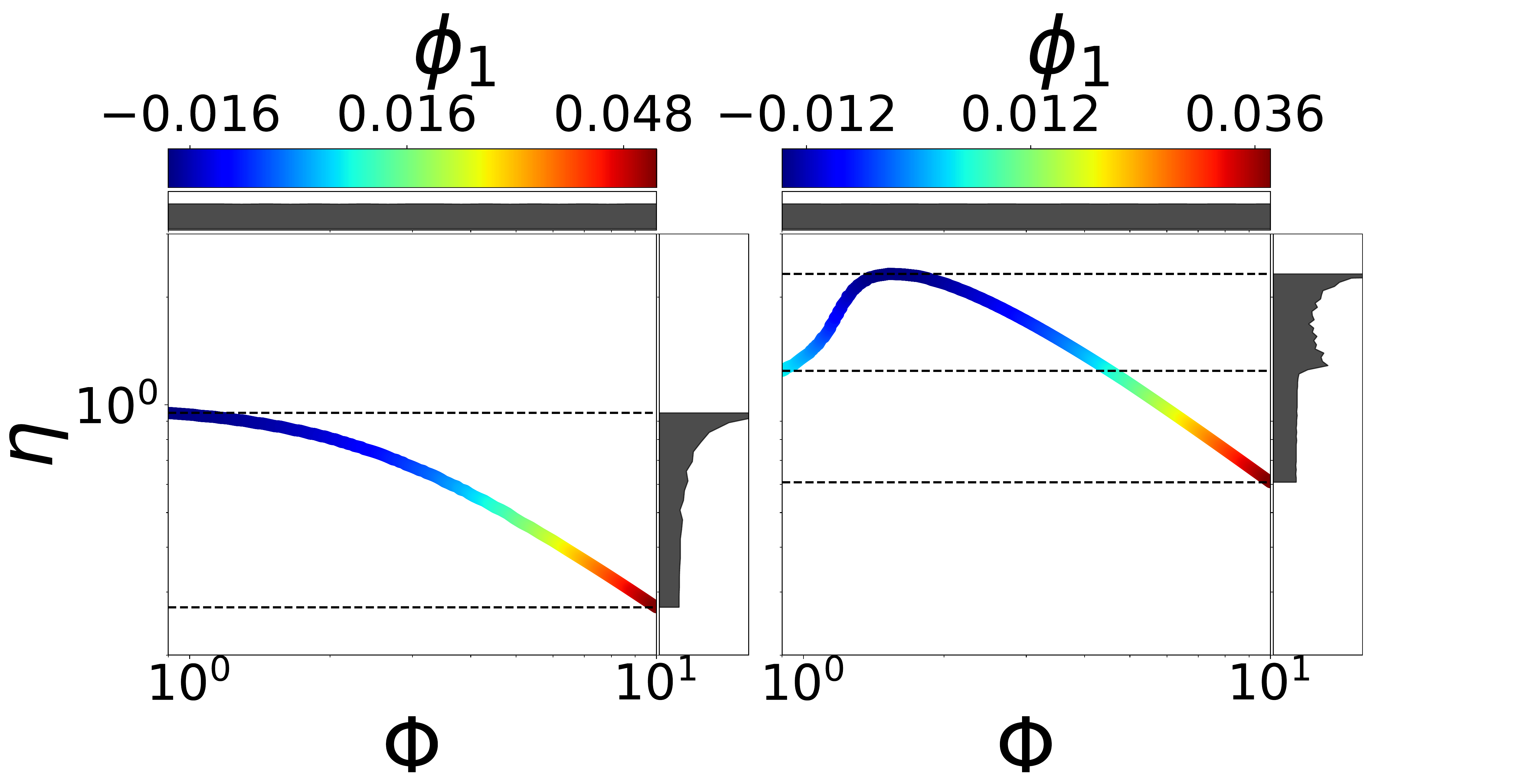}
\caption{
Left: Output $\eta$, vs. input $\Phi$ for the isothermal catalyst pellet ($\beta = 0.0$, see 
SI) colored by $\phi_1$. Right: Same plot for a nonisothermal ($\beta = 0.2$) pellet.
For uniform input $\Phi$ sampling, the observed output pdf's (see text) 
are plotted alongside each panel. Using output-only informed DMAPS is unable, 
as the coloring shows, to accurately parameterize the noninvertible case:
Widely different sections of the curve on the right take on the same color. 
Note also the discontinuity in the density along the $\eta$ axis in the right figure, a hallmark of noninvertibility.}
\label{fig:ni-eta}
\end{figure}

A well-known instance of this situation arises in the study of reaction--diffusion in porous catalysts and is illustrated in Fig.~\ref{fig:ni-eta}. For isothermal reactions, the output -- the dimensionless ``effectiveness factor'' $\omr(\prm) \equiv \eta$ -- is a monotonic function (with known asymptotic limits) of the input -- the Thiele modulus $\prm \equiv \Phi$ \cite{weisz_behaviour_1995}
(Fig.~\ref{fig:ni-eta}, left). For {\em exothermic reactions}, however, $\eta$ may depend on $\Phi$ nonmonotonically and the relation becomes noninvertible; alternatively, points on the model manifold are revisited, as the input sweeps the positive real axis, Fig.~\ref{fig:ni-eta} (right). Sampling the input 
$\Phi$ uniformly on the horizontal axis naturally results in a nonuniform density for the output $\eta$ on the vertical axis 
(plotted on the right of each panel in Fig.~\ref{fig:ni-eta}). This observed output probability density function (pdf)  embodies the input-output relation and brings to mind an analogy
with Bayesian measure transport from a prior density to a posterior one. It is worth noting that, noninvertibility causes pronounced discontinuities on the output pdf on the right.

\begin{figure}[!h]
\centering
\begin{tabular}{cc}
\includegraphics[width=0.45\textwidth]{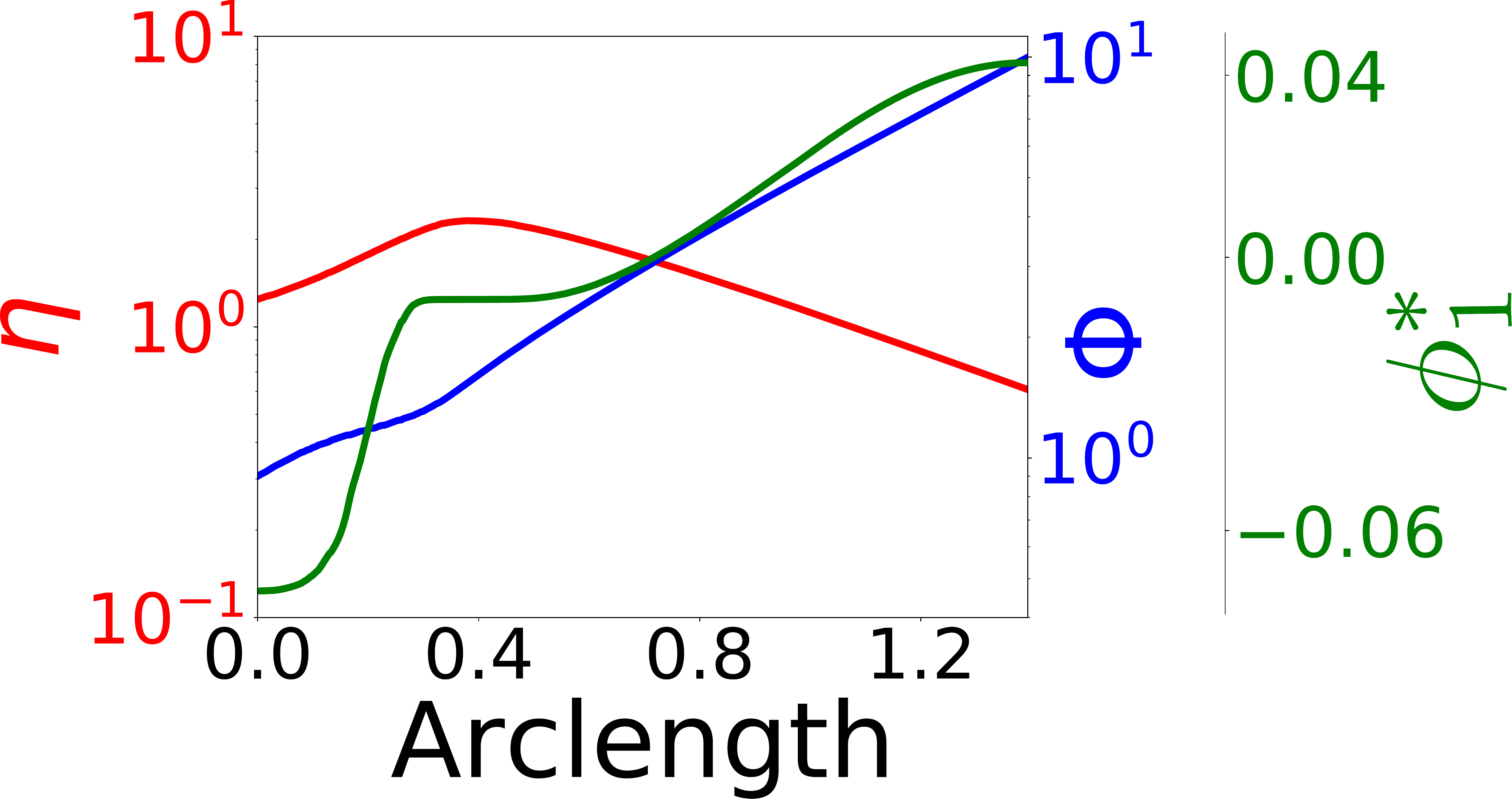}&
\includegraphics[width=0.43\textwidth]{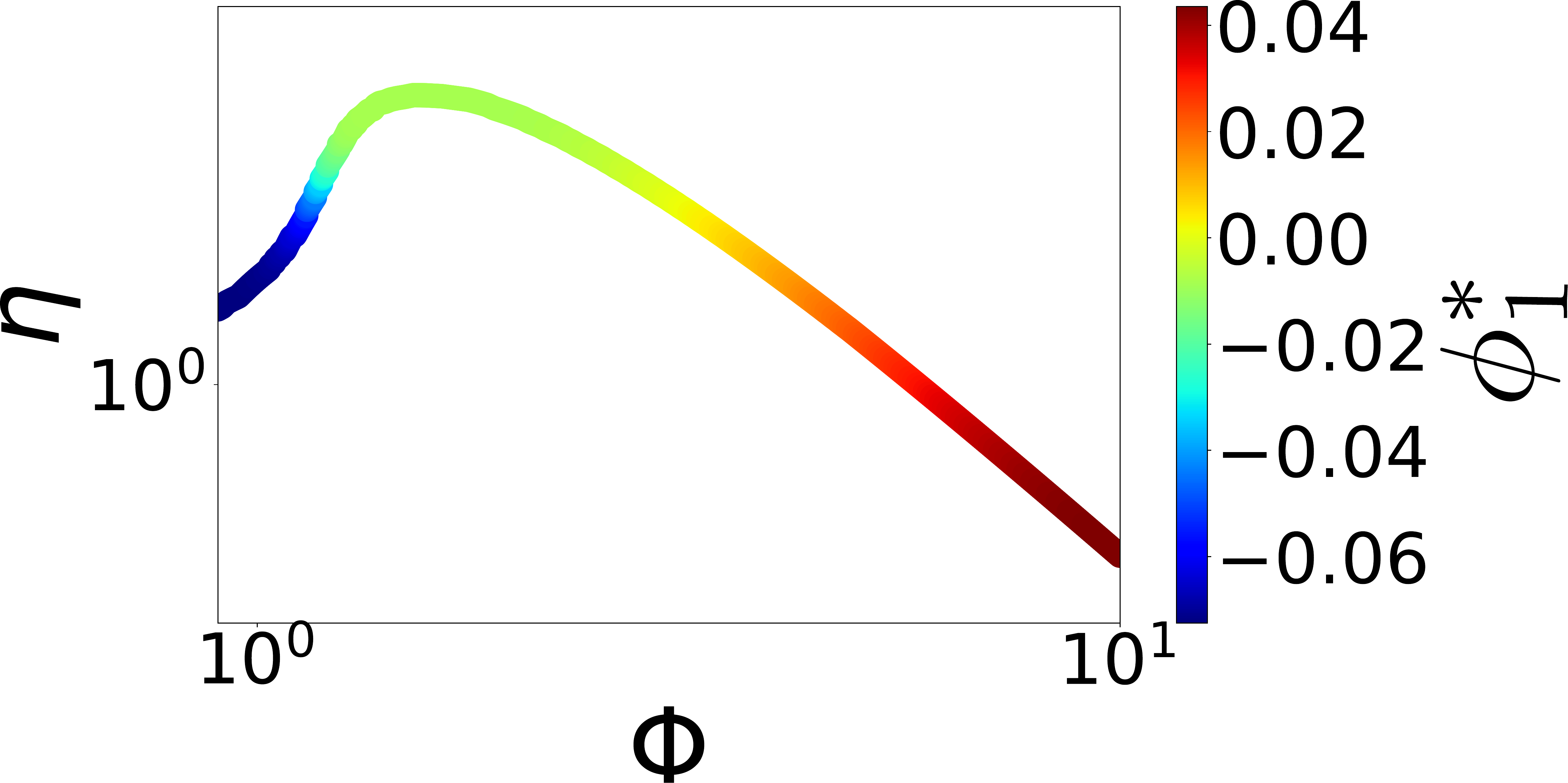}\\
(a) & (b) \\
\includegraphics[width=0.45\textwidth]{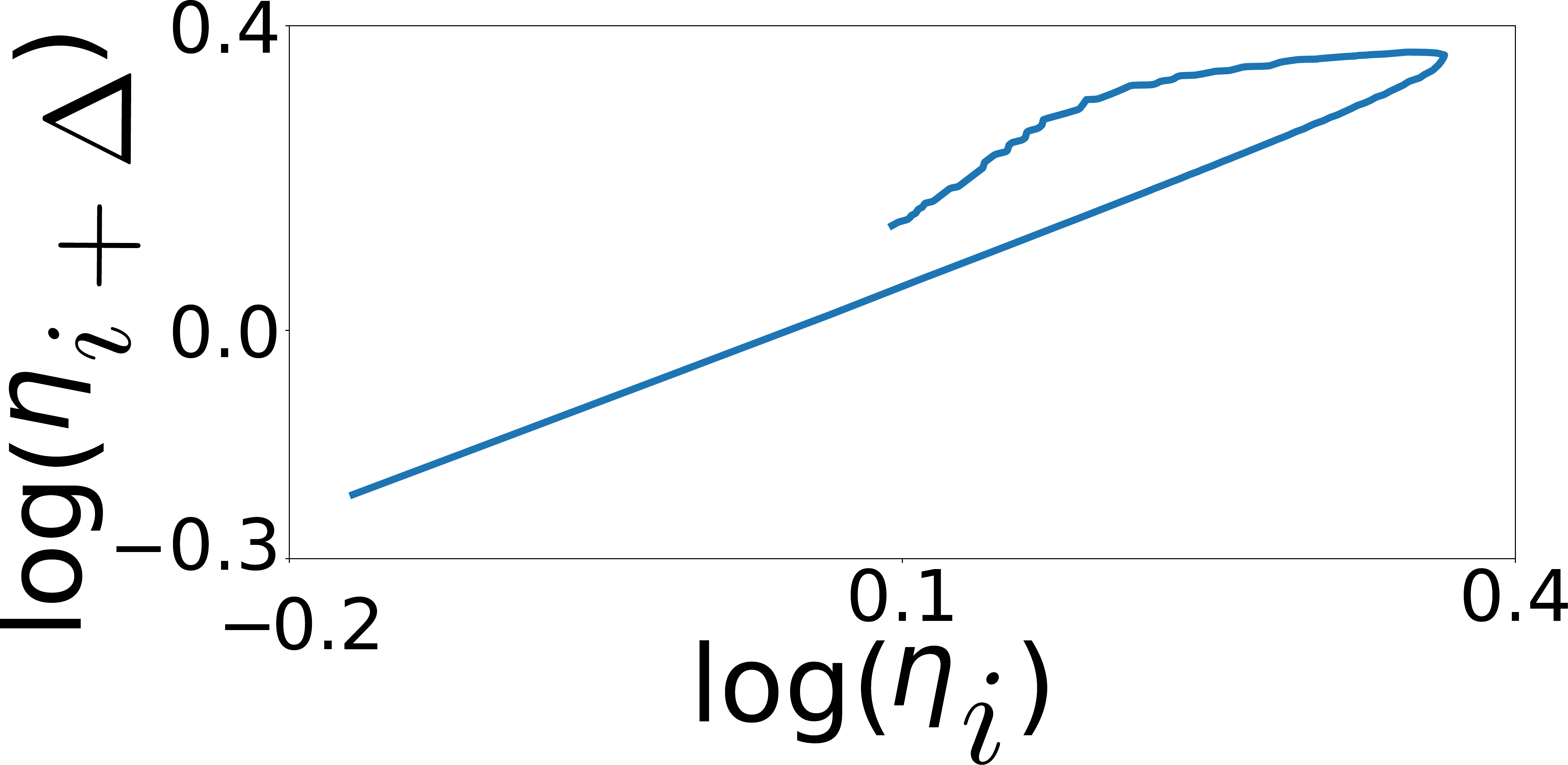} &
\includegraphics[width=0.43\textwidth]{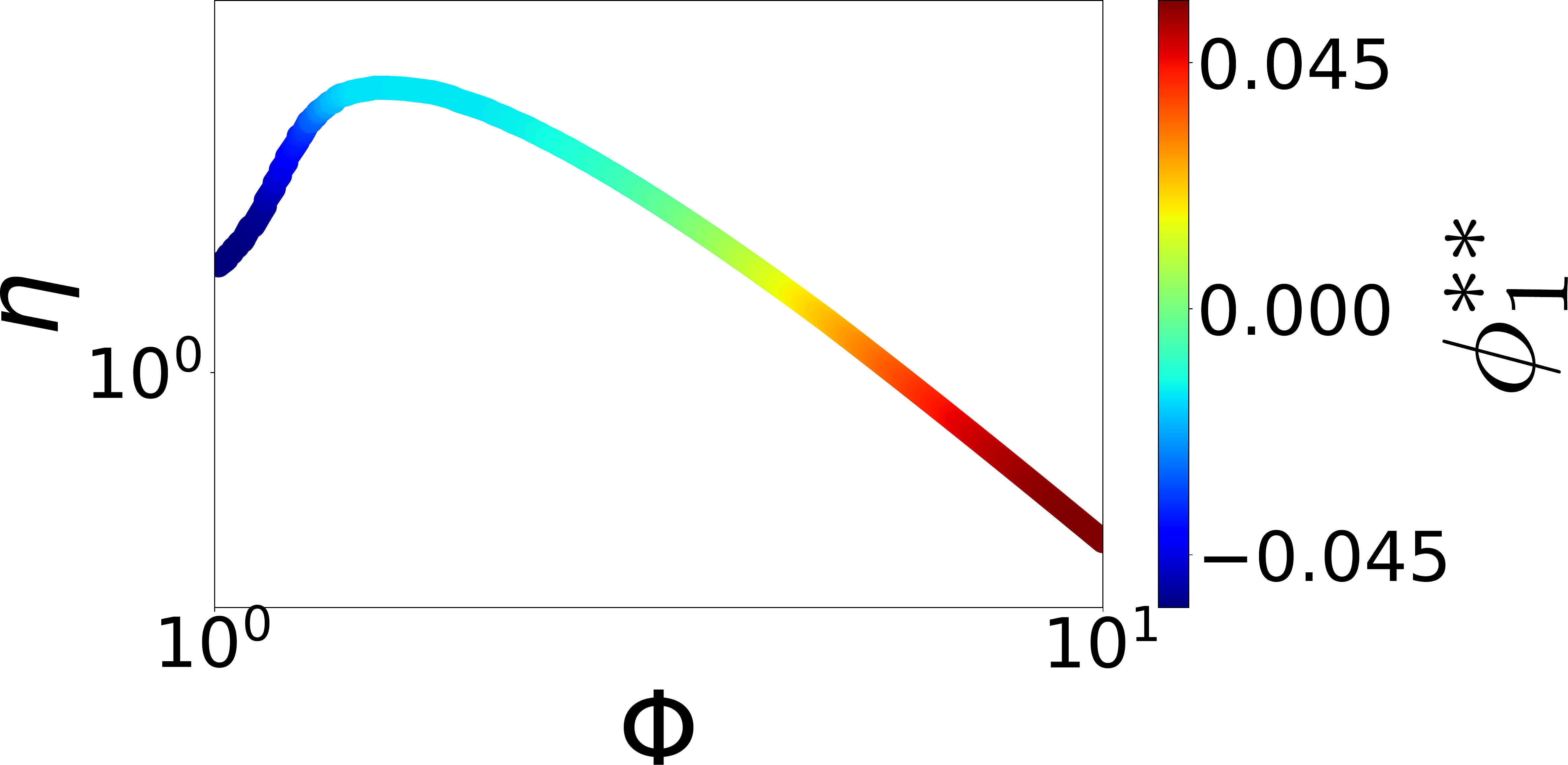}  \\
(c) & (d)
\end{tabular}
\caption{
(a) Input $\Phi$, output $\eta$ and DMAPS eigenfunction $\phi_1^*$ 
corresponding to \eqref{noninv:kernel}
plotted against the arclength of the curve in Fig.~\ref{fig:ni-eta}.
The eigenfunction clearly parameterizes {\em both} input and output. 
(b) Input--output response curve colored by the eigenfunction $\phi_1^*$ using \eqref{noninv:kernel} with $\epsilon = 0.0125$.  In both (b) and (c), $\Phi$ was sampled on a uniform grid between $0.9$ and $10$ for a total of $1043$ points.
(c) Plot of $\log(\eta_{i+\Delta})$ (corresponding to $\log(\Phi_i)+\Delta$ where $\Delta \approx 0.05$.)
against $\log(\eta_i)$ (corresponding to $\log(\Phi_i)$); this ``delay
embedding'' is one-to-one with the original curve. 
(d) Input--output 
response curve colored by the eigenfunction $\phi_1^{**}$ corresponding to the ``augmented 
output'' DMAPS kernel 
$K^{**} \equiv \exp\left( -(\bar{\eta} - \bar{\eta}')^2/\epsilon^2 \right)$, 
where 
$\bar{\eta}=(\eta_i, \eta_{i+\Delta})$ 
and similarly for $\bar{\eta}'$.  Here we used $\epsilon = 0.01$.
Either of these new, modified kernels parameterizes the non-invertible response curve successfully.
} 
\label{fig:ni-combined}
\end{figure}

Coloring input--output ($\eta - \Phi$) profiles by the leading DMAPS eigenfunction of an output-only-informed kernel shows that the data-driven coordinate, $\phi_1$, which successfully recovered (parameterized) the input $\Phi$ on the left {\em fails to do so on the right}. The problem lies with the output-only metric employed, and its resolution requires a new, more informative DMAPS kernel such as
\begin{align}
K^*(\prm, \prm') = \exp \left( - \frac{\|\prm - 
\prm'\|^2}{\epsilon^2} - \frac{\| \omr(\prm) - 
\omr(\prm') \|^2}{\epsilon^a} \right) .
\label{noninv:kernel}
\end{align}
Taking into account both inputs ($\prm - \prm'$) and outputs ($\omr(\prm) - \omr(\prm')$), this kernel manages to differentiate inputs having the same output.
Figure~\ref{fig:ni-combined}(b) corroborates the appropriateness of this kernel for $a=4$:  its primary eigenvector varies monotonically over the model manifold. 
This is also evident in Fig.~\ref{fig:ni-combined}(a), in which we have plotted input, output and the data-driven parameter $\phi^{*}_1$ 
against arclength of the input--output response curve. 
In effect, $\phi^{*}_1$ is in an one-to-one correspondence with the arclength,
and thus ``discovers'' a good parameterization of the curve. 
This particular ($a=4$) kernel -- originally proposed by Lafon \cite{lafon2004diffusion} in a different context -- prioritizes output over input;
due to the $\epsilon-$scalings, the input-term only becomes significant when needed, i.e. for {\em nearby} inputs producing similar outputs.

The use of appropriately scaled input {\em and} output similarities can thus resolve input--output noninvertibility. 
Can such noninvertibility be resolved when we do not know the inputs, yet have
some control over the measurement process? The answer is, remarkably, in the affirmative.
A data-driven parameterization of input space can 
be obtained {\em even in the absence of actual recorded input measurements} 
by using a little local history of 
output measurements in the spirit of the Whitney, Nash and Takens embedding
theorems \cite{nash1956imbedding,takens_detecting_1981,whitney_self-intersections_1994}. 
Figure~\ref{fig:ni-combined}(c--d) illustrates how unmeasured inputs can, 
in a sense, be recovered by recording {\em pairs of output measurements} rather than single output measurements.
Specifically, we formulate a measurement protocol in which the output $\eta$ is measured sequentially, first for a random input $\Phi$ 
and then for $\Phi = \Phi + \Delta$ 
(for some unknown but fixed $\Delta$). 
Using this analogy to Takens delay embeddings in nonlinear dynamics, 
redefining the model manifold in terms of such measurement pairs, 
and reverting to the output-only-informed metric based on such pairs 
yields a single data-driven effective parameter $\phi^{**}_1$ which consistently parametrizes both the (unkown) input $\Phi$ as well as the output $\eta$ pairs. 
Using a little measurement history can thus also resolve model 
noninvertibilities, and allow us to parametrize input-output {\em relations}.

\section*{Discussion}

We presented and illustrated a data-driven approach to effective parameter identification in dynamic ``sloppy'' models -- model descriptions containing more parameters than minimally required to describe their output variability. Our manifold-learning tool of choice was Diffusion Maps (DMAPS), and we applied it to datasets that typically consisted of input--output combinations generated by dynamical systems. The inputs were mostly model parameters, but we also viewed initial conditions as inputs to differentiate between (what traditionally would be referred to as) singularly and regularly perturbed multiscale models. The outputs were ensembles of temporal observations of (some of) the state variables. By modifying the customary DMAPS kernel to rely \emph{mainly} on -- or, in most of the paper, {\em only} on -- the observed outputs, we were able to ``sense'' the sloppy directions and automatically unravel nonlinear effective model re-parameterizations. 

It is important to note that, as often the case with numerical procedures, this approach does {\em not} characterize the effective parameters through explicit algebraic formulas. In fact, we saw in our treatment of the ABC model that an off--the--shelf, algebraically formulated effective parameter ($\keff^\mathrm{QSSA}$) predicted system output worse than the parameter found by DMAPS. This approach (a) helps test hypotheses about the number and physical interpretability of effective parameters, see our in-context discussion of the MMH model; (b) provides a natural context in which to make predictions for new inputs, through ``smart'' interpolation (matrix/manifold completion); and (c) assists experimental design through intelligent sampling of input space (see e.g. the biasing of computational experiments in \cite{ChiavazzoE5494,georgiou2017exploration}).
Clearly, what was achieved here by the sampling of ODE model outputs can in principle be extended to PDE models by sampling in time and space.
The leading eigenfunctions of our DMAPS-based approach (effective parameters) are, in general, {\em nonlinear} combinations of the system parameters. Actively changing the value of these combinations --``moving transversely to level sets'' of the eigenfunction -- leads to appreciable changes of the model output. 

\begin{figure}
\centering
\includegraphics[width=0.6\textwidth]{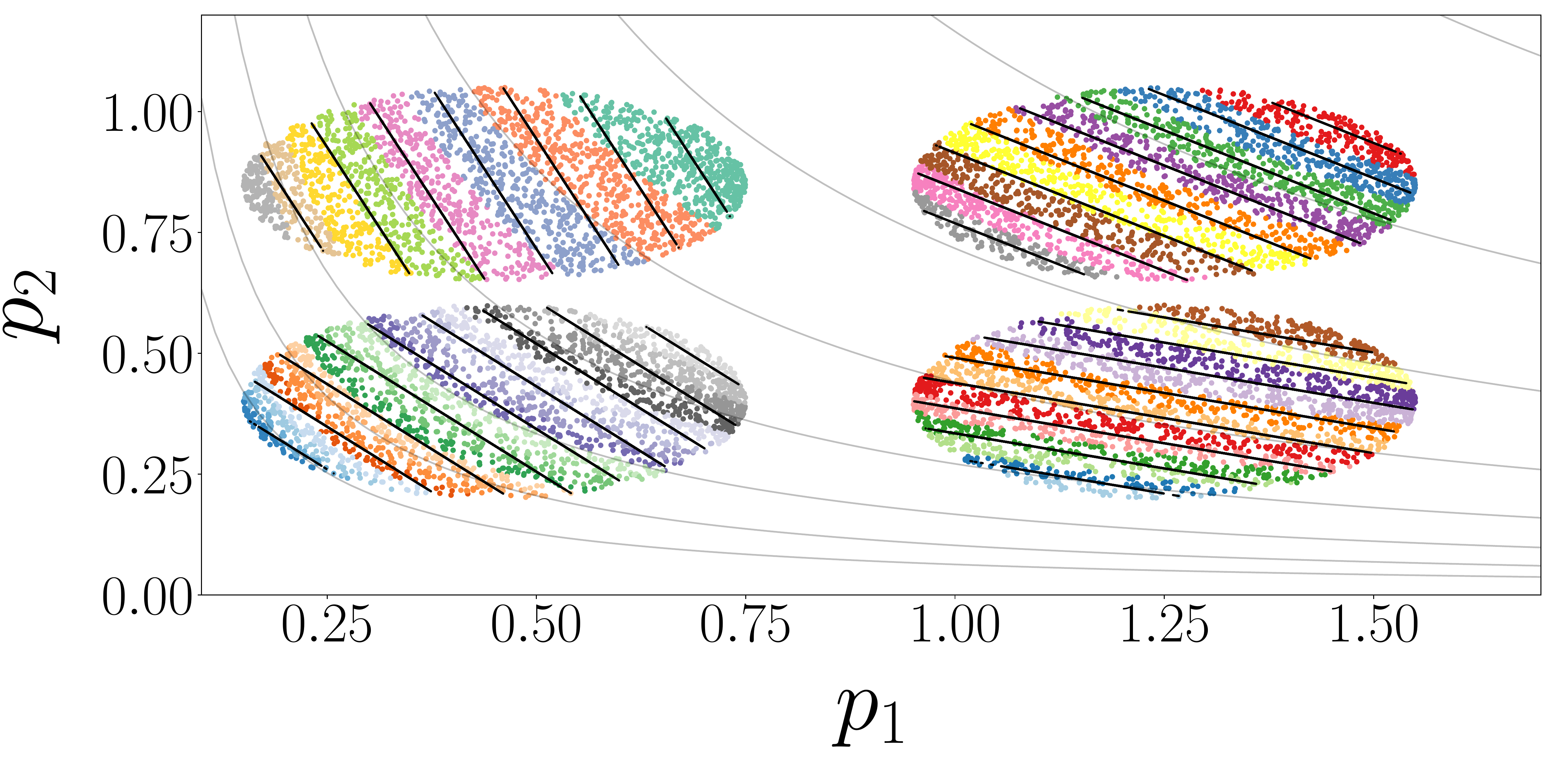}
\caption{
Comparison of Active Subspaces and DMAPS in local patches for a model similar to that shown in Fig.~\ref{fig:non-id}. See text and also SI. Here, $\omr(\textbf{p}) = \left( \log(p_1 p_2) \right)$.
}
\label{fig:active_sub_ovals}
\end{figure}

It is interesting to draw an analogy between identifying these effective parameters and the {\em linear} parameter combinations of Constantine and coworkers \cite{constantine_active_2014} affecting {\em scalar} model predictions: what they call ``active subspaces'' (see SI for a more detailed comparison). The analogy is illustrated here in Fig.~\ref{fig:active_sub_ovals}, for which we used our first, simple model that gave rise to Fig.~\ref{fig:non-id} but with the {\em scalar} output $\omr(p) = \log(p_1 p_2)$. The active subspace approach, applied independently to each of the datasets shown as oval patches, yields the solid black direction as ``neutral'' and their normal as the {\em active subspace} (per patch). To enable comparison, each dataset is also colored by the value of the leading DMAPS eigenfunction $\phi_1$ obtained with the output-only-informed metric. DMAPS plainly gives nonlinear ``neutral'' level sets (gray lines), with $\phi_1$ providing a nonlinear version of an ``active'' parameter combination: an effective parameter. Combining the data across patches leaves our curved level sets consistent; a linear approach would encounter problems, as these level sets start curving appreciably.
Two scenarios were discussed in this paper: the first, involving an output-only-informed kernel, proved useful in the data-driven study of sloppiness. Coordinates from the intrinsic model manifold geometry pulled back on the input (i.e. parameter) space provided our ``effective parameters''. The second, less explored scenario involved the non-invertible case where the same model output is observed for different isolated inputs and, more generally, one has input-output {\em relations}. The simple modifications of the DMAPS metric we used to resolve this, and the connection we drew to a ``measurement process history'' and embedding theorems, is a simple first research step in the data-driven elucidation of complex input--output relations by designing appropriate measurement protocols. We expect that similarity measures exploiting a measurement {\em process}, rather than a single measurement (e.g. ``Mahalanobis-like'' pairwise similarity measures \cite{Mahalanobis}) may well prove fruitful along these lines.
The physical interpretability of data-discovered effective parameters can be established in a postprocessing step, by testing whether they are one-to-one, on the data,  with subsets of equally many of the physical parameters.
\section{acknowledgments}
MK acknowledges funding by SNSF grant P2EZP2\_168833. AZ graciously acknowledges the hospitality of IAS/TUM, Princeton and Johns Hopkins. The work of IGK, JBR and AH was partially supported by the US National Science Foundation and by DARPA.
\section*{\refname}
\bibliographystyle{elsarticle-num-names}
\bibliography{manifold_learning}

\begin{thebibliography}{25}
\expandafter\ifx\csname natexlab\endcsname\relax\def\natexlab#1{#1}\fi
\providecommand{\url}[1]{\texttt{#1}}
\providecommand{\href}[2]{#2}
\providecommand{\path}[1]{#1}
\providecommand{\DOIprefix}{doi:}
\providecommand{\ArXivprefix}{arXiv:}
\providecommand{\URLprefix}{URL: }
\providecommand{\Pubmedprefix}{pmid:}
\providecommand{\doi}[1]{\href{http://dx.doi.org/#1}{\path{#1}}}
\providecommand{\Pubmed}[1]{\href{pmid:#1}{\path{#1}}}
\providecommand{\bibinfo}[2]{#2}
\ifx\xfnm\relax \def\xfnm[#1]{\unskip,\space#1}\fi
\bibitem[{Gutenkunst et~al.(2007)Gutenkunst, Waterfall, Casey, Brown, Myers,
  and Sethna}]{gutenkunst2007universally}
\bibinfo{author}{R.~N. Gutenkunst}, \bibinfo{author}{J.~J. Waterfall},
  \bibinfo{author}{F.~P. Casey}, \bibinfo{author}{K.~S. Brown},
  \bibinfo{author}{C.~R. Myers}, \bibinfo{author}{J.~P. Sethna},
\newblock \bibinfo{title}{Universally sloppy parameter sensitivities in systems
  biology models},
\newblock \bibinfo{journal}{PLoS Comput. Biol.} \bibinfo{volume}{3}
  (\bibinfo{year}{2007}) \bibinfo{pages}{1--8}.
  \DOIprefix\doi{10.1371/journal.pcbi.0030189}.
\bibitem[{Raue et~al.(2009)Raue, Kreutz, Maiwald, Bachmann, Schilling,
  Klingm{\"u}ller, and Timmer}]{raue_structural_2009}
\bibinfo{author}{A.~Raue}, \bibinfo{author}{C.~Kreutz},
  \bibinfo{author}{T.~Maiwald}, \bibinfo{author}{J.~Bachmann},
  \bibinfo{author}{M.~Schilling}, \bibinfo{author}{U.~Klingm{\"u}ller},
  \bibinfo{author}{J.~Timmer},
\newblock \bibinfo{title}{Structural and practical identifiability analysis of
  partially observed dynamical models by exploiting the profile likelihood},
\newblock \bibinfo{journal}{BMC bioinformatics.} \bibinfo{volume}{25}
  (\bibinfo{year}{2009}) \bibinfo{pages}{1923--1929}.
\bibitem[{Constantine et~al.(2014)Constantine, Dow, and
  Wang}]{constantine_active_2014}
\bibinfo{author}{P.~G. Constantine}, \bibinfo{author}{E.~Dow},
  \bibinfo{author}{Q.~Wang},
\newblock \bibinfo{title}{Active subspace methods in theory and practice:
  applications to kriging surfaces},
\newblock \bibinfo{journal}{SIAM J. Sci. Comput.} \bibinfo{volume}{36}
  (\bibinfo{year}{2014}) \bibinfo{pages}{A1500--A1524}.
\bibitem[{Transtrum et~al.(2010)Transtrum, Machta, and
  Sethna}]{transtrum2010nonlinear}
\bibinfo{author}{M.~K. Transtrum}, \bibinfo{author}{B.~B. Machta},
  \bibinfo{author}{J.~P. Sethna},
\newblock \bibinfo{title}{Why are nonlinear fits to data so challenging?},
\newblock \bibinfo{journal}{Phys. Rev. Lett.} \bibinfo{volume}{104}
  (\bibinfo{year}{2010}) \bibinfo{pages}{060201}.
\bibitem[{Transtrum and Qiu(2014)}]{transtrum2014model}
\bibinfo{author}{M.~K. Transtrum}, \bibinfo{author}{P.~Qiu},
\newblock \bibinfo{title}{Model reduction by manifold boundaries},
\newblock \bibinfo{journal}{Phys. Rev. Lett.} \bibinfo{volume}{113}
  (\bibinfo{year}{2014}) \bibinfo{pages}{098701}.
\bibitem[{Coifman and Lafon(2006)}]{coifman_diffusion_2006}
\bibinfo{author}{R.~R. Coifman}, \bibinfo{author}{S.~Lafon},
\newblock \bibinfo{title}{Diffusion maps},
\newblock \bibinfo{journal}{Appl. Comput. Harmon. Anal.} \bibinfo{volume}{21}
  (\bibinfo{year}{2006}) \bibinfo{pages}{5 -- 30}.
  \DOIprefix\doi{https://doi.org/10.1016/j.acha.2006.04.006},
  \bibinfo{note}{special Issue: Diffusion Maps and Wavelets}.
\bibitem[{Lafon(2004)}]{lafon2004diffusion}
\bibinfo{author}{S.~S. Lafon}, \bibinfo{title}{Diffusion maps and geometric
  harmonics}, Ph.D. thesis, Yale University PhD dissertation,
  \bibinfo{year}{2004}.
\bibitem[{Jolliffe(1986)}]{jolliffe1986principal}
\bibinfo{author}{I.~T. Jolliffe},
\newblock \bibinfo{title}{Principal component analysis and factor analysis},
\newblock in: \bibinfo{booktitle}{Principal component analysis},
  \bibinfo{publisher}{Springer}, \bibinfo{year}{1986}, pp.
  \bibinfo{pages}{115--128}.
\bibitem[{Coifman and Hirn(2014)}]{CoifmanHirn}
\bibinfo{author}{R.~R. Coifman}, \bibinfo{author}{M.~J. Hirn},
\newblock \bibinfo{title}{Diffusion maps for changing data},
\newblock \bibinfo{journal}{Appl. Comput. Harmon. Anal.} \bibinfo{volume}{36}
  (\bibinfo{year}{2014}) \bibinfo{pages}{79--107}.
\bibitem[{Silver et~al.(2017)Silver, Hubert, Schrittwieser, Antonoglou, Lai,
  Guez, Lanctot, Sifre, Kumaran, Graepel et~al.}]{silver2017mastering}
\bibinfo{author}{D.~Silver}, \bibinfo{author}{T.~Hubert},
  \bibinfo{author}{J.~Schrittwieser}, \bibinfo{author}{I.~Antonoglou},
  \bibinfo{author}{M.~Lai}, \bibinfo{author}{A.~Guez},
  \bibinfo{author}{M.~Lanctot}, \bibinfo{author}{L.~Sifre},
  \bibinfo{author}{D.~Kumaran}, \bibinfo{author}{T.~Graepel}, et~al.,
\newblock \bibinfo{title}{Mastering chess and shogi by self-play with a general
  reinforcement learning algorithm},
\newblock \bibinfo{journal}{arXiv preprint arXiv:1712.01815}
  (\bibinfo{year}{2017}).
\bibitem[{Bodenstein(1913)}]{Bodenstein}
\bibinfo{author}{M.~Bodenstein},
\newblock \bibinfo{title}{Eine theorie der photochemischen
  reaktionsgeschwindigkeiten},
\newblock \bibinfo{journal}{Zeitschrift f{\"u}r physikalische Chemie}
  \bibinfo{volume}{85} (\bibinfo{year}{1913}) \bibinfo{pages}{329--397}.
\bibitem[{Bello-Rivas(2017)}]{juan_m_bello_rivas_2017_581667}
\bibinfo{author}{J.~M. Bello-Rivas}, \bibinfo{title}{jmbr/diffusion-maps
  0.0.1}, \bibinfo{year}{2017}. \DOIprefix\doi{10.5281/zenodo.581667}.
\bibitem[{Rawlings and Ekerdt(2002)}]{rawlings_chemical_2004}
\bibinfo{author}{J.~B. Rawlings}, \bibinfo{author}{J.~G. Ekerdt},
  \bibinfo{title}{Chemical reactor analysis and design fundamentals},
  \bibinfo{publisher}{Nob Hill Pub, Llc}, \bibinfo{year}{2002}.
\bibitem[{Achard and De~Schutter(2006)}]{achard_complex_2006}
\bibinfo{author}{P.~Achard}, \bibinfo{author}{E.~De~Schutter},
\newblock \bibinfo{title}{Complex parameter landscape for a complex neuron
  model},
\newblock \bibinfo{journal}{PLoS Comput. Biol.} \bibinfo{volume}{2}
  (\bibinfo{year}{2006}) \bibinfo{pages}{1--11}.
  \DOIprefix\doi{10.1371/journal.pcbi.0020094}.
\bibitem[{Johnson and Goody(2011)}]{johnson_original_2011}
\bibinfo{author}{K.~A. Johnson}, \bibinfo{author}{R.~S. Goody},
\newblock \bibinfo{title}{The original {M}ichaelis constant: translation of the
  1913 {M}ichaelis--{M}enten paper},
\newblock \bibinfo{journal}{Biochemistry} \bibinfo{volume}{50}
  (\bibinfo{year}{2011}) \bibinfo{pages}{8264--8269}.
\bibitem[{Michaelis and Menten(1913)}]{michaelis_kinetik_1913}
\bibinfo{author}{L.~Michaelis}, \bibinfo{author}{M.~Menten},
\newblock \bibinfo{title}{Die kinetik der inwertinwirkung.},
\newblock \bibinfo{journal}{Biochemestry}  (\bibinfo{year}{1913})
  \bibinfo{pages}{333--369}.
\bibitem[{Heineken et~al.(1967)Heineken, Tsuchiya, and
  Aris}]{heineken_mathematical_1967}
\bibinfo{author}{F.~Heineken}, \bibinfo{author}{H.~Tsuchiya},
  \bibinfo{author}{R.~Aris},
\newblock \bibinfo{title}{On the mathematical status of the pseudo-steady state
  hypothesis of biochemical kinetics},
\newblock \bibinfo{journal}{Math. Biosci.} \bibinfo{volume}{1}
  (\bibinfo{year}{1967}) \bibinfo{pages}{95--113}.
\bibitem[{Segel and Slemrod(1989)}]{segel_quasi-steady-state_1989}
\bibinfo{author}{L.~A. Segel}, \bibinfo{author}{M.~Slemrod},
\newblock \bibinfo{title}{The quasi-steady-state assumption: a case study in
  perturbation},
\newblock \bibinfo{journal}{SIAM Rev. Soc. Ind. Appl. Math.}
  \bibinfo{volume}{31} (\bibinfo{year}{1989}) \bibinfo{pages}{446--477}.
\bibitem[{Weisz and Hicks(1962)}]{weisz_behaviour_1995}
\bibinfo{author}{P.~Weisz}, \bibinfo{author}{J.~Hicks},
\newblock \bibinfo{title}{The behaviour of porous catalyst particles in view of
  internal mass and heat diffusion effects},
\newblock \bibinfo{journal}{Chem. Eng. Sci.} \bibinfo{volume}{17}
  (\bibinfo{year}{1962}) \bibinfo{pages}{265 -- 275}.
  \DOIprefix\doi{https://doi.org/10.1016/0009-2509(62)85005-2}.
\bibitem[{Nash(1956)}]{nash1956imbedding}
\bibinfo{author}{J.~Nash},
\newblock \bibinfo{title}{The imbedding problem for riemannian manifolds},
\newblock \bibinfo{journal}{Ann. Math.}  (\bibinfo{year}{1956})
  \bibinfo{pages}{20--63}.
\bibitem[{Takens(1981)}]{takens_detecting_1981}
\bibinfo{author}{F.~Takens},
\newblock \bibinfo{title}{Detecting strange attractors in turbulence},
\newblock in: \bibinfo{booktitle}{Dynamical systems and turbulence, Warwick
  1980}, \bibinfo{publisher}{Springer}, \bibinfo{year}{1981}, pp.
  \bibinfo{pages}{366--381}.
\bibitem[{Whitney(1944)}]{whitney_self-intersections_1994}
\bibinfo{author}{H.~Whitney},
\newblock \bibinfo{title}{The self-intersections of a smooth n-manifold in
  2n-space},
\newblock \bibinfo{journal}{Ann. Math.}  (\bibinfo{year}{1944})
  \bibinfo{pages}{220--246}.
\bibitem[{Chiavazzo et~al.(2017)Chiavazzo, Covino, Coifman, Gear, Georgiou,
  Hummer, and Kevrekidis}]{ChiavazzoE5494}
\bibinfo{author}{E.~Chiavazzo}, \bibinfo{author}{R.~Covino},
  \bibinfo{author}{R.~R. Coifman}, \bibinfo{author}{C.~W. Gear},
  \bibinfo{author}{A.~S. Georgiou}, \bibinfo{author}{G.~Hummer},
  \bibinfo{author}{I.~G. Kevrekidis},
\newblock \bibinfo{title}{Intrinsic map dynamics exploration for uncharted
  effective free-energy landscapes},
\newblock \bibinfo{journal}{Proc. Natl. Acad. Sci.} \bibinfo{volume}{114}
  (\bibinfo{year}{2017}) \bibinfo{pages}{E5494--E5503}.
  \DOIprefix\doi{10.1073/pnas.1621481114}.
\bibitem[{Georgiou et~al.(2017)Georgiou, Bello-Rivas, Gear, Wu, Chiavazzo, and
  Kevrekidis}]{georgiou2017exploration}
\bibinfo{author}{A.~S. Georgiou}, \bibinfo{author}{J.~M. Bello-Rivas},
  \bibinfo{author}{C.~W. Gear}, \bibinfo{author}{H.-T. Wu},
  \bibinfo{author}{E.~Chiavazzo}, \bibinfo{author}{I.~G. Kevrekidis},
\newblock \bibinfo{title}{An exploration algorithm for stochastic simulators
  driven by energy gradients},
\newblock \bibinfo{journal}{Entropy} \bibinfo{volume}{19}
  (\bibinfo{year}{2017}) \bibinfo{pages}{294}.
\bibitem[{Mahalanobis(1936)}]{Mahalanobis}
\bibinfo{author}{P.~C. Mahalanobis},
\newblock \bibinfo{title}{On the generalized distance in statistics},
\newblock \bibinfo{organization}{Nat. Inst. of Sci. India},
  \bibinfo{year}{1936}.

\end{thebibliography}

\end{document}


\begin{frontmatter}
\title{Manifold learning for parameter reduction,\\
Supplementary Materials}
\author[add1]{Alexander Holiday\corref{cor2}}
\author[add2]{Mahdi Kooshkbaghi\corref{cor2}}
\author[add2]{Juan M. Bello-Rivas}
\author[add1]{C. William Gear}
\author[add3]{Antonios Zagaris\corref{cor1}}
\ead{antonios.zagaris@asml.com}
\author[add1,add2,add4]{Ioannis G. Kevrekidis\corref{cor1}}
\ead{yannisk@jhu.edu}
\address[add1]{Chemical and Biological Engineering, Princeton University, USA}
\address[add2]{The Program in Applied and Computational Mathematics (PACM), 
Princeton University, USA}
\address[add3]{Wageningen Bioveterinary Research, Wageningen UR, The Netherlands}
\address[add4]{Department of Chemical and Biomolecular Engineering, John Hopkins University, USA}
\cortext[cor2]{These two authors contributed equally to this work.}
\cortext[cor1]{Corresponding author}

\end{frontmatter}

\section{Mathematical setting \label{setup}}
We start our discussion with a simple framework covering all models encountered in this manuscript; a more general setting and some theoretical comments regarding sloppiness are presented later, in Section~\ref{s-origins} herein. Specifically, we consider multivariable vector functions $\mathbf{x}(t \vert \p) \in \R^D$, with $D$ finite but arbitrary. Typically, $\mathbf{x} = [x_1,\ldots,x_D]$ is only known implicitly as the solution to some problem, e.g. to an initial-value ODE problem with $t > 0$ representing time and $\p = (p_1,\ldots,p_M) \in \ps \subset \R^M$ model parameters. Also typically, one monitors \emph{not} the entire trajectory $\fmr(\p)$, for all $t>0$, but merely a number of \emph{functionals} on it, $f_1 , \ldots , f_N$. We call the $N-$tuple $\omr(\p) = [f_1(\fmr(\p)) , \ldots , f_N(\fmr(\p))] \in \R^N$ model \emph{output} or \emph{response}. In our discussion, these functionals concretely correspond to a state variable $x_d$ observed at specific preset times, i.e. $\omr(\p) = [x_d(t_1) , \ldots , x_d(t_N)]$ for certain times $t_1,\ldots,t_N$. As $\p$ ranges over $\ps$, $\omr(\p)$ traces out a (generically) $M-$dimensional manifold called the \emph{model manifold} $\omm$. In general, we will understand that manifold as
\[
 \omm
=
 \mathrm{graph}(\omr)
=
 \{(\p,\omr(\p)) \,\vert\, \p \in \ps\}
\subset
 \R^{M+N} .
\]
The projection of $\omm$ on $\R^M\times\{0\}^N$ is injective and yields the \emph{input domain} $\ps$. The projection on $\{0\}^M\times\R^N$, on the other hand, is not guaranteed to be injective even in the important case $N > M$; we consider such a \emph{non-invertible} model in the main text. During our discussion, $\omm$ will be endowed with a metric suggested by the application (in a sense by us, the observers, and our measurement capabilities), which will turn it into a Riemannian manifold. For the time being, we postpone a discussion of this important subject to a later section, where we discuss diffusion maps (DMAPS) and their relation to the Laplace-Beltrami operator.

\section{Toy example: singular perturbation}
The system discussed in the main text is a supremely basic prototype for singularly perturbed dynamical systems,
\begin{align}
   \frac{dx}{dt}&=2-x-y , \\
   \frac{dy}{dt}&= \frac{1}{\varepsilon}(x-y) .
\label{xy-model}
\end{align}
For $\varepsilon \ll 1$, trajectories are attracted to the slow manifold which, at leading order, has the form $y=x$. In our numerical experiments, we kept the initial condition for the slow variable at $x_0=-1$, fixed $t_1=0.5$, $t_2=1.0$ and $t_3=1.5$ and sampled the input space $(\varepsilon,y_0)$ for small parameter values ($10^{-3}\leq \varepsilon \leq 1$) and fast variable initializations ($3\leq y_0 \leq5$). The map from input to output space assumes the form
\begin{equation}
	[3,5]\times[10^{-3},1]\supset\Theta\ni(\varepsilon,y_0) \mapsto  \omr(\varepsilon,y_0)=
    [y(t_1) , y(t_2) , y(t_3)]\in \mathbb{R}^3 ,
\end{equation}
with a closed-form expression easily derived by quadrature on~\eqref{xy-model}.

\begin{figure}[!htb]
\centering
\includegraphics[width=0.50\textwidth]{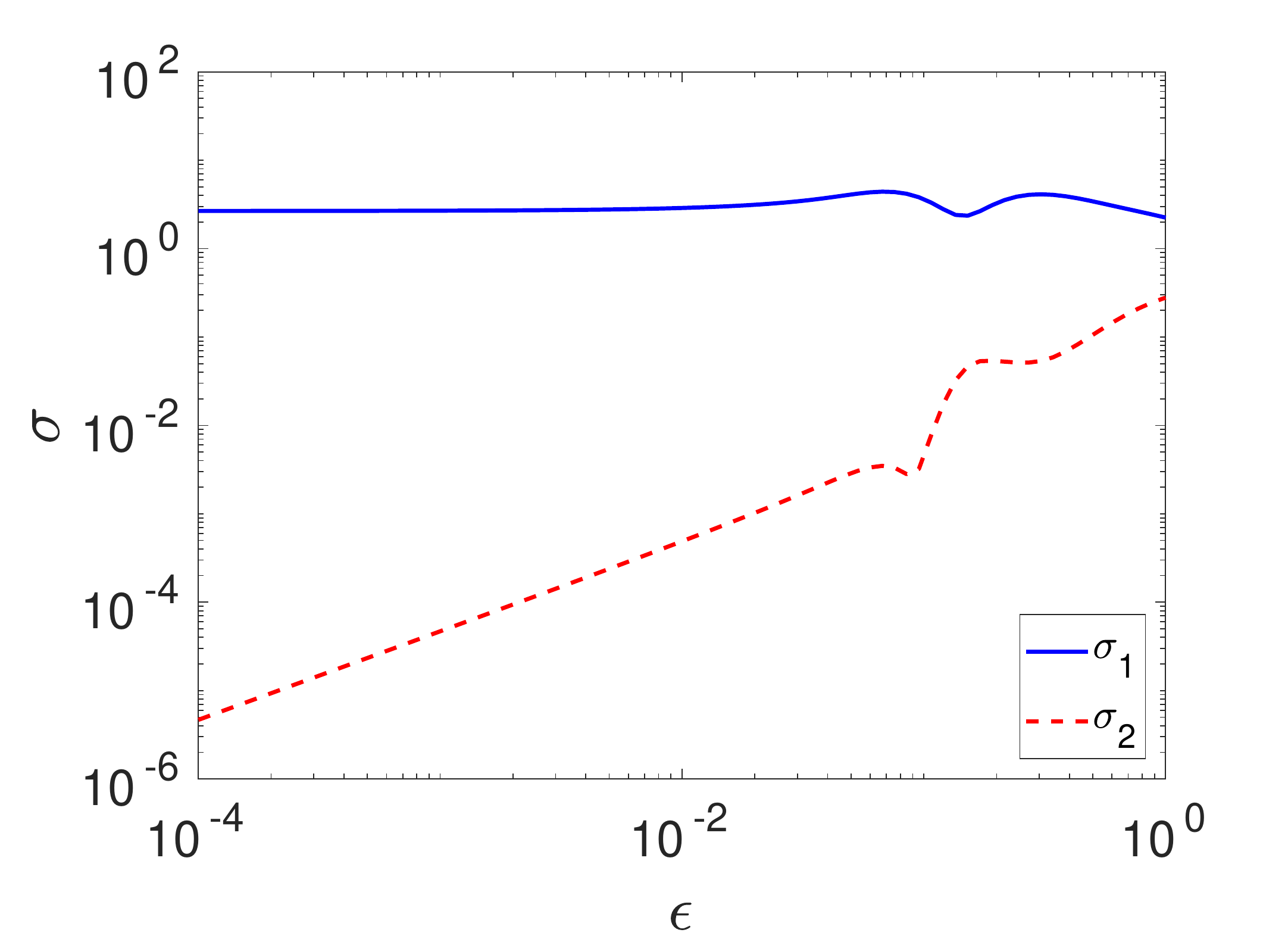}
\caption{Singular values $\sigma$ of the transformation 
map from the parameter space to the model response manifold for the singular perturbation prototype~\eqref{xy-model} and the trajectory starting at the initial condition $(x_0, y_0)=(-1,4)$. Model output $y$ is sampled at the time instants $t_1=0.5$, $t_2=1.0$ and $t_3=1.5$.}
\label{fig:singular_values}
\end{figure}

The $3\times2$ Jacobian of the transformation is
\begin{equation}
	\boldsymbol{D}\omr = \boldsymbol{J}(\varepsilon, y_0, t)=
	\begin{bmatrix}
	\frac{\partial y(t_1)}{\partial \varepsilon} & 
	\frac{\partial y(t_1)}{\partial y_0} \\
	\frac{\partial y(t_2)}{\partial \varepsilon} &
	\frac{\partial y(t_2)}{\partial y_0} \\
	\frac{\partial y(t_3)}{\partial \varepsilon} &
    \frac{\partial y(t_3)}{\partial y_0} \\		
	\end{bmatrix} ,
\end{equation}
and its singular value decomposition (SVD) yields the transformation rank and, eventually, the model manifold dimensionality. In Fig.~\ref{fig:singular_values}, we have plotted the singular values of the transformation against $\varepsilon$ for the trajectory initialized at $(x_0, y_0)=(-1,4)$. At larger values of $\varepsilon$, the transformation has rank two (evidenced by two $\mathcal{O}(1)$ singular values) and consequently the model manifold is, in principle, two-dimensional. As we decrease $\varepsilon$, however ($\varepsilon<10^{-1}$), the smallest singular value clearly approaches zero. One may impose a threshold for that value, below which the matrix is effectively rank deficient and the model manifold one-dimensional. For even smaller $\varepsilon-$values, close to the boundary $\varepsilon = 0$, the model manifold becomes effectively zero-dimensional as also
evidenced by our Diffusion Maps calculations.

\section{Spectral geometry, diffusion maps, and multi-criteria optimization}

\subsection{The Laplace--Beltrami operator and diffusion maps}

A Riemannian manifold $(\Theta,g)$ is a smooth $m$-dimensional manifold $\Theta \subseteq \mathbb{R}^n$ endowed with a metric $g$. We can regard the metric as a device that allows us to measure distances and angles on $\Theta$. At each point of $\Theta$, the Riemannian metric $g$ can be represented as a symmetric and positive definite matrix. If $u$ is a smooth real-valued function on $\Theta$, the Laplace--Beltrami operator $\Delta$~\cite{Morita2001} is a linear operator given, in local coordinates $(z_1, \dotsc, z_m)$, by the formula
\begin{equation*}
  \Delta u = \frac{1}{\sqrt{\det g}} \sum_{i = 1}^m \sum_{j = 1}^m \frac{\partial}{\partial z_i} \left( \sqrt{\det g} \, (g^{-1})_{ij} \frac{\partial u}{\partial z_j} \right) .
\end{equation*}
Here, $\det g$ is the determinant of the matrix associated with the Riemannian metric and $g^{-1}$ is the corresponding inverse matrix. Intuitively, the eigenfunctions of the Laplace--Beltrami operator determine how heat propagates on the manifold $\Theta$.

Let $S_N = \{x_1, \dotsc, x_N\} \subset \mathbb{R}^n$ be a set of points sampled from an arbitrary probability distribution on the manifold $\Theta$. The computational complexity of traditional approximation schemes for the Laplace--Beltrami operator, such as finite differences or finite elements, scales exponentially in the dimension $n$. By contrast, the diffusion maps method (DMAPS; to be discussed below) is a data-driven approximation of $\Delta$ with a computational complexity that scales quadratically in the number of samples $N$.

To approximate $\Delta$ in the DMAPS sense, fix a scale $\epsilon > 0$ and consider the affinity matrix $A \in \mathbb{R}^{N \times N}$ with entries
\begin{equation*}
  A_{ij} = \exp\left\{ -\frac{\| x_i - x_j \|^2}{2 \epsilon} \right\}.
\end{equation*}
For each $i = 1, \dotsc, N$, set further $q_i = \sum_j A_{ij}$ and introduce the matrix $W \in \mathbb{R}^{N \times N}$ with entries
\begin{equation*}
  W_{ij} = \frac{A_{ij}}{q_i \, q_j}.
\end{equation*}
We define the (random walk) graph Laplacian $L$ as~\cite{Chung1996}
\begin{equation*}
  L = I - D^{-1} W,
\end{equation*}
where $I \in \mathbb{R}^{N \times N}$ is the identity matrix and $D \in \mathbb{R}^{N \times N}$ is the diagonal matrix with diagonal entries
\begin{equation*}
  D_{ii} = \sum_{j} W_{ij}.
\end{equation*}
Note that $L$ is determined by $\epsilon$ and the set of point-samples $S_N$. Let us now represent the smooth function $u \colon \Theta \to \mathbb{R}$ by a vector $U \in \mathbb{R}^N$ with components $U_i = u(x_i)$, for $i = 1, \dotsc, N$. It is known~\cite{Lafon2004} that
\begin{equation*}
  \lim_{\epsilon \to 0} \epsilon^{-1} \lim_{N \to \infty} \sum_{j} L_{ij} U_j = \tfrac{1}{2} \Delta u(x_i).
\end{equation*}
More precisely, it can be shown~\cite{Singer2006} that the choice
\begin{equation}
  \label{eq:optimal-epsilon}
  \epsilon = C \, N^{\frac{-2}{6 + m}},
\end{equation}
with $C > 0$ a constant depending on the geometry of $M$ but not on $N$, leads to the minimal error bound
\begin{equation}
  \label{eq:error-bound}
  \epsilon^{-1} \sum_{j} L_{ij} U_j = \tfrac{1}{2} \Delta u(x_i) + \mathcal{O}(N^{-\frac{2}{m + 6}})
  \quad \text{as} \ N \to \infty.
\end{equation}
For example, for a 2D manifold ($m = 2$), \eqref{eq:optimal-epsilon}--\eqref{eq:error-bound} establish that the error decays \emph{slowly} as $N^{-1/4}$.

Let $\psi_0, \psi_1, \psi_2, \dotsc$ be the eigenfunctions of $\Delta$ corresponding to the eigenvalues $0 = \lambda_0 \le \lambda_1 \le \lambda_2 \le \dotsc < +\infty$ and let $d \ge m$.
The diffusion map $\Psi_\epsilon \colon \mathbb{R}^n \to \mathbb{R}^d$, defined by $\Psi_\epsilon = (\mathrm{e}^{\lambda_1 \epsilon} \psi_1, \dotsc, \mathrm{e}^{\lambda_d \epsilon} \psi_d)$, is used for dimensionality reduction in manifold learning. The leading eigenfunction $\psi_0$ is not included in the definition of $\Psi_\epsilon$, because it is always a constant function and carries no information on $(\Theta,g)$.

\begin{figure}[!htp]
  \centering
  \subfigure{\includegraphics{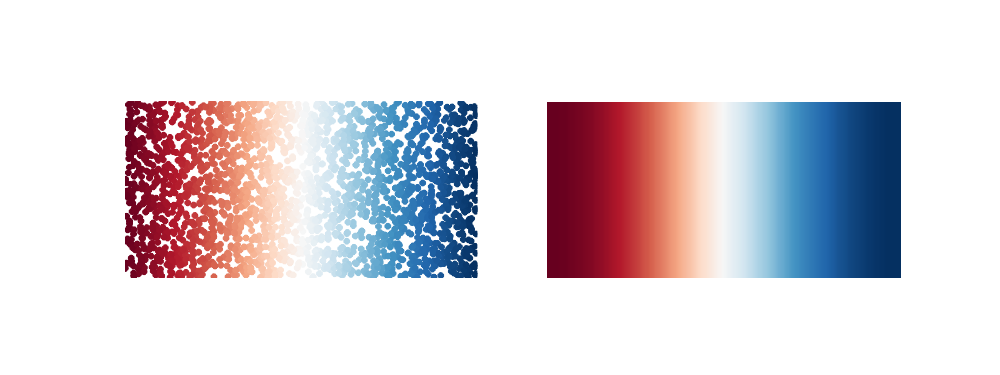}}
  \subfigure{\includegraphics{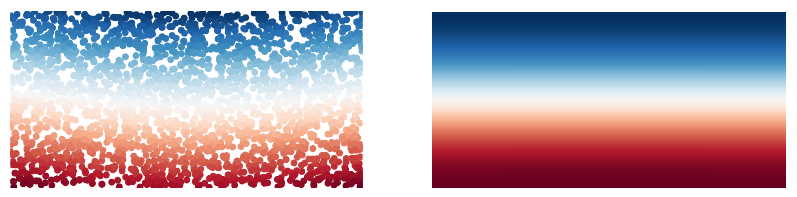}}
  \subfigure{\includegraphics{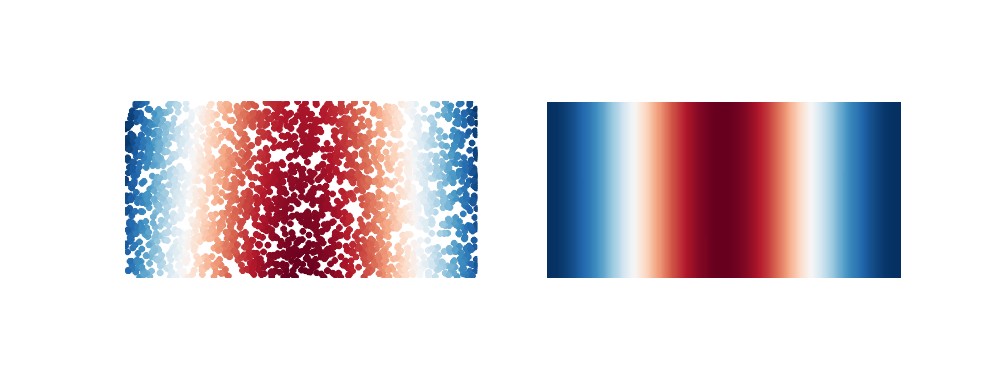}}
  \subfigure{\includegraphics{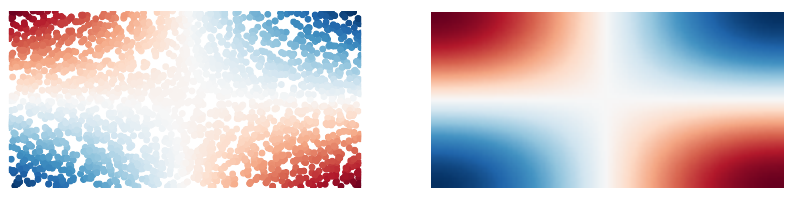}}
  \subfigure{\includegraphics{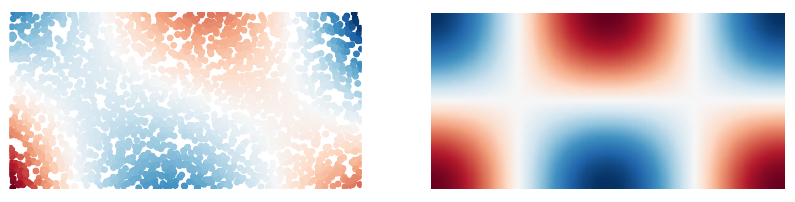}}
  \subfigure{\includegraphics{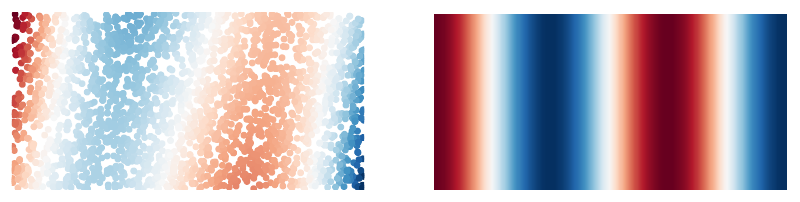}}
  \caption{Comparison of some eigenfunctions of the Laplace--Beltrami operator for a rectangle with aspect ratio 2. Columns $1$ and $3$ show eigenfunctions obtained by applying DMAPS on $10^4$ points sampled uniformly and with $\epsilon = 10^{-1}$. Columns $1$ and $3$ show the corresponding eigenfunctions obtained by the finite element method, using an adaptive scheme with Lagrange elements of degree 2.}
  \label{fig:eigenfunctions-rectangle}
\end{figure}

In practice, we calculate the spectral decomposition of the Laplace--Beltrami operator $\Delta$ by the eigenvalues and eigenvectors of the graph Laplacian $L$. However, in cases in which the dimensionality of the problem is sufficiently low, using the finite element method~\cite{Brenner2008} is faster and more accurate. This additional accuracy is crucial, if we are interested in computing a large number of eigenfunctions. To demonstrate the method, we consider the simple example of a rectangle, $\Theta = [0, 1] \times [0, \ell] \subset \mathbb{R}^2$, endowed with the Euclidean metric. In Fig.~\ref{fig:eigenfunctions-rectangle}, we show the first few eigenfunctions of the Laplace--Beltrami operator, computed using both diffusion maps and a high-accuracy finite element method for fixed aspect ratio $1/\ell = 2$. The exact spectral decomposition is given by the family of eigenfunctions $\psi_{ij}(x, y) = \cos(i \pi x) \, \cos(j \pi y / \ell)$ and eigenvalues $\lambda_{ij} = \pi^2 (i^2 + j^2/\ell^2)$, for $i, j = 0, 1, 2, \dotsc$. On one hand, the spectral decomposition of the graph Laplacian was computed using a set of $N = 10^4$ points, sampled uniformly from $\Theta$, with $\epsilon = 10^{-1}$. On the other hand, the corresponding spectral decomposition of the Laplace-Beltrami operator was computed using the FEniCS finite element library~\cite{Logg2012} with quadratic Lagrange elements on an adaptively refined mesh.

\subsection{Interpretation in terms of continuum mechanics}

Let $\Theta$ and $\mathcal{M}$ be two smooth manifolds. Given a smooth map $\mathbf{f} \colon \Theta \to \mathcal{M}$, we consider the deformation gradient, given by the Jacobian $D \mathbf{f}$, and the associated Green deformation tensor~\cite{Marsden1994} (also known as the right Cauchy--Green tensor), defined by $D \mathbf{f}^T D \mathbf{f}$. The deformation tensor can be regarded as a tool quantifying the change in $\mathbf{f}(\Theta) \subseteq \mathcal{M}$ that results from a change in $\Theta$. We elaborate on this idea in what follows.

\begin{figure}[!htp]
  \centering
  \subfigure{\includegraphics{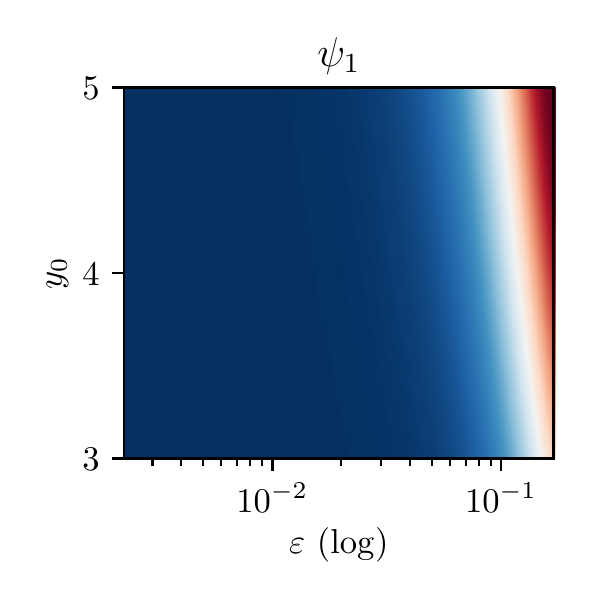}}
  \subfigure{\includegraphics{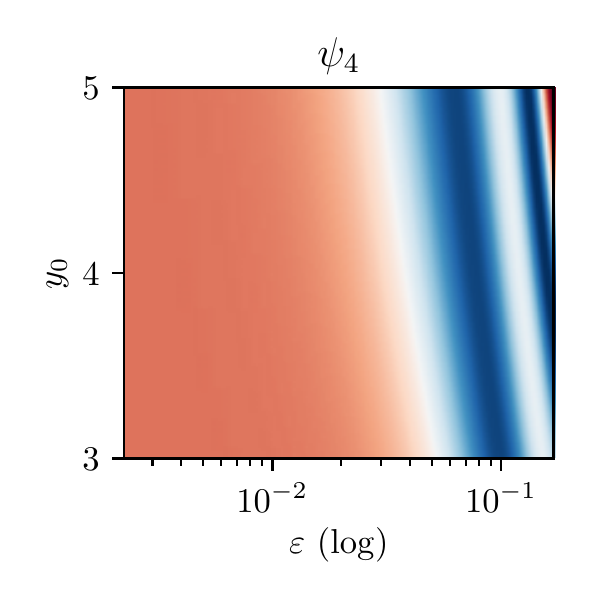}}
  \subfigure{\includegraphics{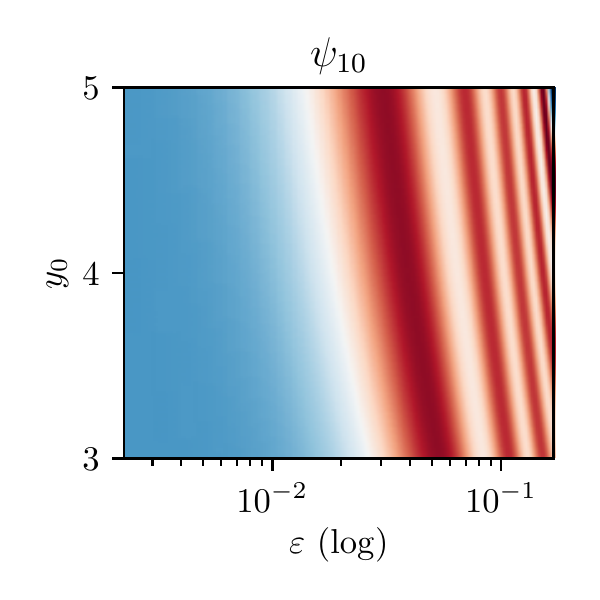}}
  \subfigure{\includegraphics{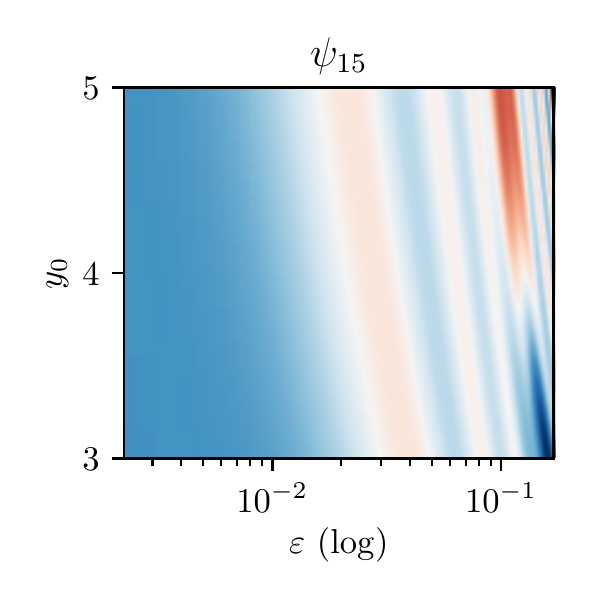}}
  \subfigure{\includegraphics{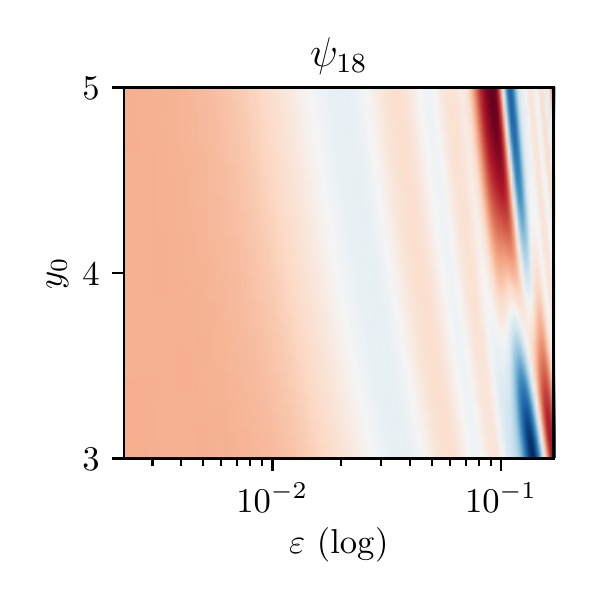}}
  \subfigure{\includegraphics{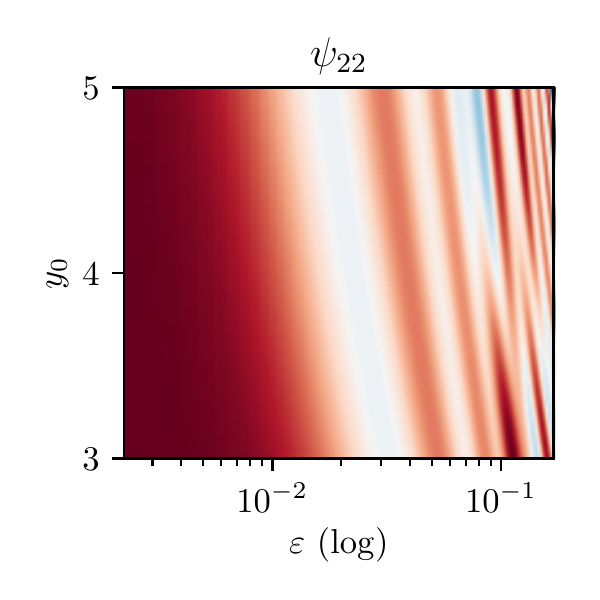}}
  \caption{Some eigenfunctions of the singularly perturbed problem.
    The eigenfunctions $\psi_4$ and $\psi_{10}$ (top) are harmonics of the first non-constant eigenfunction $\psi_1$ (also top), whereas the rest (bottom) are not.}
  \label{fig:eigenfunctions-singularly-perturbed}
\end{figure}

For the sake of concreteness, we look at a particular example: the analytically solvable, linear, singular perturbation example studied in the paper,
\begin{equation}
  \label{eq:ivp}
  \left\{
    \begin{aligned}
      &\dot{x} = 2 - x - y, \\
      \varepsilon &\dot{y} = x - y, \\
      &x(0) = -1, \\
      &y(0) = y_0.
    \end{aligned}
  \right.
\end{equation}
Using the solution of~\eqref{eq:ivp} evaluated at $t_1 = 1/2$, $t_2 = 1$, and $t_3 = 3/2$, we obtain a diffeomorphism $\mathbf{f} \colon \Theta  \to \mathcal{M}$ given by
\begin{equation*}
  (\varepsilon , y_0) \mapsto \mathbf{f}(\varepsilon, y_0) = [y(t_1 \vert \varepsilon,y_0), y(t_2 \vert \varepsilon,y_0), x(t_3 \vert \varepsilon,y_0)] .
\end{equation*}
Here, $\Theta = [3,5] \times [0.00225,0.1] \subset \mathbb{R}^2$ and $\mathcal{M} = \mathbf{f}(\Theta) \subset \mathbb{R}^3$. In our current context, the manifold $\Theta$ represents \emph{inputs} for the transformation $\mathbf{f}$ and, correspondingly, $\mathcal{M}$ will be the manifold of \emph{outputs} of $\mathbf{f}$. A perturbation of an input $x$ in a neighborhood $V \subseteq \Theta$ by a tangent vector $\Delta x \in T_x V$, with $\| \Delta x \| = h>0$, results in a change in the output of size 
\begin{equation*} 
  \| \mathbf{f}(x + \Delta x) - \mathbf{f}(x) \|^2 = \| D \mathbf{f}(x) \, \Delta x \|^2 + o(\varepsilon) = (\Delta x)^T (D \mathbf{f}^T D \varphi) \, \Delta x + o(h),
\quad\mbox{as}\
 h \downarrow 0 .
\end{equation*}
Therefore, we can endow $\Theta$ with the metric determined by the deformation tensor $g = D \mathbf{f}^T D \mathbf{f}$. Indeed, $g$ measures how close to each other are the images, under $\mathbf{f}$, of the points $x = (\varepsilon, y_0)$ and $x + \Delta x = (\varepsilon + \Delta \varepsilon, y_0 + \Delta y_0)$. It is then natural to study the eigenfunctions of the Laplace--Beltrami operator on $(\Theta,g)$. These eigenfunctions yield another parametrization~\cite{Berard1994} of $\Theta$ that reflects the sensitivity of the outputs to changes in the inputs. Some of the relevant eigenfunctions are plotted in Fig.~\ref{fig:eigenfunctions-singularly-perturbed}.

\subsection{Thoughts on multi-criteria optimization}
\label{sec:thoughts-on-multi-criteria-optimization}
\begin{figure}[ht]
  \centering
  \subfigure{\includegraphics{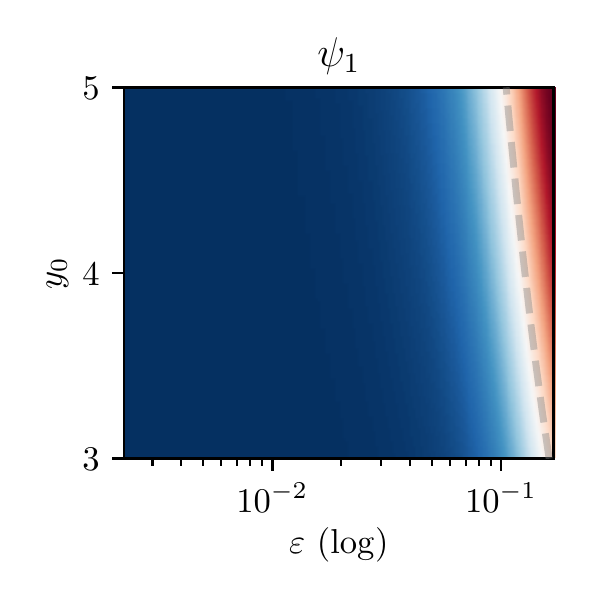}}
  \subfigure{\includegraphics{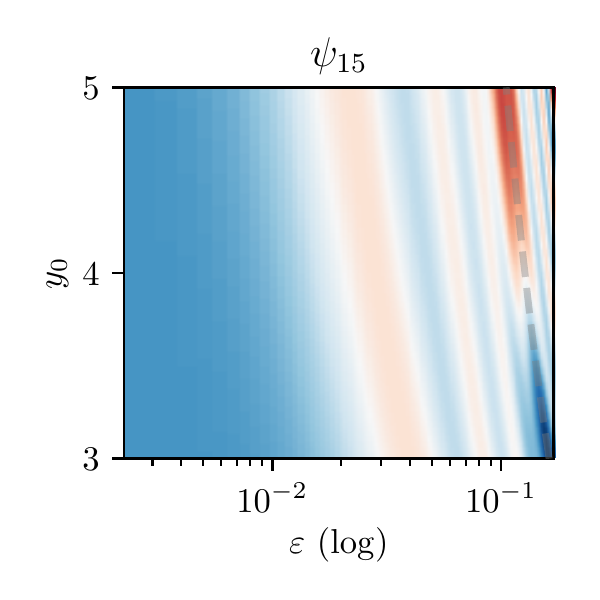}}
  \subfigure{\includegraphics{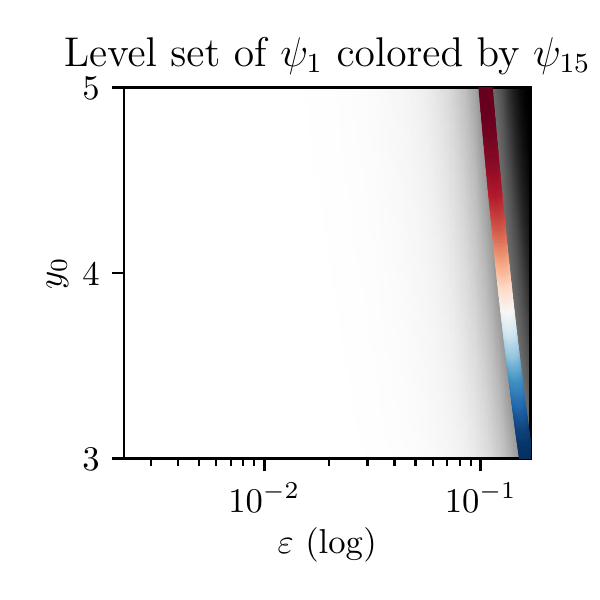}}
  \caption{Parameterization of the first non-constant eigenfunction of the singularly perturbed problem using another eigenfunction. Left: a level set of $\psi_1$ (dashed) overimposed on a heatmap of $\psi_1$, the first non-constant eigenfunction. Middle: same but with a heatmap of $\psi_{15}$, the first eigenfunction that is not a higher harmonic of $\psi_1$. Plainly, $\psi_{15}$ is monotonic on that level set. Right: a concatenation of the top panels. The gray heatmap covers the entire panel and represents $\psi_1$, whereas the colored one is confined in a neighborhood of the level set and represents $\psi_{15}$.}
  \label{fig:singularly-perturbed-level-sets}
\end{figure}

The curve $\gamma \colon [0, 1] \to \Theta$ depicted (dashed) in Fig.~\ref{fig:singularly-perturbed-level-sets} is a level set of the first non-constant eigenfunction $\psi_1$ of the Laplace--Beltrami operator on $(\Theta,g)$, i.e.
\[
 \psi_1(\gamma(t)) = c ,
\ \mbox{for some constant} \
 c \in \mathbb{R}
\ \mbox{and all} \
 t \in [0, 1] .
\]
The image of $\gamma$ under $\psi_{15}$, which is not a harmonic of $\psi_1$ (cf. Fig.~\ref{fig:eigenfunctions-singularly-perturbed}) is  increasing monotonically function \emph{along} the curve; indeed, $\mathrm{d}\psi_{15}(\gamma(t))/\mathrm{d}t$ neither changes sign nor vanishes for any $t \in [0, 1]$. As a result, we can parameterize $\gamma$ using the values of $\psi_{15}$. The fact that level sets of $\psi_1$ can be parameterized using another eigenfunction $\psi_k$ has obvious applications to multi-criteria optimization~\cite{Ehrgott2005}. In particular, if our optimization criterion follows a lexicographic ordering, in which we first seek an optimal level set $\gamma$ of $\psi_1$, the parametrization of $\gamma$ in terms of $\psi_k$ is a natural way to subsequently seek an optimal point \emph{on} $\gamma$.

\section{Regular perturbation}

The system under investigation is $\dot{x} = \varepsilon x^3 - x$, with
$\varepsilon \ll 1$ a small parameter and $x$ the system state. Using that this differential equation is separable and restricting ourselves to $x \ge 0$, we can write the analytical solution
\begin{align}
  x(t \vert \varepsilon,x_0) = \left({\varepsilon + e^{2 t} \left(\frac{1}{x_0^2} -
  \varepsilon \right)}\right)^{-1/2} .
\label{eq:regpert}
\end{align}
We are interested in the behavior of the system in the perturbation limit $\varepsilon \rightarrow 0$. The limiting solution can be found simply by omitting the cubic term from the differential equation, giving us
\begin{align}
  \lim_{\varepsilon \rightarrow 0} x(t) = x_0 e^{-t}.
\label{eq:regpert-limit}
\end{align}
The key distinction between this system and the singularly perturbed model above is that the initial conditions \emph{does} influence the trajectories, even in the limiting case of small $\varepsilon$. Varying both $\varepsilon$ and $x_0$ for large $\varepsilon-$values, one obtains the model response~\ref{eq:regpert} that depends nontrivially on {\em both} $\varepsilon$ and $x_0$. For small values of $\varepsilon$, on the other hand, trajectories converge to the limiting solution~\ref{eq:regpert-limit}, where $\varepsilon$ does not affect the model response. This is precisely what we see in the model manifold depicted in Fig.~5 of the main text. For large values of $\varepsilon$, the model manifold is two dimensional and variations in both $x_0$ and $\varepsilon$ affect the model response. At smaller values of $\varepsilon$, the model manifold converges to a 1D object parameterized by $x_0$. The exact form of that object is simply the ray $f(x_0)= x_0 (e^{-t_1},e^{-t_2},e^{-t_3})$, with $x_0>0$, bearing no dependence on $\varepsilon$.

To create Fig.~5 in the paper, we fixed $t_1 = 0.25$, $t_2 = 1.0$ and $t_3 = 1.75$ and set the model response to $\mathbf{f}(x_0,\varepsilon) = [x(t_1 \vert x_0,\varepsilon) , x(t_2 \vert x_0,\varepsilon) , x(t_3 \vert x_0,\varepsilon)] \in \mathbb{R}^3$, with $x(\cdot \vert x_0,\varepsilon)$ given by~\eqref{eq:regpert}. We then drew $2500$ points from the rectangle $x_0 \in [1.0, 2.5]$ and $\log(\varepsilon) \in [-3, -1]$, uniformly in these two parameters, and performed DMAPS on the input--output combinations with $\epsilon = 5.0$ as kernel scale. Coloring parameter space by the resulting $\phi_1$ values gives Fig.~5 in the paper, showing that $\phi_1$ ``discovers'' the regularly perturbed nature of the problem in a data-driven manner.
\section{The ABC system}
In this section, we identify the singularly perturbed regime for the dynamics associated with the reaction scheme
\[
 \mathrm{A}
\xrightleftharpoons[k_{-1}]{k_1}
 \mathrm{B}
\xrightarrow[]{k_2}
 \mathrm{C} .
\]
As we will see, the dynamics in that regime is driven by a {\em single} effective parameter $\keff$.
\subsection{Exact solution}
The evolution of the molar concentrations is dictated by the linear ODEs
\[
 \frac{d}{dt}
\left[\begin{array}{c}
 A \\ B \\ C
\end{array}\right]
=
\left[\begin{array}{ccc}
 -k_1 & k_{-1} & 0 \\ k_1 & -(k_{-1}+k_2) & 0 \\ 0 & k_2 & 0
\end{array}\right]
\left[\begin{array}{c}
 A \\ B \\ C
\end{array}\right] ,
\ \mbox{subject to} \
\left[\begin{array}{c}
 A(0) \\ B(0) \\ C(0)
\end{array}\right]
=
\left[\begin{array}{c}
 A_0 \\ B_0 \\ C_0
\end{array}\right] .
\]
This system has the explicit solution
\be
\begin{array}{rcl}
\left[\begin{array}{c}
 A(t) \\ B(t) \\ C(t)
\end{array}\right]
&=&
\displaystyle
 \frac{\lp A_0 + (\lp+k_2)B_0}{\dl}
\left[\begin{array}{ccc}
 -(1+\lm/k_2) \\ \lm/k_2 \\ 1
\end{array}\right]
 \mathrm{e}^{\lm t}
\\
&{}&
\displaystyle
+
 \frac{\lm A_0 + (\lm+k_2)B_0}{\dl}
\left[\begin{array}{ccc}
 -(1+\lp/k_2) \\ \lp/k_2 \\ 1
\end{array}\right]
 \mathrm{e}^{\lp t}
\\
&{}&
\displaystyle
+
 (A_0+B_0+C_0)
\left[\begin{array}{ccc}
 0 \\ 0 \\ 1
\end{array}\right] ,
\end{array}
\label{ABC-soln}
\ee
where $\dl = \lp-\lm > 0$ is the difference between the nonzero eigenvalues
\be
 \lpm
=
 \frac{1}{2}
\left[
 -(k_{-1}+k_1+k_2)
\pm
 \sqrt{(k_{-1}+k_1+k_2)^2 - 4 k_1 k_2}
\right]
<
 0 .
\label{ABC-lambda}
\ee
Among them, $\lm$ controls the fast time scale (\emph{transient dynamics}) and $\lp$ the rate of approach to the steady state $(A_*,B_*,C_*)=(0,0,A_0+B_0+C_0)$ (\emph{slow dynamics}). The third eigenvalue is identically zero and associated with the conservation law $A(t)+B(t)+C(t)=A_0+B_0+C_0$. For our purposes, we assume a fixed initial state $(A_0,B_0,C_0)=(1,0,0)$ (in arbitrary units) and variable reaction rates $\p = (k_{-1},k_1,k_2)$. Observing the system takes the form of monitoring the product concentration at preset times, $\omr(\p) = [C(t_1 \vert \p) , \ldots , C(t_5 \vert \p)]$. Here, the times $t_1, \dots, t_5$ form a uniform grid on the interval $[-0.5/\lambda_+ \,,\, -5/\lambda_+]$, chosen to resolve the slow dynamics.

\subsection{Singularly perturbed regime}
Rearranging terms from slowest to fastest in \eqref{ABC-soln} yields
\be
 C(t \vert \p)
=
 1
+
 \frac{\lm}{\dl}
 \mathrm{e}^{\lp t}
\left(
 1
+
 \frac{\lp}{\lm}
 \mathrm{e}^{r \lp t}
\right) ,
\quad\mbox{with}\
 r = \frac{\dl}{|\lp|} > 0 .
\label{ABC-C-soln}
\ee
To resolve the slow dynamics, we select $t_1 = \alpha/|\lp|$ for $\alpha \approx 0.5$.
In the presence of time scale disparity, the fast component should be negligible by that time already,
which necessitates that $\vert\lp/\lm\vert \mathrm{e}^{-\alpha r} \ll 1$.
This order relation determines the asymptotic regime in parameter space.
Since the ratio $\vert\lp/\lm\vert$ depends only algebraically on $\p$,
time scale separation must arise from the exponential term
for all $r> r_*$ with $\alpha r_* \gg 1$.
As an indication, the value $r_* = 6$ reduces the exponential term to approximately $0.05$.\\

\begin{figure}[!htp]
\centering
\includegraphics[width=0.5\textwidth]{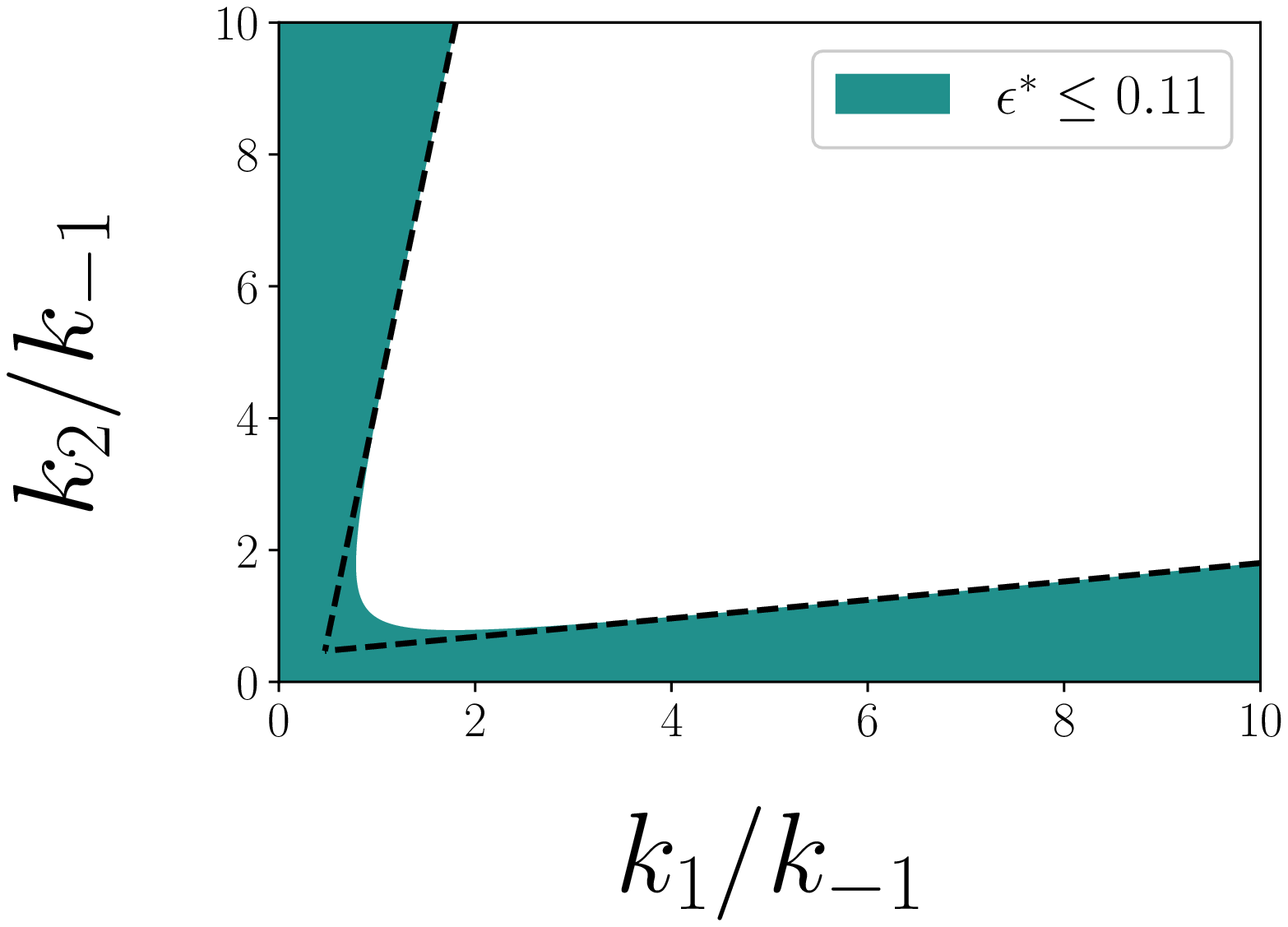}
\caption{\label{f-regim}
The region on the $(k_1/k_{-1} \,,\, k_2/k_{-1})-$plane within which there is timescale separation.
}
\end{figure}

To identify the asymptotic regime explicitly, we work with the compactification
\be
 \eps
=
 \frac{1}{4}
\left[
 1 - \left(\frac{r}{r+2}\right)^2
\right]
=
 \frac{k_1 k_2}{(k_{-1}+k_1+k_2)^2}
<
 \eps_*
=
 \frac{1}{4}
\left[
 1 - \left(\frac{r_*}{r_*+2}\right)^2
\right] ,
\label{ABC-asympt-bound}
\ee
where we have used \eqref{ABC-lambda} to express $r$ in terms of the kinetic parameters.
The composite parameter $\eps$ decreases with $r$ from $1/4$ to zero,
hence timescale disparity exists for $\eps < \eps_* \ll 1$: $\eps$ acts as a \emph{small parameter}.
As an indication, the value $r_*=6$ above yields $\eps_* \approx 0.11$.
To understand the quadratic curve~\eqref{ABC-asympt-bound} bounding the singularly perturbed regime,
we introduce the transformation
\[
\left[\begin{array}{c}
 \kappa_1 \\ \kappa_2
\end{array}\right]
=
 \frac{1}{\sqrt2}
\left[\begin{array}{rr}
 1+\sqrt{1-4\eps_*} & -1+\sqrt{1-4\eps_*} \\ -1+\sqrt{1-4\eps_*} & 1+\sqrt{1-4\eps_*}
\end{array}\right]
\left[\begin{array}{c}
 k_1/k_{-1} - 2 \eps_*/(1 - 4 \eps_*) \\ k_2/k_{-1} - 2 \eps_*/(1 - 4 \eps_*)
\end{array}\right] ,
\]
which factorizes \eqref{ABC-asympt-bound} as $\kappa_1 \kappa_2 < 2 \eps_*/(1 - 4 \eps_*)$.
In the $(\kappa_1,\kappa_2)-$plane, the asymptotic regime is bounded by two hyperbolas in the first and third quadrants.
Since the $\kappa_1-$ and $\kappa_2-$axes align with the axes in the $(k_1/k_{-1},k_2/k_{-1})$ plane, at leading order,
the regime $\eps_* \downarrow 0$ becomes the narrow sliver shown in Fig.~\ref{f-regim}.

\subsection{Effective parameter}
By construction, data generated by parameter values in the asymptotic regime
are well-described by the slow component alone, whose decay rate is
\[
 \vert \lp \vert
=
 \frac{1 - \sqrt{1-4\eps}}{2}
 (k_{-1}+k_1+k_2)
\sim
 \frac{k_1 k_2}{k_{-1}+k_1+k_2}
=
 \keff .
\]
For parameter values in that regime, the observable effectively reads (cf.~\eqref{ABC-C-soln})
\[
 C(t \vert \p)
=
 1
-
 \mathrm{e}^{-\keff t} ,
\quad\mbox{since}\
 \frac{\lp}{\lm}
 \mathrm{e}^{-r \lp t}
\ll
 1
\ \mbox{and} \
 \frac{\lm}{\dl}
=
-
\frac{
 1 + \sqrt{1 - 4 \eps}
}{
 2 \sqrt{1 - 4 \eps}
}
\sim
 -1 .
\]
This is the formula we based our discussion in the main text on, together with its domain of applicability, cf.~\eqref{ABC-asympt-bound}. This leading order result shows the observable to depend solely on the decay rate $\keff$ in the singularly perturbed regime. That regime is foliated by (subsets of) the level sets of $\keff$, with the model response remaining largely constant on each such surface.\\

It is interesting to note that, next to sloppiness, the setup above also exhibits \emph{structural non-identifiability}. Indeed, the observable only depends on the two parameter combinations $\lpm$, see the exact formula \eqref{ABC-C-soln}. It is important to understand that this effect is distinct from time scale disparity and thus not limited to the singularly perturbed regime. This further entails that the parameter space is foliated by curves along which $\lpm$ remain constant, with all points on any such curve yielding \emph{identical} model responses.
Since $\keff$ is merely another name for $\lp$,
each surface $\keff = const.$ is itself foliated by these curves. In other words, the curves of \emph{identical} model response (constant $\lpm$)
foliate surfaces of \emph{leading-order identical} model response (constant $\lp$) which, in turn, foliate the singularly perturbed regime. We remark once again that sloppiness and non-identifialibity are pertinent to both the system and the monitoring protocol employed. Allowing $B_0>0$ in the initial condition or observing $A(t)$ and/or $B(t)$, next to $C(t)$, suffices to lift the non-identifiability.

\section{The Michaelis--Menten--Henri system}
We consider the prototypical chemical pathway \cite{JG11,MM13}
\[
 \mathrm{S + E}
\
 \xrightleftharpoons[k_{-1}]{k_1}
\
 \mathrm{C}
\
 \xrightarrow{k_2}
\
 \mathrm{P + E} ,
\label{s2c2p}
\]
modeling the two-step conversion of a substrate $\mathrm{S}$ into product $\mathrm{P}$ through the mediation of an enzyme $\mathrm{E}$.
The constituent concentrations evolve under
\be
\begin{array}{rclcl}
 S' &=& -k_1 E S + k_{-1} C ,
\\
 C' &=& \ \ \, k_1 E S - (k_{-1} + k_2) C ,
\\
 E' &=& -k_1 E S + (k_{-1} + k_2) C ,
\\
 P' &=& \ \ \, k_2 C ,
\end{array}
\label{SCEP-ODE}
\ee
supplemented by the initial concentrations $S_0$, $E_0$, $C_0$ and $P_0$.
This system has two exact conservation laws expressing mass balance,
\[
 S+C+P = S_0+C_0+P_0 =: S_T
\quad\mbox{and}\quad
 C+E = C_0+E_0 =: E_T .
\]
Typically, one uses these to eliminate the last two ODEs, thus obtaining 
\be
\begin{array}{rclcl}
 S' &=& -k_1 (E_T-C) S + k_{-1} C ,
\\
 C' &=& \ \ \, k_1 (E_T-C) S - (k_{-1} + k_2) C .
\end{array}
\label{SC-ODE}
\ee
This is the classical Michaelis--Menten--Henri system in \emph{dimensional} form \cite{SS89}.
In a typical experimental setting, $C_0=P_0=0$ so that $S_T=S_0$.
We follow this setup here and consider a problem with five parameters,
the three kinetic constants $k_{\pm1}$ and $k_2$ and the total concentrations $S_T$ and $E_T$.
To further emulate an experimental setting,
we set our observable to be the product concentration,
whose time course is determined by the IVP
\be
 P' = k_2 C ,
\quad\mbox{subject to}\
 P_0 = 0 .
\label{P-IVP}
\ee
Equations~\eqref{SC-ODE} and \eqref{P-IVP} represent the \emph{original form} of the Michaelis--Menten--Henri system.

\subsection{System dynamics}
The multiscale dynamics of this system has been analyzed in a series of publications,
among which the landmark articles \cite{HTA67,SS89}.
The salient features of these two articles is the system nondimensionalizations they propose and, in particular,
the small parameter they use to define the singularly perturbed regime.
Specifically, the authors of \cite{HTA67} work with the small parameter $\bar{\eps} = E_T/S_T$,
whereas those of \cite{SS89} use $\eps = E_T/(S_T+K_M)$;
the asymptotic regime is defined as $\bar{\eps} \ll 1$ and $\eps \ll 1$, respectively.\\

To nondimensionalize the system, we draw inspiration from \cite{SS89} and define the new parameter set
\be
 (\sigma , K_M , V_M , \kappa , \eps)
=
\left(
 \frac{S_T}{K_M}
\,,\,
 \dst\frac{k_{-1} + k_2}{k_1}
\,,\,
 k_2 E_T
\,,\,
 \dst\frac{k_{-1}}{k_2}
\,,\,
 \dst\frac{E_T}{S_T + K_M}
\right) .
\label{params-new}
\ee
Here, $S_T$ and $K_M$ have units of concentration, $V_M$ of reaction speed and $\kappa$ and $\eps$ are non-dimensional.
The inverse of this bijection is
\[
 S_T
=
 \sigma K_M ,
\
 E_T
=
 \eps K_M (\sigma + 1) ,
\
 (k_{-1} , k_1 , k_2)
=
 \frac{V_M}{\eps K_M (\sigma + 1)}
\left(
 \kappa , \frac{\kappa+1}{K_M} , 1
\right) .
\]
Additionally, we nondimensionalize time and reactant concentrations through
\[
 \tau
=
 \frac{t}{t_s}
\quad\mbox{and}\quad
 (s , c , p)
=
\left(
 \frac{S}{S_T} , \frac{C}{\bar{C}} , \frac{P}{S_T}
\right) ,
\]
where the slow timescale $t_S$ and complex concentration estimate $\bar{C}$ are \cite{SS89}
\[
 t_s
=
 \frac{S_T + K_M}{V_M}
\quad\mbox{and}\quad
 \bar{C}
=
 \frac{E_T S_T}{S_T + K_M} .
\]
The new system of ODEs becomes
\be
\begin{array}{rclcl}
 \dot{s}
&=&
\dst
 (\kappa+1)
\left[
-
 \left(1 + \sigma\right)
 s
+
 \sigma
 c s
+
 \frac{\kappa}{\kappa+1}
 c
\right] ,
\vspace{2mm}\\
 \eps
 \dot{c}
&=&
\dst
 (\kappa+1)
\left[
\ \ \,
 \left(1 + \sigma\right)
 s
-
 \sigma
 c s
-
 c
\right] ,
\end{array}\!\!
\label{SC-ODE-new}
\ee
with initial conditions $s_0 = 1$ and $c_0 = 0$.
The observable $p$ evolves under
\be
 \dot{p}
=
 c ,
\quad\mbox{subject to}\
 p_0 = 0 .
\label{P-IVP-new}
\ee
This is the \emph{rescaled form} of the system, studied in the main text.

\section{The catalyst pellet system}
When a heterogeneous chemical reaction takes place in a porous catalyst pellet, transport of the reactants
to the pellet surface and then through its pores to the catalyst sites plays an important role in determining the overall reaction rate. The Thiele modulus ($\Phi$) is a dimensionless parameter determining the relative strengths of diffusive transport and reaction. For a spherical particle, $\Phi$ can be defined as
\be
\Phi = R \sqrt{k/D} ,
\ee
with $R$ the particle radius, $k$ the first order reaction rate constant and $D$ the effective diffusivity.

The overall performance of the catalyst pellet is traditionally expressed in terms of an {\em effectiveness factor} ($\eta$), which compares the average reaction rate throughout the catalyst to the reaction rate had the conditions (concentration, temperature) in it been uniformly the same as on its surface. For an isothermal reaction, $\eta \downarrow 0$ in the presence of severe diffusion limitations ($\Phi \ll 1$), while $\eta \uparrow 1$ in the absence of transport limitations ($\Phi \approx 1$ and the conditions inside the pellet match those on its surface). For a spherical particle with first-order kinetics,
\begin{equation}
\eta = \left.\frac{3}{\Phi^2}\frac{\delta C}{\delta r} \right\vert_{r=R} .
\end{equation}

In the non-isothermal case where one allows for heat generation during the reaction, however, one may find that $\eta$ exceeds unity. This is so because the hot pellet interior accelerates the reaction, relative to the cool region close to its surface. In this case, both mass and heat transfer play a role. Writing $H$ for the molar heat of reaction, $K$ for the heat conductivity, $T_0$ for the boundary temperature, and 
$C_0$ for the boundary concentration, one can define the parameter $\beta = C_0HR/(KT_0)$ characterizing non-isothermal behavior. In the isothermal case ($\beta=0$), $\eta$ is in an one-to-one correspondence with $\Phi$. In non-isothermal cases ($\beta \neq 0$), the relation between $\Phi$ and $\eta$ loses injectivity and a single $\eta-$value can correspond to various $\Phi-$values. Another significant parameter when considering non-isothermal reactions is $\gamma = Q/(RT_0)$, with $Q$ the activation energy of the reaction. This parameter can be interpreted as the sensitivity of the reaction to temperature changes, since its logarithm corresponds to the Arrhenius expression of the reaction rate outside the pellet. For our discussion in the main text, we fixed the $\gamma-$value to $20$ and varied the value of $\beta$ to produce different model response curves. For detailed information regarding the numerical solution of $\eta$ vs. $\Phi$ curves in non-isothermal catalysis, we refer the interested reader to~\cite{weisz_behaviour_1995}.

With this data in hand, we can proceed to apply DMAPS as usual. To use ``offset'' data ($\eta_{i + \Delta}$, cf.~main text) as part of our model response, we generated the $(\Phi,\eta)-$curve using a regular grid of step $\Delta$ in $\log{\Phi}$. That way, each model response $\eta_i$ corresponding to the grid point $\log{\Phi_i}$ could easily be combined with $\eta_{i + \Delta}$ corresponding to the adjacent grid point $\log{\Phi_i}+\Delta$.

\section{Characterizing the ``good parameter set''}
Traditionally, parameter sensitivities have been analyzed by
inspecting the eigenvalues of the Hessian of some objective function 
near a reported minimum. 
Vanishingly small eigenvalues suggest directions in parameter space
in which the goodness of fit remains nearly invariant. 
Such directions provide us with a sense of the dimensionality of the 
``good parameter set'' (or set of good fits) -- the set of parameter 
values leading to an objective function value practically indistinguishable from its value at the reported minimum.
For the ABC model presented in the main text, this ``good set'' is 
visibly 2D (see Fig.~6 of the manuscript), and one might be 
tempted to exploit this feature to determine the number of 
effective parameters (see also the recent work of \cite{lamont2017correspondence} for a connection
to ideas from statistical mechanics).
Specifically, the ABC model has a total of three parameters 
and two neutral directions parameterizing the ``good set,'' which
suggests the existence of a {\em single} effective parameter. 
This is confirmed by our output-informed kernel DMAPS computations for that model in the 
main text.

This short section uses a somewhat contrived example to showcase a caveat: 
that nonlinearity {\em in the way the inputs enter the model} may obscure the true
dimensionality of the ``good set'' and, through this, lead to an erroneous 
estimation of the degree of model sloppiness. 
We first introduce an ODE model followed by a transformation of the states which in
total contains four parameters
$(\lambda, \varepsilon, a, b)$.
The two-dimensional ODE system 
\begin{align}
  \begin{aligned}
    X' &= -\lambda X , \\
    \varepsilon Y' &= -Y ,
    \label{eqn:sp}
  \end{aligned}
\end{align}
is followed by the transformation $(X, Y) \mapsto (x,y)$ given by
\begin{align}
  \begin{aligned}
    x &= X + b y^2 , \\
    y &= Y + a x^2 .
    \label{eqn:sp-t}
  \end{aligned}
\end{align}
where $\lambda$ and $a$ control slow contraction rate and slow manifold topology
respectively (these two are effective parameters). 
The parameter $\varepsilon$ dictates the fast transient and $b$ controls
fast fiber shape ($\varepsilon$ and $b$ are here sloppy parameters).
When $\varepsilon \ll \lambda$, Eq.~\eqref{eqn:sp} becomes 
singularly perturbed and $Y$ quickly decays to zero. 
The transformation in (\ref{eqn:sp-t}) serves to create \emph{nonlinear} 
fast and slow manifolds in the $(x,y)$ plane, $x = X(t_0) + b y^2$ 
and $y = a x^2$, respectively.
To make our point, we now transform the parameters $a$ and $\lambda$
to two other parameters $u_2,w_2$ that are invertible functions of them.
This is accomplished through two iterations of the H\'enon map (for $A=1.4$ and $B=0.3$)
(that provide the nonlinear invertible transformation) as
\begin{align}
  \begin{aligned}
    u_2 &= 1 - A (1 - A \lambda^2 + a)^2 + B (1 - A \lambda^2 +
    a) , \\
    w_2 &= b (1 - A \lambda^2 + a).
  \end{aligned}
\end{align}
The new model (which can be thought of as an observation of \eqref{eqn:sp} through a ``curved mirror'') 
has two parameters $\mathbf{p}=(u_2, w_2)$, which are in an one--to--one 
correspondence with the original parameters $(\lambda, a)$. 
We now fix a base value of $\mathbf{p}^* = (0.7956, 1.8)$, 
corresponding to $(\lambda^*, a^*) = (1, 1)$ and compute the model response
\begin{align}
  \omr(\mathbf{p}) = \begin{bmatrix} x(t_0 \vert \mathbf{p}) & y(t_0 \vert \mathbf{p}) \\
    x(t_1 \vert \mathbf{p}) & y(t_1 \vert \mathbf{p}) \\ \vdots & \vdots \\ x(t_{N} \vert \mathbf{p}) & y(t_{N} \vert \mathbf{p}) . \end{bmatrix}
\end{align}
\begin{figure*}[!htp]
 \centering
 \begin{tabular}{cc}
  \includegraphics[width=0.49\textwidth,height=0.2\textheight]{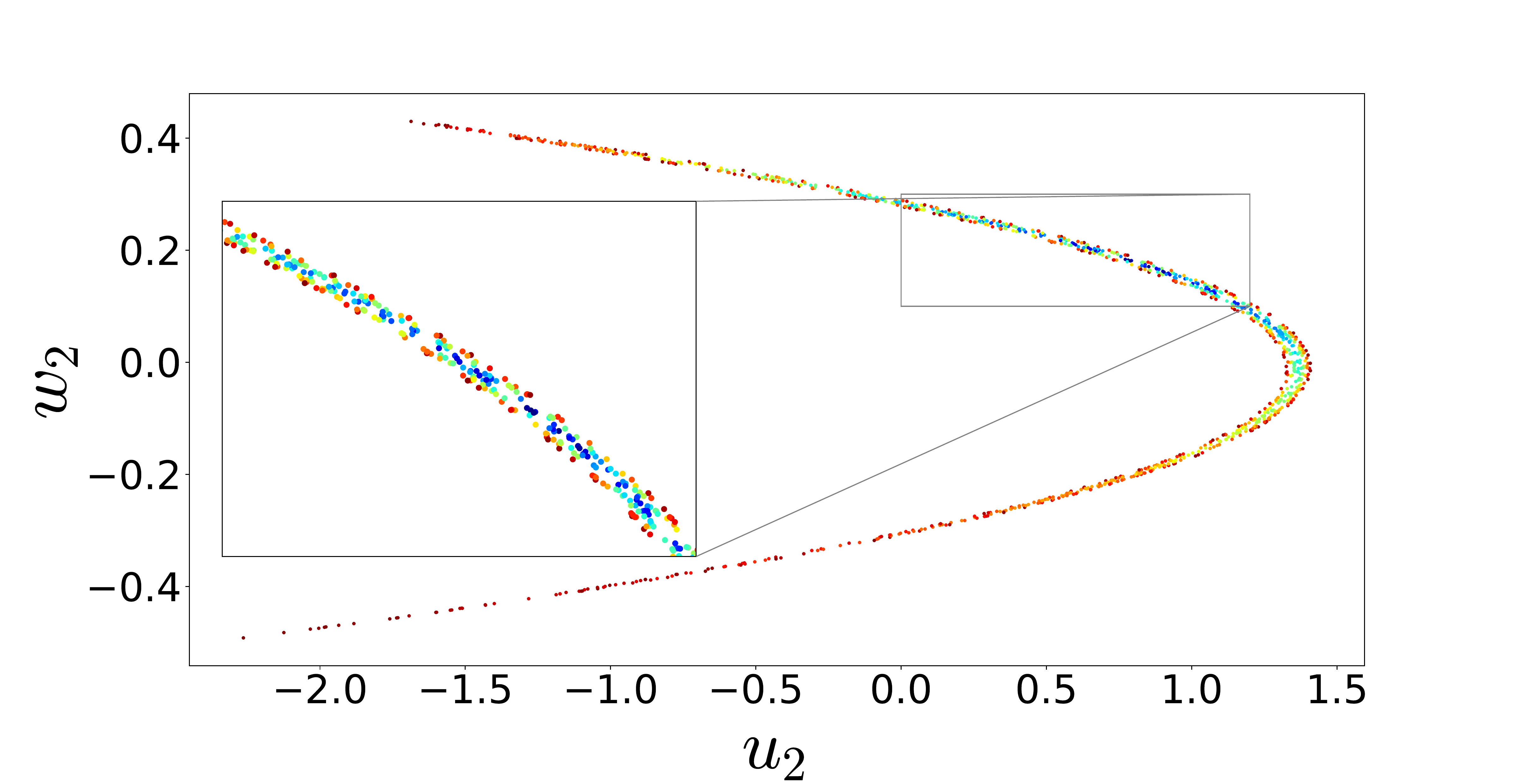}
 &
  \includegraphics[width=0.49\textwidth,height=0.2\textheight]{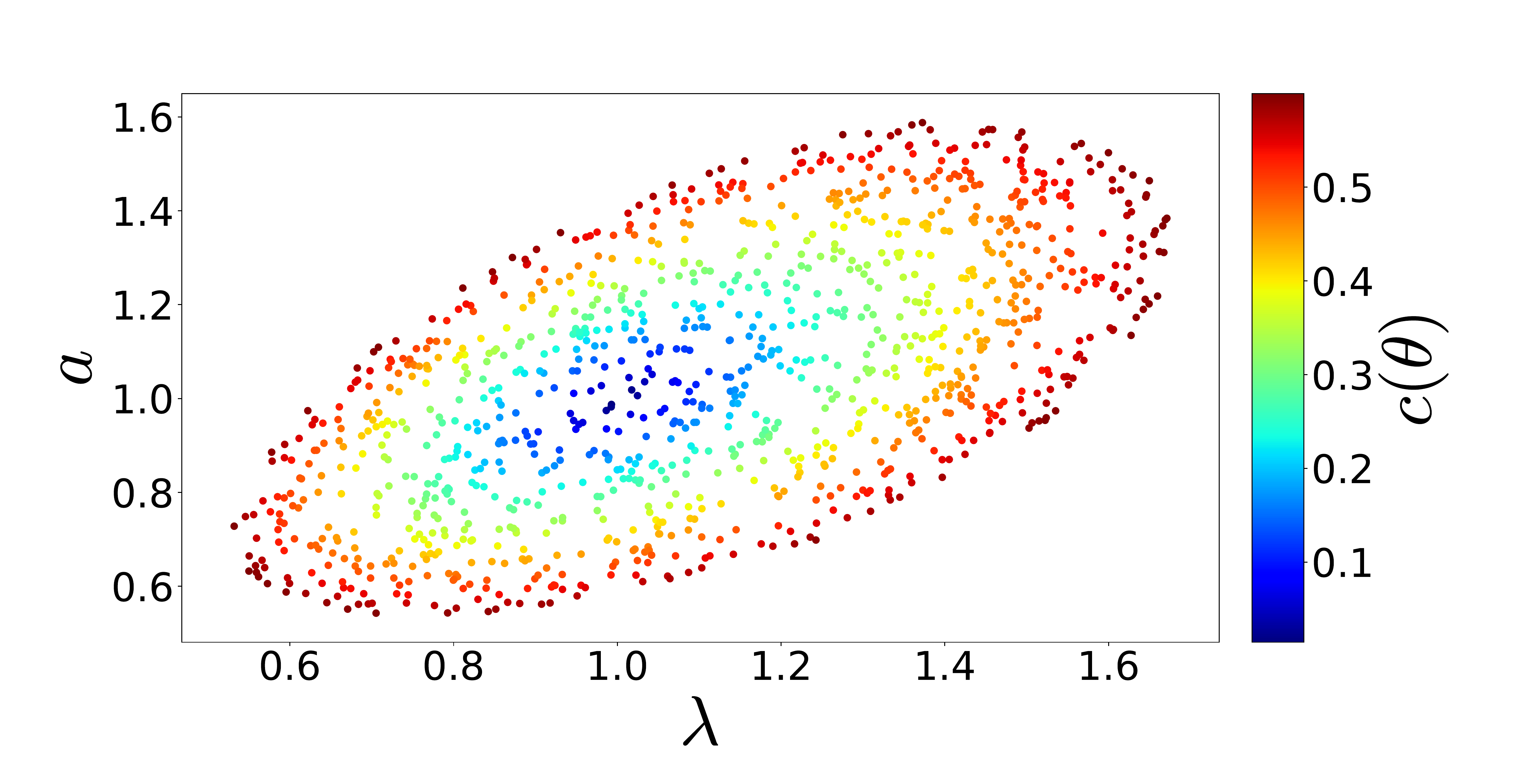}
 \\
  (a) & (b)
 \\
 \end{tabular}
 \caption{(a) Sample of the transformed ($u_2,w_2$) parameter space colored by $c$. Significant deviations from the expected ellipsoidal ``good'' parameter set are observed. (b) Original ($a, \lambda$) parameter space, in which the expected ellipse is recovered. Both figures share the color bar on the right. \label{fig:trans-params}}
\end{figure*}
Here, we fixed $\varepsilon=10^{-3}$, $b=10^{-2}$ and 
$t_1,\ldots,t_N$ to be $N=10$ evenly spaced points in $[0.1,1.0]$. 
To investigate which parameter values $\mathbf{p}$
generate points on the model manifold close to $\omr^* = \omr(\mathbf{p}^*)$, 
we first sampled $\mathbf{p}$ 
uniformly on the rectangle $u_2 \in (-2, 30)$, $w_2 \in (-1.5, 0.7)$; we 
then used each point as an initial value for a least squares minimization 
routine with objective function
\begin{align}
  c(\mathbf{p}) = \| \omr(\mathbf{p}) - \omr^* \|^2_F ;
\end{align}
here, $\| \cdot \|_F$ denotes the Frobenius norm. 
Since noise is not part of our setup, our objective function has a unique minimum at
$\mathbf{p}^*$ making it zero; to reflect that, we terminate our iterative minimization
routine at any point $\theta=\mathbf{p}_i$ satisfying $c(\theta) < 0.8$. 
Since different initializations for our gradient descent algorithm terminate at 
different points, our procedure samples the ``good set'' of parameter settings giving predictions
close to $\omr^*$. 
This set of parameter combinations is shown in Fig.~\ref{fig:trans-params}~(a), 
with each point $\mathbf{p}_i$ colored by its objective function value $c(\theta)$. 
Model nonlinearity is evident in that the set deviates markedly from the expected elliptical shape
close to an isolated minimum.
Transforming this set back to the original parameters $(a, \lambda)$, as in 
Fig.~\ref{fig:trans-params}~(b), 
we recover the typical, elliptical structure expected around the perfect fit.

\begin{figure*}
\centering
\begin{tabular}{cc}
\includegraphics[width=0.49\textwidth,height=0.2\textheight]{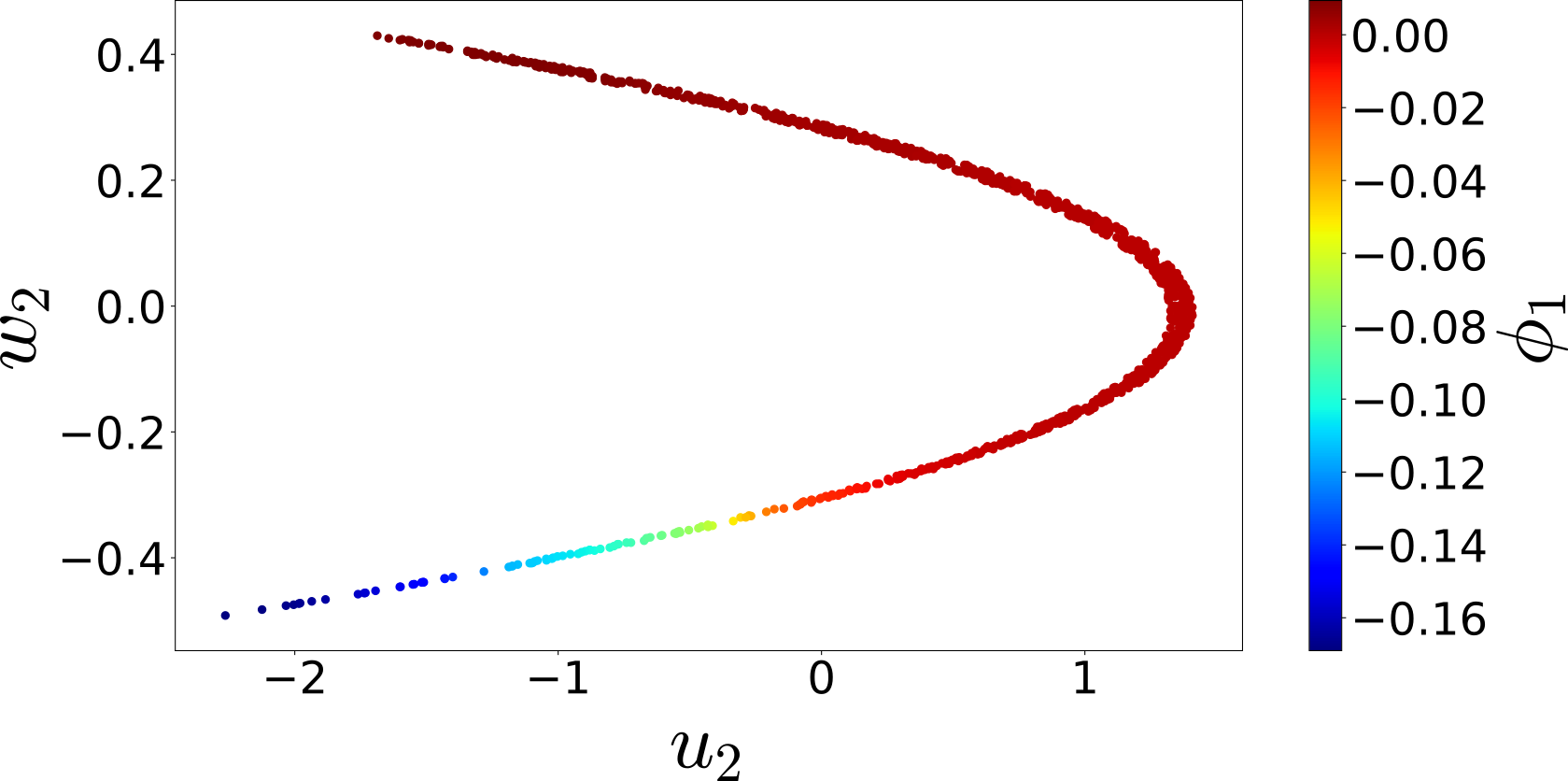} &
\includegraphics[width=0.5\textwidth]{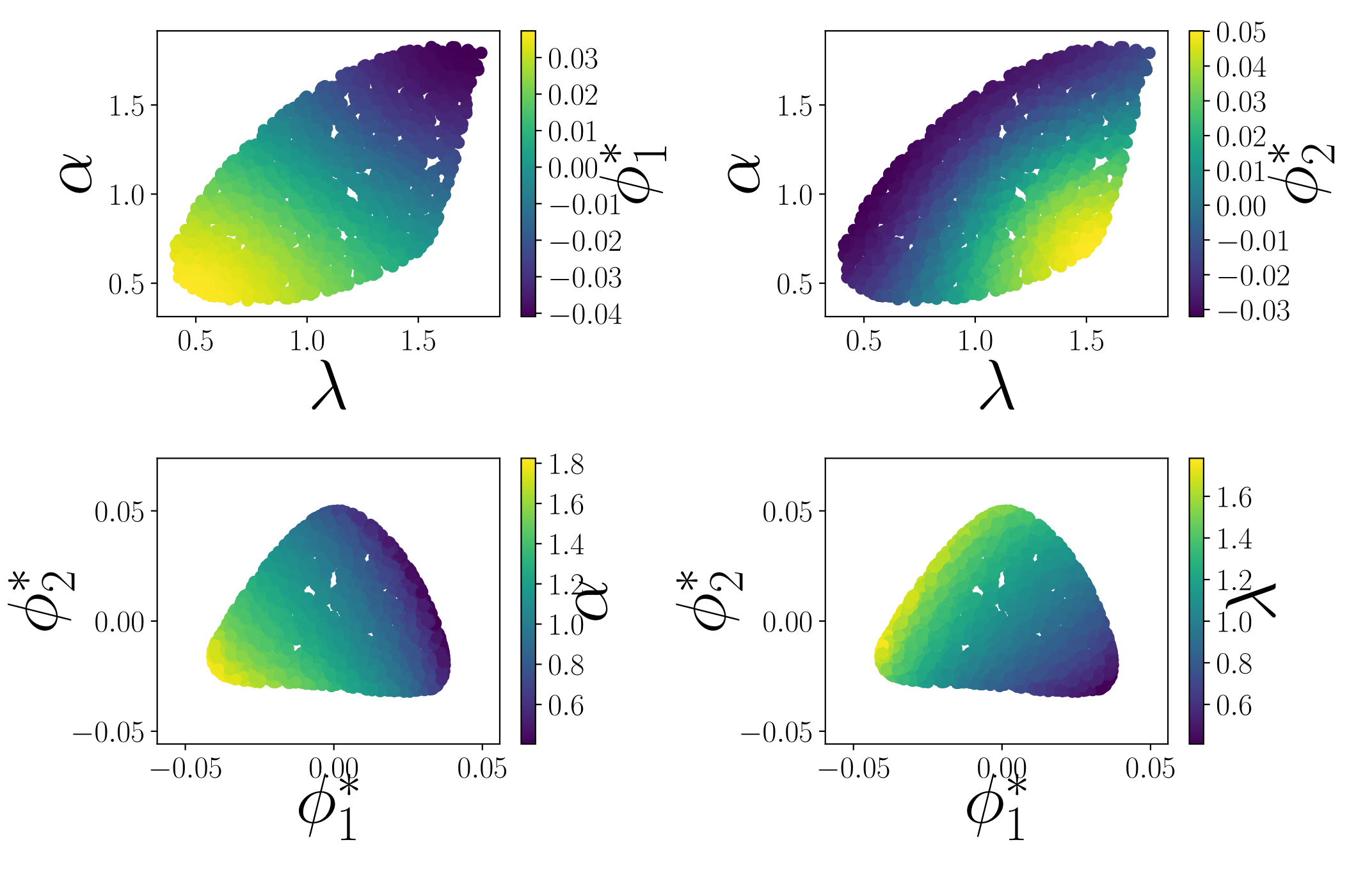} \\
(a) & (b)
\end{tabular}
\caption{
(a) Coloring the $(u_2, w_2)$ plane by $\phi_1$ from a DMAPS analysis based solely on
input--only (parameter) informed kernel. 
The long curve is captured but the thin dimension (the width) is not resolved (the next thirty eigenvectors did not capture this new direction!). 
(b) Original parameters $\alpha$ and $\lambda$
colored by the mixed, input--output kernel DMAPS coordinates 
$(\phi_1^*, \phi_2^*)$ together with diffusion space colored 
by the original parameters. The two-dimensional effective nature of the ``good''
input set (and its correlation with the original inputs $(a, \lambda)$ is clearly visible.}
\label{fig:henon-params-dmaps}
\end{figure*}

A standard DMAPS analysis of that ``good set'' with the \emph{input}-only informed 
kernel, i.e with the Euclidean norm in $\mathbf{p}$ (i.e. in $u_2,w_2$), 
suggests an apparent dimensionality of one. 
The first DMAPS eigenvector parameterizes the long, 
curved, thin direction ``along'' the cloud, see Fig.~\ref{fig:henon-params-dmaps}(a),
while the second, thin dimension is ``lost'' in subsequent, higher-order eigenvectors.

As discussed in Eq.~8 of the main text, applying a more informative kernel,
that includes both input and outputs, in our data-driven DMAPS analysis 
can give a more informative result.
Figure~\ref{fig:henon-params-dmaps}~(b) shows the {\em original} parameter space $(a, \lambda)$
parametrized by the first two non-trivial DMAPS eigenvectors $(\phi_1^*,\phi_2^*)$ (top);
and the converse (bottom) using the mixed kernel; as we expected the ``good'' parameter set 
now appears visibly two dimensional.

This phenomenon is the result of our (intentionally) poor choice of the new model parameters ($u_2,w_2$).
Data-driven approaches, such as DMAPS, can thus help us 
reparameterize (i.e. appropriately transform)
parameter space $(u_2, w_2)$ to a new one $(\phi_1^*, \phi_2^*)$ 
that has a much better bi-Lipschitz 
relation with the original parameter set $(a, \lambda)$ 
better resolving model variability.

\section{A quick discussion of Active Subspaces}

\begin{figure}[!h]
\centering
\begin{tabular}{cc}
\includegraphics[width=0.4\textwidth]{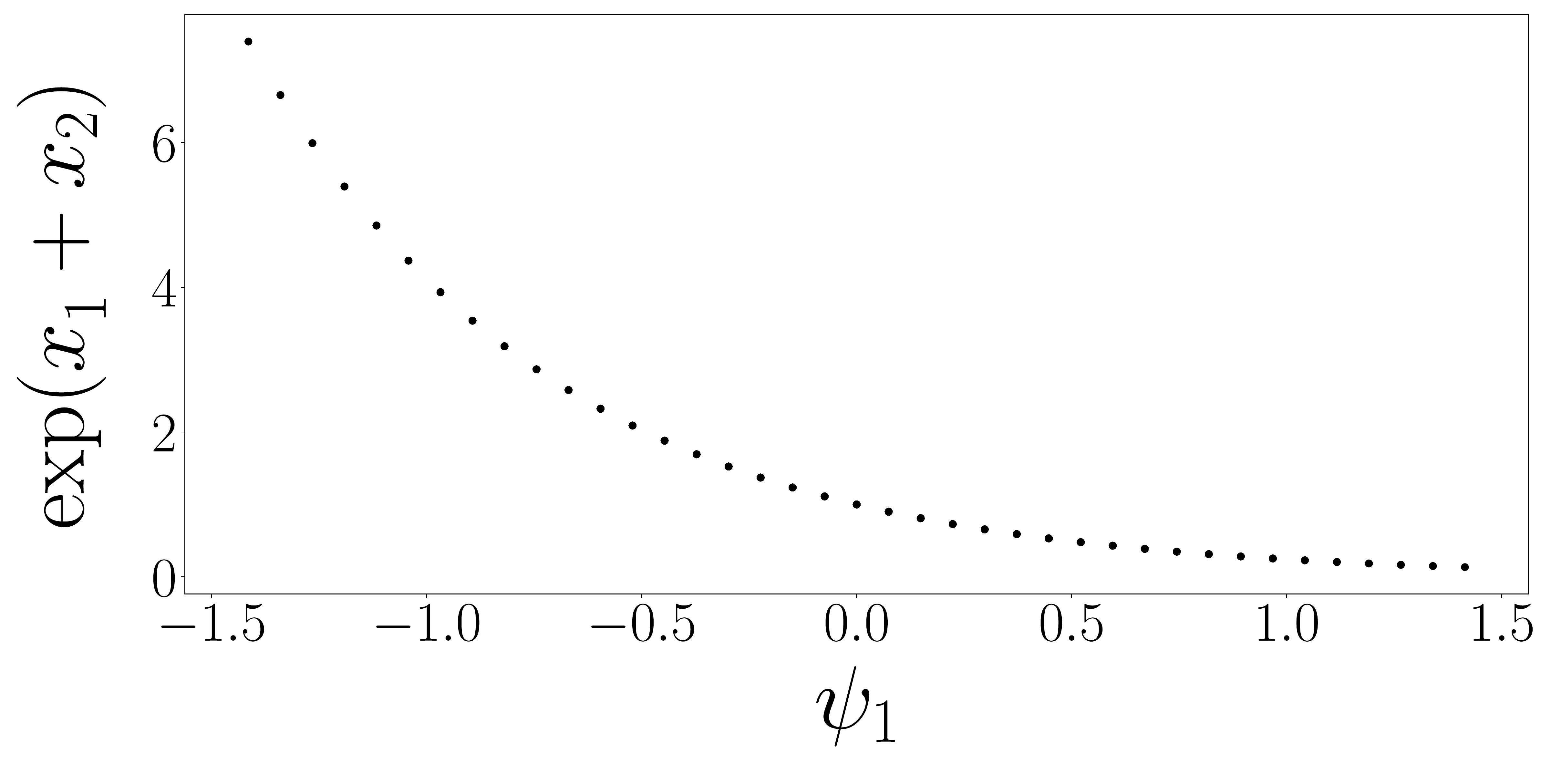} &
\includegraphics[width=0.4\textwidth]{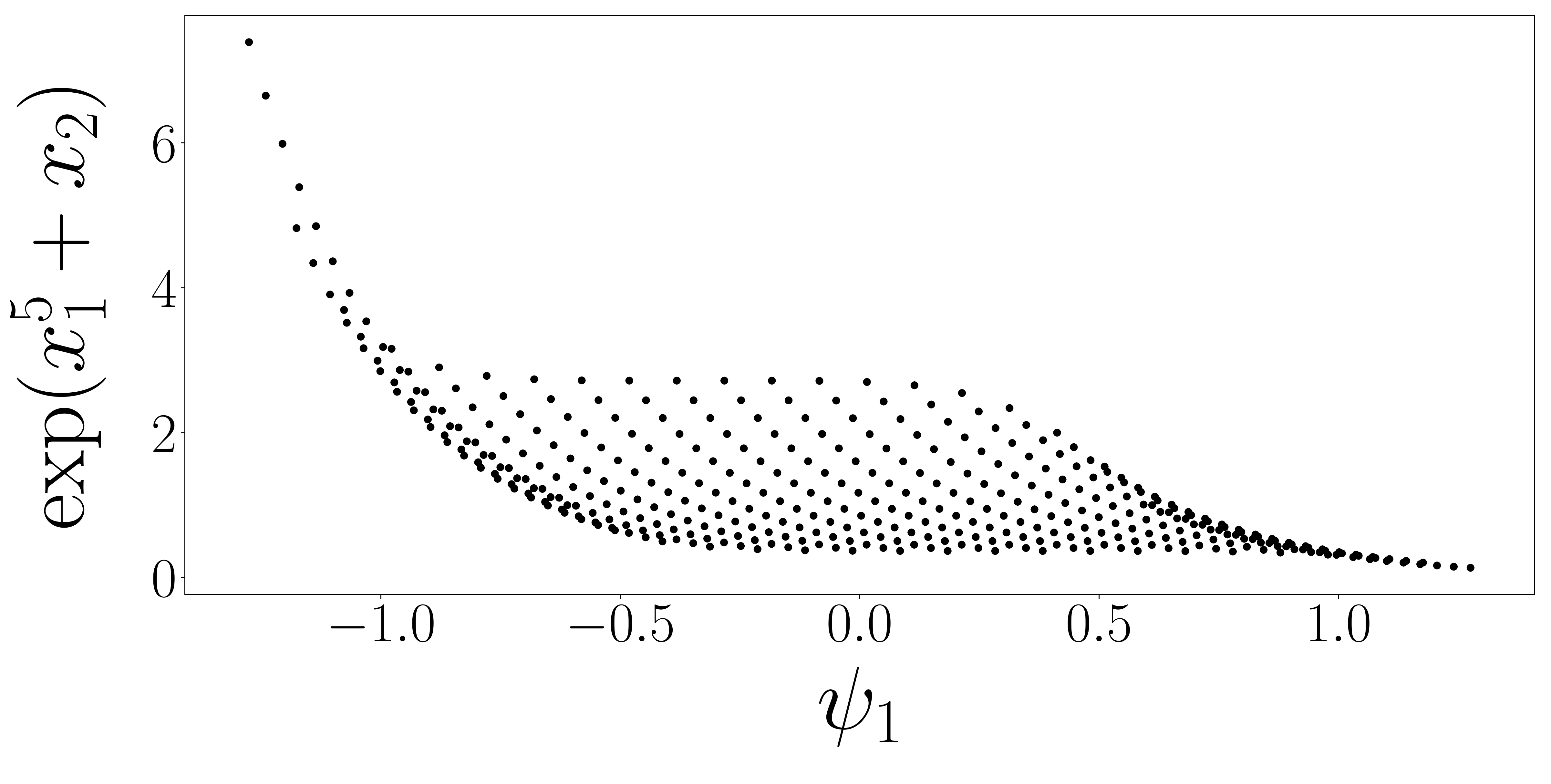}\\
(a) & (b)
\end{tabular}
\caption{The observable $f$, plotted against the first active subspace coordinate for (a) $\alpha=1$ and (b) $\alpha=5$.}
\label{fig:active_1}
\end{figure}

The Active Subspaces algorithm of P. Constantine and coworkers \cite{constantine_active_2014} has been developed based on the idea of finding the ``important directions'' in the space of all inputs of a nonlinear scalar function $f:\mathbb{R}^m\rightarrow\mathbb{R}$, a map from $m$-dimensional parameter space to the real line. The aforementioned directions are assumed to be weighted \textit{linear} combinations of the input parameters. These directions are called ``active subspaces'' and point towards the direction of most intense change of an ``observable.''
\begin{enumerate}
	\item We consider $N$ sample points in parameter space. For each sample point $\mathbf{x_n} \in \mathbb{R}^m$, we observe $f(\mathbf{x_n})$ and the gradient $\nabla f(\mathbf{x_n})$.
	\item The average of the outer product of the gradient with itself on the sample data is computed through
	\begin{equation*}
		\hat{C}=\frac{1}{M}\Sigma_{n=1}^N
        \nabla f(\mathbf{x_n})
        \nabla f(\mathbf{x_n})^T.
	\end{equation*}
	      Here, $\nabla f$ is seen as an $m-$dimensional column vector, hence $\hat{C}$ is an $m \times m$ matrix.
	\item We find the eigendecomposition of $\hat{C} = \hat{W} \hat{\Lambda} \hat{W}^T$.
	\item 
    The observable exhibits the greatest change in the space spanned by the leading eigenvectors of $\hat{C}$, which are stored column-wise in $\hat{W}$ \cite{constantine_active_2014}.
\end{enumerate}
To illustrate this procedure and its outcome on a simple example, we consider the map
\begin{equation}
f : \mathbb{R}^2 \rightarrow \mathbb{R}
\quad \mbox{given by} \
 (x_1,x_2)\mapsto f(x_1,x_2) = \mathrm{e}^{x_1^\alpha + x_2} .
\label{active-example}
\end{equation}
where $\alpha$ can be fixed at an arbitrary value. Although the parameter space $(x_1,x_2)$ is 2-D, the map effectively only depends on the single effective parameter $x_1^\alpha + x_2$. As a result, the parameter space can be re-parameterized by any one-to-one function of $x_1^\alpha + x_2$. To discover active subspaces for \eqref{active-example}, we consider a uniform grid on $[-1,1] \times [-1,1]$ and evaluate the map $f$ on each grid point. Since the parameter space is effectively 1-D (composed of level curves of the effective parameter), the active subspace parameterization is given by $\psi_1 = w_1^T\cdot [x_1,x_2]^T$, at each point. We compute the active subspaces for $\alpha=1$ and $\alpha=5$:
\begin{itemize}
\item For $\alpha=1$: 
As expected, the effective 
parameter $\psi_1$ is the linear combination $x_1+x_2$ of the input parameters. As shown in 
Fig.~\ref{fig:active_1}(a), the observable is in an one-to-one correspondence with $\psi_1$, meaning that active subspaces recover the effective parameter.

\item For $\alpha=5$: 
In this case, the effective parameter is not a linear combination of input parameters. As shown in Fig.~\ref{fig:active_1}(b), the basic algorithm does not discover the effective one-dimensional (but nonlinear) relation between the input parameters.

\end{itemize}

\begin{figure}[H]
\centering
\begin{tabular}{cc}
\includegraphics[width=0.42\textwidth]{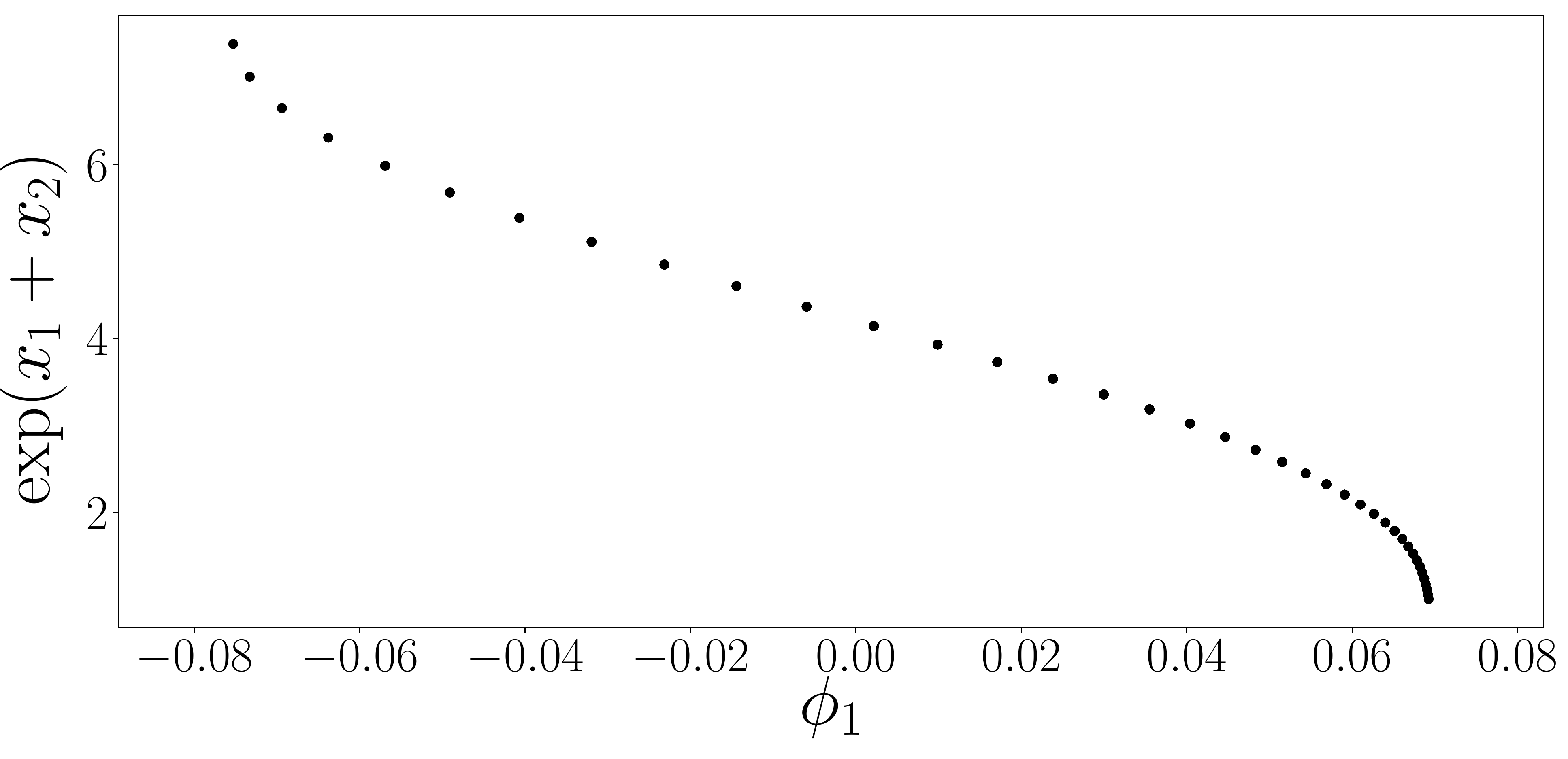} &
\includegraphics[width=0.4\textwidth]{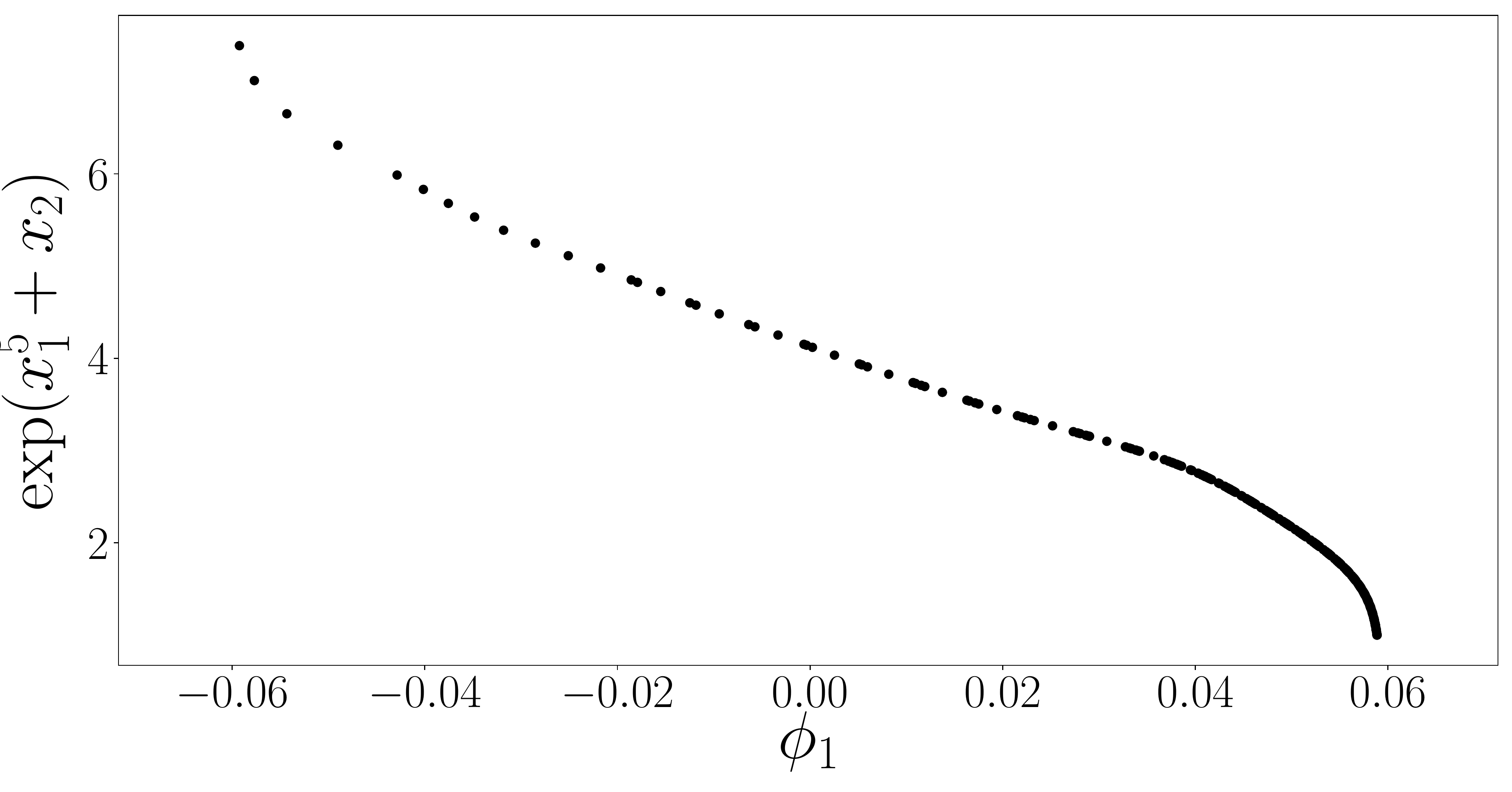} \\
(a) & (b)
\end{tabular}
\caption{The observable $f$ plotted against
the first DMAPS coordinate for (a) $\alpha=1$ and (b) $\alpha=5$.}
\label{fig:dmap_toy}
\end{figure}

\noindent
To conclude this discussion and compare methodologies, we applied the output-only DMAPS algorithm to each of these cases and plotted the observable against the first non-trivial eigenfunction of the graph Laplacian; see Fig.~\ref{fig:dmap_toy}. Plainly, th DMAPS algorithm captures the effective parameter both when that is linear and when it is nonlinear in the input parameters.
The idea of active subspaces has been employed by Constantine and coworkers to decompose input space using products of powers of the inputs, since these become linear combinations in a logarithmic scale (``ridge functions'' \cite{constantine2017near}). It will be interesting to explore how more general nonlinear relations can be used in discovering effective nondimensionalizations.

\section{On the origins of sloppiness \label{s-origins}}
We now attempt to generalize the setting presented in Section~\ref{setup} and explore the origins of parameter (input) sloppiness. Here also, we consider multivariable vector functions $\mathbf{x}(\mathbf{t} \vert \p) \in \R^D$, where $D$ is arbitrary. The independent variables are partitioned in what one might call ``bona fide variables" $\mathbf{t} = (t_1,\ldots,t_K) \in I \subset \R^K$ and \emph{parameters} $\p = (p_1,\ldots,p_M) \in \ps \subset \R^M$. 
For each $\p\in\ps$, we term the function $\fmr(\p) = \mathbf{x}(\,\cdot\, \vert \p) \in \fms$ the \emph{full model response} and the function space $\fms$ containing it the \emph{full model space}. To know the mapping $\fmr : \ps \to \fms$ is to know fully the solution to the problem at hand, and such knowledge is typically unattainable. As $\p$ ranges over $\ps$, $\fmr$ traces out the \emph{full model manifold} $\fmm = \{ (\p,\fmr(\p)) \}_{\p \in \ps} \subset \ps\times\fms$. That manifold is generically $M-$dimensional and, if $\mathbf{x}$ depends on $\p$ in a $C^1$ manner, also continuously differentiable. Here also, as in Section~\ref{setup}, we only monitor one or more \emph{functionals}, $f_1 , \ldots , f_N : \fms \to \R$. We term each individual $f_n(\fmr(\p))$ a \emph{partial observation} and the $N-$tuple
\[
 \omr(\p)
=
\left[ f_1(\fmr(\p)) , \ldots , f_N(\fmr(\p)) \right]
\in
 \R^N
\]
the \emph{observed} model response. Under the action of $\omr$, the full model manifold $\fmm$ is projected to the \emph{(observed) model manifold} $\omm = \{ (\p,\omr(\p)) \}_{\p \in \ps} \subset \ps \times \R^N$. We will assume $\omr : \ps \to \omr(\ps)$ to be a homeomorphism, unless explicitly stated, and the linear map $\D_\p \omr$ to have rank $M$. These assumptions make $\omr$ an atlas for $\omm$ and are satisfied by all models we consider in the main text, save for the non-invertible model therein.\\

To discuss sensitivity of the observed model response to parameter variations, one must assess shifts in that response relative to such variations. In \cite{TMS11}, the authors defined distances on $\omm$ by omitting the $\p-$component and working with the projected manifold $\pi\omm$, where $(\p,\omr(\p)) \xmapsto{\pi} \omr(\p)$. Distances were measured using a Riemannian metric induced by a norm in the carrier space of $\pi\omm$, i.e. $\R^N$. Similarly, in this section we work with the standard Euclidean norm, reserving important questions on weighing and correlation of partial observations for the future. An infinitesimal displacement $\df\p = (\df p_1 , \ldots , \df p_M)^{\rm T}$ in $\ps$
yields the infinitesimal displacement $\df z = (\df z_1 , \ldots , \df z_N)^{\rm T} = (\D_\p\omr)\,\df\p \in \mathrm{T}_\p(\pi\omm)$ with length $\norm{\df z}^2 = (\df\p)^{\rm T} \, g \, \df\p$.
The $M \times M$ positive definite matrix $g = (\D_\p\omr)^{\rm T}(\D_\p\omr)$ is the \emph{metric tensor} for $\pi\omm$ for the specific atlas $\omr$. Although the positive definiteness of $g$ makes the system responsive to \emph{all} parameter variations, the observed model response locally around some point $\p \in \ps$ may vary greatly with the direction of $\df\p$ and be disproportionately small along certain directions. This is the phenomenon termed \emph{sloppiness}, and it manifests itself in the spectrum of the metric. Small eigenvalues yield small observed responses, with sloppy directions in $\ps$ being the pull-backs under $\mu$ of the associated eigendirections in $\mathrm{T}_\p(\pi\omm)$.\\

In terms of the full model manifold, parameter values are first mapped to $\fmm \subset \ps\times\fms$ equipped with the Riemannian metric $h = (\D_\p\fmr)^{\rm T}(\D_\p\fmr)$, then projected to the model manifold $\omm \subset \R^N$, equipped with the metric $g = (\D_\p\omr)^{\rm T}(\D_\p\omr)$. We can quantify the relation between these two metrics and examine how that projection can generate non-identifiability as well as sloppiness. For clarity of presentation, we restrict our attention to linear functionals $f_1,\ldots,f_N$ on a Hilbert space $\fms$ equipped with inner product $\langle\cdot,\cdot\rangle$ and induced norm $\norm{\cdot}$. In that setting, $\fms$ is isomorphic to its dual $\fms^*$ and hence $f_n = \langle e_n,\cdot \rangle$, for some $e_n \in \fms$ and all $n = 1,\ldots,N$.
Therefore, $\omr(\p) = [ \langle e_1,\fmr(\p) \rangle , \ldots , \langle e_N,\fmr(\p) \rangle ]^\mathrm{T}$,
with push-forward
\[
 \D_\p\omr
=
\left[\begin{array}{ccc}
 \langle e_1,\partial_{p_1}\fmr\rangle
&
 \ldots
&
 \langle e_1,\partial_{p_N}\fmr\rangle
\\
\vdots
&
 \ddots
&
\vdots
\\
 \langle e_N,\partial_{p_1}\fmr\rangle
&
 \ldots
&
 \langle e_N,\partial_{p_N}\fmr\rangle
\end{array}\right]
=
\left[\begin{array}{c}
 \langle e_1,\cdot\rangle
\\
\vdots
\\
 \langle e_N,\cdot \rangle
\end{array}\right]
 \D_\p\fmr .
\]
Based on this, the metric on $\omm$ is written as
\be
 g
=
 (\D_\p\omr)^{\rm T}
 (\D_\p\omr)
=
 (\D_\p\fmr)^{\rm T}
 F
 (\D_\p\fmr) .
\label{g-vs-h}
\ee
The linear operator $F : \fms \to \fms^*$, here, is given for each $v \in \R^N$ by
\[
 Fv
=
\sum_{n=1}^N
 \langle e_n , v\rangle
 \langle e_n , \cdot\rangle
=
\vspace*{4mm}
\begin{array}{c}
\Big[
 \langle e_1,\cdot\rangle
,
 \ldots
,
 \langle e_N,\cdot\rangle
\Big]
\\
\vspace*{4mm}
\end{array}
\hspace*{-3mm}
\left[\hspace*{-2mm}\begin{array}{c}
 \langle e_1,v\rangle
\\
 \vdots
\\
 \langle e_N,v\rangle
\end{array}\hspace*{-2mm}\right]
\in
 \fms^* .
\]
To understand \eqref{g-vs-h} better,
we use the isomorphism $\fms^* \cong \fms$ to interpret $F$ as a symmetric endomorphism on $\fms$.
Its spectrum consists of the zero eigenvalue,
linked to the co-dimension $N$ kernel $\mathrm{Ker}F = \bigcap_{n=1}^N e_n^\perp$,
and of a nontrivial part linked to the invariant subspace
$\mathrm{Im}F = \mathrm{span}(e_1,\ldots,e_N)$.
In the basis $\{e_1,\ldots,e_N\}$ for $\mathrm{Im}F$,
the restriction $F\vert_{\mathrm{Im}F}$ is represented by the matrix
\be
 G
=
\left[\begin{array}{ccc}
 \langle e_1,e_1\rangle
&
 \ldots
&
 \langle e_1,e_N\rangle
\\
 \vdots
&
 \vdots
&
 \vdots
\\
 \langle e_N,e_1\rangle
&
 \ldots
&
 \langle e_N,e_N\rangle
\end{array}\right] .
\label{Gram-mat}
\ee
This proves that the nontrivial part of the spectrum consists of the eigenvalues of the $N \times N$ Gram matrix $G$.
To rewrite \eqref{g-vs-h} using this information,
we decompose the columns of $\D_\p\fmr$ along the invariant subspaces $\mathrm{Ker}F$ and $\mathrm{Im}F$,
\[
 \partial_{p_n}\fmr
=
 \mathrm{N}_n
+
 \mathrm{I}_n
=
 \mathrm{N}_n
+
 \left[ e_1 , \ldots , e_N \right]
 C_n ,
\quad \mbox{with} \
 \mathrm{N}_n \in \mathrm{Ker}F
\ \mbox{and} \
 \mathrm{I}_n \in \mathrm{Im}F .
\]
The corresponding matrix decomposition is $\D_\p\fmr = \mathrm{N} + \left[ e_1 , \ldots , e_N \right] C$,
with $C_{mn}$ the component of $\partial_{p_n}\fmr$ along $e_m$,
and thus
\be
 g
=
 C^\mathrm{T}
 G^2
 C .
\label{g-ito-proj}
\ee
This equation expresses the metric on $\omm$ in terms of a Gram matrix, determined by the functionals $f_1,\ldots,f_N$,
and of a matrix quantifying the projection of $\mathrm{T}_{\fmr(\p)}\fmm$ on the combined range of those functionals along their joint kernel. \\

Equation~\eqref{g-ito-proj} shows that an ill-conditioned $G$ or $C$ leads to sloppiness, manifested in disparities within $\sigma(g)$. An ill-conditioned matrix $G$ points to functionals that are either badly scaled or nearly dependent. This was the case in \cite{WCGBMBES06}, where sloppiness was traced to a Vandermonde matrix specific to Taylor polynomials. Replacing those polynomials by an orthonormal set would have sufficed to remove sloppiness. An ill-conditioned matrix $C$, instead, is due to directions in $\mathrm{T}_\p\fmm$ that align well with $\bigcap_{n=1}^N \mathrm{Ker}f_n$. In that case, parameter variations can generate negligible model responses on $\pi\omm$, although the observation functionals are proportionate and the \emph{full} model response on $\fmm$ appreciable. This is the case with multiscale systems, in which certain parameters combinations affect behavior at unobserved scales (fast transients).

\providecommand{\bysame}{\leavevmode\hbox to3em{\hrulefill}\thinspace}
\providecommand{\MR}{\relax\ifhmode\unskip\space\fi MR }
\providecommand{\MRhref}[2]{%
  \href{http://www.ams.org/mathscinet-getitem?mr=#1}{#2}
}
\providecommand{\href}[2]{#2}
\bibliographystyle{elsarticle-num-names}